\documentclass[12pt]{amsart}
\usepackage{amsmath,amssymb,amsfonts,latexsym}
\usepackage{mathrsfs}%pour les caracteres en calligraphie
\usepackage[dvips]{geometry}
\usepackage{array}
\usepackage{epsf}
\usepackage{epsfig}
\newtheorem*{lemma}{Lemma}
\newtheorem*{prop}{Proposition}
\newtheorem*{thm}{Theorem}
\newtheorem*{cor}{Corollary}

\newcommand{\twoheaddownarrow}{\overset{\sim}{\twoheaddownarrow}}

\newcommand{\nc}{\newcommand}
\nc{\Ker}{\operatorname{Ker}} \nc{\rker}{\operatorname{rKer}}
\nc{\im}{\operatorname{Im}}
\nc{\stab}{\operatorname {Stab}}
\nc{\ann}{\operatorname {Ann}}
\nc{\Id}{\operatorname {Id}}
\nc{\Prim}{\operatorname {Prim}}
\nc{\Real}{\operatorname {Re}}
\nc{\Ext}{\operatorname {Ext}}
\nc{\rad}{\operatorname {rad}}
\nc{\rk}{\operatorname {rank}}
\nc{\Aut}{\operatorname {Aut}}
\nc{\supp}{\operatorname {supp}}
\nc{\height}{\operatorname {ht}}

\usepackage{array}
\usepackage{epsf}
\usepackage{epsfig}
\usepackage{xcolor}

\usepackage{tikz,tikz-cd}

\usetikzlibrary{matrix}

\usetikzlibrary{fit}

\usetikzlibrary{matrix,backgrounds}
\pgfdeclarelayer{myback}
\pgfsetlayers{myback,background,main}

\tikzset{mycolor/.style = {line width=1bp,color=#1}}%
\tikzset{myfillcolor/.style = {draw,fill=#1}}%
\usetikzlibrary{arrows,matrix,positioning}
\usepackage{pstricks}

\usetikzlibrary{graphs,graphs.standard}
\usepackage{changepage}
\usepackage{tikz-cd}
\usepackage{mathrsfs}
\usetikzlibrary{positioning}
\makeatletter

\usepackage{circledsteps}

\newcommand*{\encircled}[1]{\relax\ifmmode\mathpalette\@encircled@math{#1}\else\@encircled{#1}\fi}
\newcommand*{\@encircled@math}[2]{\@encircled{$\m@th#1#2$}}

\newcommand*{\@encircled}[1 ]{%
  \tikz[baseline,anchor=base]{\node[draw,circle,outer sep=0.1pt,inner sep=.18ex,  line width = 1 pt ] {#1};}}

\newcommand{\cir}[1]{\tikz[baseline]{%
    \node[anchor=base, draw, circle, inner sep=0, minimum width=1.2em]{#1};}}
\usepackage{float}
\usepackage{xcolor}
\newcommand{\violet}[1]{\textcolor{violet}{#1}}
\usepackage{lscape}

\begin{document}

\title [Component Tableaux]{The Magic and Mystery of Component Tableaux.}
\author [Yasmine Fittouhi and Anthony Joseph]{Yasmine Fittouhi and Anthony Joseph}
\date{\today}
\maketitle
\vspace{-.9cm}\begin{center}

Department of Mathematics\\
The Weizmann Institute of Science\\
Rehovot, 7610001, Israel\\
fittouhiyasmine@gmail.com
\end{center}\

\

\

\vspace{-.9cm}\begin{center}
Donald Frey Professional Chair\\
Department of Mathematics\\
The Weizmann Institute of Science\\
Rehovot, 7610001, Israel\\
anthony.joseph@weizmann.ac.il
\end{center}\

\

\date{\today}
\maketitle



 \

\textbf{Abstract}.

Let $G$ be a simple algebraic group over the complex field $\mathbb C$, $B$ a fixed Borel subgroup, $P$ a parabolic subgroup containing $B$, $P'$ its derived group and $\mathfrak m$ the Lie algebra of its nilradical.

The nilfibre $\mathscr N$ for this action is the zero locus  of the augmentation $\mathscr I_+$ of the semi-invariant algebra $\mathscr I:=\mathbb C[\mathfrak m]^{P'}$.  %Its components are studied here in type $A$.

 Although the Lasker-Noether theorem allows one to define the components of an algebraic variety, their determination is notoriously difficult, especially for $\mathscr N$ as it is  homogeneous.  Even for $G=SL(n)$ practically nothing was known previously.  The only result of comparable, but lesser complexity, is for $\mathscr V:=\mathscr O\cap \mathfrak n$, with $\mathscr O$ a nilptent $G$ orbit and $\mathfrak n$  the set of strictly upper triangular matrices.  Then $\mathscr V$ is equi-dimensional, with components, known as orbital varieties, parameterised by standard tableaux whose shape is dictated by $\mathscr O$. This result is made relatively easy through Bruhat decomposition.
%
%  For $G=SL(n)$ an \textit{entirely new} construction is introduced to compute the components of $\mathscr N$.  Thee increase exponentially in $n$ and a priori do not obey any discernable pattern.

  Here the components of $\mathscr N$ are studied for $G=SL(n)$.  They increase exponentially in $n$ with no a priori discernable pattern.  Nor is it possible to use Bruhat decomposition.

     The composition of $n$ defined by the Levi block sizes in  $P$ defines a standard tableau $\mathscr T$. For each choice of numerical data $\mathcal C$, a semi-standard tableau $\mathscr T^\mathcal C$, is constructed from $\mathscr T$. A \textit{delicate and tightly interlocking} analysis constructs a set of excluded root vectors from $\mathfrak m$ such that the complementary space $\mathfrak u^\mathcal C$ has the following properties. First it is a subalgebra of $\mathfrak m$.  Secondly $\mathscr C:=\overline{B.\mathfrak u^\mathcal C}$ lies in $\mathscr N$ to which, thirdly, a Weierstrass section can be associated.  Fourthly $\dim \mathscr C = \dim \mathfrak m-\textbf{g}$, where \textbf{g} is the number of generators of the polynomial algebra $\mathscr I$. Fifthly the  Weierstrass section, is shown to imply that $\mathscr C$  an irreducible component of $\mathscr N$, yet  $\mathscr C$ is  \textit{ only sometimes} an orbital variety closure.

     The resulting ``component map'' $\mathscr T^\mathcal C \mapsto \mathscr C$ is shown to be injective. Evidence for its surjectivity is given.

     \

     \textbf{Acknowledgements}.

     \

     The work of the first author was supported by ISF grant 1957/21 jointly held by M. Gorelik, Weizmann Institute and S. Reif, Bar-Ilan University and by ISF grant 1781/23 held by D. Gourevich, Weizmann Institute.

     \

     The main results described here were the subject of two invited talks, one by Y. Fittouhi and one by A. Joseph in the program ``Algebraic and Combinatorial Methods in Representation Theory'',  held during November 13-24, 2023 in the Tata Institute, Bangalore, India.  We would like to thank the organisers for the invitation and for their hospitality.

     Following this meeting we spent the best of $8$ months writing up the proofs and so it would appropriate to add that this
     research was also significantly aided by the International Centre for Theoretical Sciences (ICTS) through participation in the program -   Algebraic and Combinatorial Methods in Representation Theory (code:  ICTS/ACMRT2023/11).

\section {Introduction} \label {1}  The base field is assumed to be algebraically closed and of characteristic zero.  It will simply be denoted as $\mathbb C$.   Given $m\leq n$ is  positive integers, set $[m,n]:=\{m, m+1,\ldots,n\}$.  If $m$, $n$, or both are omitted we note this respectively as $]m,n],[m,n[,]m,n[$.
 \subsection {The Nilfibre}\label {1.1}

 Let $Q$ be a connected algebraic group acting linearly on a vector space $\mathfrak m$.  Let $\mathscr I$ be the algebra $\mathbb C[\mathfrak m]^Q$ of $Q$ invariant functions on $\mathfrak m$.  The nilfibre $\mathscr N$ (for this action) is the zero variety of the augmentation $\mathscr I_+$ of $\mathscr I$.

 In general the nilfibre is not irreducible and it is hardly ever known how to describe its components (meaning here and henceforth - irreducible components).

  This is the problem we propose to solve for a particular family $\mathscr N(n):n \in \mathbb N^+$, of exponentially growing cardinality.

  \subsection {Associated Weierstrass Sections}\label {1.2}

  Let $\mathscr C$ be a component of $\mathscr N$. A Weierstrass section $e+V$  is a pair $e\in \mathscr C$ and a subspace $V \subset \mathfrak m$ such the restriction of $\mathscr I$ to $e+V$ induces an isomorphism of $\mathscr I$ onto $\mathbb C[e+V]$.

  %In this we shall denote $\mathscr C$ as $\mathscr N^e$.

 The existence of a Weierstrass section implies that $\mathscr I$ is polynomial.  Moreover the generators of $\mathscr I$ are then simultaneously linearized to become a basis for $V$, which is quite an accomplishment.

 Geometrically an orbit of $Q$ can meet a Weierstrass section $e+V$ at only one point, because invariants are constant on orbits and the map is surjective.

 We may wish to demand more of a Weierstrass section.  For example

 \

 $(i)$. That the pair $e,V$ determines $\mathscr C$.

 \

 $(ii)$.  That $Qe$ is dense in $\mathscr C$, so of maximal dimension $\dim \mathscr C$, thereby called a regular orbit.

 \

 $(iii)$.  That every element of $e+V$ generates a regular $Q$ orbit in $\mathfrak m$.

 \

 $(iv)$.   That every regular $Q$ orbit in $\mathfrak m$ meets $e+V$.

 \

 Generally $(iv)$ will fail unless $\mathscr N$ is irreducible.

  \subsection {Coadjoint Action}\label {1.3}

   Suppose $Q$ is a connected, simply-connected, simple algebraic group $G$ acting on its Lie algebra $\mathfrak g$.  Through the classical work of Kostant \cite {K1}, a Weierstrass section exists and $(i)-(iv)$ all hold.  This had some importance for the Toda lattice and its higher dimensional versions. The key point in Kostant's work was to use the principal s-triple, though it is enough to extract from it an ``adapted pair'', \cite {J1}. This concept has been extensively investigated by the second author in papers too numerous to mention.

  The results of Kostant were extended for example for $G=SL(3)$ acting on its $10$ dimensional representation by Popov and Vinberg \cite {PV}. Here one can view the regular orbits as being the elliptic curves in $\mathbb P^2$ and then $(iv)$ above, recovers Weierstrass canonical form. For this reason Popov and Vinberg assigned the epithet - Weierstrass sections.  Because there are two invariant generators (written down by Weierstrass himself) it is not quite obvious that  $\mathscr N$ is irreducible, but this is proved in \cite [Lemma 2.5] {J3}.

  In \cite {J1} adapted pairs were constructed for the coadjoint action of a parabolic subgroup\footnote{Even for a biparabolic subgroup.}  on the nilradical of its Lie algebra in type $A$.  However in this case $\mathscr N$ may fail to be irreducible and some components of $\mathscr N$ may fail to admit regular elements.  No attempt was made to classify the components of $\mathscr N$.  Outside type $A$, and particularly in type $E_6,E_7$ for the ``Heisenberg'' parabolic, the existence of an adapted pair leads to an ``improved upper bound'' \cite [Lemma 6.11, 6.14]{J2} showing that $\mathscr I$ is polynomial!  Notably this method has been extended By Fauquant-Millet \cite {F} for other parabolics in classical type.  The original intention had been to recover the result of Yakimova \cite {Y} (proved partly by computer calculations) that polynomiality fails for the Heisenberg parabolic in type $E_8$.

\subsection {Adjoint Action}\label {1.4}

 Let $\mathfrak g$ be the Lie algebra of $G$. Fix a Borel subgroup $B$ of $G$  and let $\mathfrak n$ be the nilradical of the Lie algebra $\mathfrak b$ of $B$.  Fix a standard parabolic subgroup $P$ of $G$, that is to say one containing $B$, and let $\mathfrak m$ be the nilradical of the Lie algebra $\mathfrak p$ of $P$. Let $H$ be a Cartan subgroup for $B$ and $\mathfrak h$ its Lie algebra. Let $W$ be the Weyl group for the pair $(G,H)$.  Let $P'$ be the derived group of $P$.

  A theorem of Richardson asserts that $P$ has a dense orbit in $\mathfrak m$.  This implies \cite [2.2.2]{FJ1} that $\mathscr I:=\mathbb C[\mathfrak m]^{P'}$ is polynomial.   One calls $\mathscr I$ the semi-invariant subalgebra for the adjoint action of $P$ (or $\mathfrak p$) on $\mathfrak m$.  It, and the nilfibre $\mathscr N$ it defines, will be the main objects of our study.

   One may take the generators of $\mathscr I$ to be $H$ eigenvectors which are irreducible as polynomials.  Thus $\mathscr I_+$ is stable under $P$ and so is $\mathscr N$. Since $P$ is connected, the components of $\mathscr N$ are also $P$ stable.  However a $P'$ orbit through a point of $e+V$, need not be a $P$ orbit, and the latter could even be dense in $\mathfrak m$.

In the present paper we restrict our attention to $G$ being of type $A_{n-1}$, that is when $G=SL(n)$.

  \subsection {Type $A$}\label {1.5}

  For the adjoint action of a parabolic subalgebra on its nilradical, the method of adapted pairs fails miserably even in Type $A$.  First of all, a component of $\mathscr N$ may fail to admit a dense $P$ orbit \cite [Lemma 6.10.7]{FJ2} (though examples are not easy to find).  Even when a dense orbit exists, an adapted pair need not exist (such examples are easy to find).  This means here that (ii) of \ref {1.2} will fail and we must be more clever in establishing (i).  This was carried out for the ``canonical component'' as a consequence of \cite [Prop. 6.10.4(iii)]{FJ2} - the possible generalisation of which we shall not attempt here.

  In \cite [Thm. 2]{FJ1} we proved the existence of a Weierstrass in type $A$.  This was clarified in \cite {FJ2},\cite {FJ3} but not made any simpler.  Yet these papers led in \cite {FJ4} to a dramatic simplification through the construction of a ``composition tableau'', which surely could not have been achieved without the help of our preliminary work.

 \subsection {Primary decomposition}\label {1.6}

  Recall  that the (irreducible) components of an affine algebraic variety $\mathscr V$ are defined as follows.

  First the ideal $\mathfrak a$ of definition of $\mathscr V$ is an ideal in a polynomial ring $R$, so in particular noetherian.

 Yet in a noetherian ring $R$ every ideal $\mathfrak a$ admits a primary decomposition\footnote {Here of course we only need to know this for a polynomial ring which is due to Lasker, whilst for a noetherian ring it is due to Noether.}, that is to say is a finite intersection of primary ideals \cite [Lemmas 7.11, 7.12]{AM}. The radicals of these primary ideals are  prime and are independent of the primary decomposition of $\mathfrak a$ \cite [Thm. 4.5]{AM}. Those which are minimal (primes over $\mathfrak a$) define the (irreducible) components of $\mathscr V$.  (Although a noetherian ring does not satisfy the descending chain condition, its prime ideals do by virtue of Krull dimension (see for example \cite[Lemma 3.5.1]{D}).  Thus for an ideal of a noetherian ring one may directly define the minimal primes over it.)

 This beautiful general theory tells us little of how to find the irreducible components - just as we do not know how to decompose a number into its prime factors.

  \subsection {Component Tableaux and Main Results}\label {1.7}

 Let $\mathscr N$ be defined as in \ref {1.4}.

 Here we propose a method for finding the irreducible components of $\mathscr N$ involving the generalisation of the composition tableau.  These are the component tableaux.  Their construction is given in Section \ref {3}.

 We believe the component tableaux to be completely new. They surely could not have found without our preliminary work in \cite {FJ1},\cite {FJ2},\cite {FJ3},\cite {FJ4}\footnote{Even using chatGPT!}.

 In Section \ref {4} we assign to any component tableau a selected set of excluded roots Then the required component is the $B$ saturation set of the space spanned by the complement of the excluded root vectors.

 Here the choice of excluded roots is rather delicate for otherwise the said $B$ saturation set would be too small, but it also must be large enough so that the Benlolo-Sanderson invariants vanish when these excluded roots are set equal to zero.

 This last result requires showing that a triply infinite sequence of determinants, growing exponentially in $n$, all vanish, when the excluded roots are set equal to zero.

 \

  \textit {The construction of the excluded roots and the resulting vanishing is mainly due to the first author, with the second author following her intuition with blind faith.}

  \

 %Moreover the set of excluded roots has to be rather small making the analysis extremely tight.

 In Section \ref {5} we show that each component tableau also gives rise to a Weierstrass section.

  In \ref {6.4}, we construct as indicated above, a map from component tableaux to components, which we call the component map.

   %This uses the fact that  the set of excluded roots is rather small making the proof of vanishing of the Benlolo-Sanderson invariants extremely tight.

   It Section \ref {7} we prove an Exchange Lemma  \ref {7.2}.  Combined with partial linearity, Corollary \ref {7.4}, this proves that the component map is injective, Proposition \ref {7.5}.
   %In Section \ref {8}, surjectivity is also verified in a number of special cases.

   In subsequent work we tackle surjectivity which will complete the description of components.  Though this has yet to be finished we can already say it will involve a new family of tableaux (Reverse Tableaux) which will sequentially describe a successive factorisation of invariants, so that Krull's theorem can be applied to deduce surjectivity in a manner mirroring the special cases described in Section \ref {8}.

  \subsection {Orbital Varieties}\label {1.8}

  A co-adjoint $G$ orbit $\mathscr O$ for which $\mathscr O \cap \mathfrak n$ is non-empty is called a nilpotent orbit. After Spaltenstein \cite {S}, the intersection $\mathscr O \cap \mathfrak n$ is equidimensional of dimension $\frac{1}{2}\dim \mathscr O$.  Its components are called orbital varieties.  Orbital varieties are Lagrangian subvarieties for the Kirillov-Kostant symplectic form on $\mathscr O$ and as such can in principle be quantized as highest weight modules, an open problem of some difficulty, partly because quantization can sometimes fail (see \ref {2.3.4}).

% Through the Steinberg triple variety one may show that there is a surjective map $w \mapsto \overline{B.(\mathfrak n\cap \mathfrak n)^w}$ of $W$ onto the set of orbital variety closures.  In type $A$, the Robinson-Stensted correspondence gives a rather precise description of the fibres of this map.  This description was given by Steinberg though rather belatedly as the fibres had already been described by Spaltenstein via result of Knuth, the latter having admitted to having proved just for fun.

  \subsection {Hypersurface Orbital Varieties}\label {1.9}

    When we take the generators of $\mathbb C[\mathfrak m]^{P'}$ to be $H$ eigenvectors which are irreducible as polynomials, they determine a hypersurface in $\mathfrak m$ which is an orbital variety closure \cite [2.3.4] {FJ1}. We call this a hypersurface orbital variety.  Further details may be found in \cite [2.3.4] {FJ1} and references therein.

   Through the Steinberg triple variety one may show that there is a surjective map $w \mapsto \overline{B.(\mathfrak n\cap \mathfrak n)^w}$ of $W$ onto the set of orbital variety closures.  In type $A$, the Robinson-Schensted correspondence gives a rather precise description of the fibres of this map.  This description was given by Steinberg \cite {St} (see also references therein) though rather belatedly as the fibres had already been described by Spaltenstein via a result of Knuth who admitted (privately to my colleague A. Melnikov) to having proved this last result just for fun.

   Thus in type $A$ orbital varieties can be viewed as being known.

   Those which are hyper-surfaces orbital varieties are less well-known, particularly outside type $A$ and certainly in type $E$.

\subsection {Components as Orbital Varieties}\label {1.10}

\textit{Not} all the components of the nilfibre are orbital variety closures.  Perhaps the simplest example arises for the composition $(2,1,1,1,2)$.  In this case (see \ref {8.2}) the nilfibre has exactly three components, one is an orbital variety closure, the other two are not!

\subsection {Comments}\label {1.11}

We claim that our proofs are as tight and intricate as Chinese wooden box puzzles   which seem impossible to take apart\footnote{Even making the proofs fit for human consumption was a major challenge.}.  Indeed it was already noted in \ref {1.7} that the choice of the set of excluded roots is rather delicate. Again to define a Weierstrass section (see Section \ref {5}) for a given component, we select  a rather special monomial in each Benlolo-Sanderson invariant.  Each factor must be labelled by a $1$ or by a $\ast$ and just the latter are in the set of excluded root vectors.  Yet again every excluded root vector which is not a root labelled by a $\ast$ must have (see \ref {6.2.7}) to its left a root vector labelled by a $1$, whilst a certain root vector labelled by a $1$ must have no excluded roots in the upper right hand quadrant it defines (see \ref {7.3}), so in particular no excluded roots to its right.

It is surely magical that all this can be accomplished for all (standard) parabolics in Type $A_n$, which grow exponentially in $n$ and then rather mysteriously through a set of tableaux, duly earning the designation ``magical and mysterious''. In view of the remarks in \ref {4.7} we could well add ``miraculous''.

%We hope that our work will arouse the admiration and curiosity of every reader.  Besides we believe it can inspire much further work on components of the nilfibre to tackle surjectivity (of the component map), general type, biparabolics, co-adjoint action and reductive action, where almost no-one has dared to tread.

The construction of component tableaux came quite out of the blue and surely could not have been achieved without our preliminary wok in \cite {FJ2} culminating in the composition tableau \cite {FJ4} of which the component tableaux are a generalisation. Thus through no fault of our own we seem to have stumbled on a construction which we believe will arouse the admiration and curiosity of every reader.  Moreover many may be eager to extend these results to all simple Lie groups, which for the present, does not seem to have any obvious approach.

\subsection {Wider Horizons}\label {1.12}

A component of the nilfibre $\mathscr N \subset \mathfrak m$ is $P$ stable and so its ideal of definition is a $P$ stable prime ideal of $\mathbb C[\mathfrak m]=S(\mathfrak m^*)$.  This is a much smaller set than the set of all prime ideals and so in principle amenable to classification.  Had we just considered a single Lie group $G$, with Lie algebra $\mathfrak g$, being in addition self-dual (so semisimple) then the $G$ invariant prime ideals of $\mathbb C[\mathfrak g]=S(\mathfrak m)$ are related to the prime ideals of the enveloping algebra $U(\mathfrak g)$, which are in turn related its primitive ideals, which are simply annihilators of simple modules.  This brings representation theory into play and not surprisingly yields vast quantities of further results - see \cite {D} and \cite {Ja} for starters.

This puts the present study into a much wider context.

%Besides we believe it can inspire much further work on components of the nilfibre to tackle surjectivity (of the component map), general type, biparabolics, co-adjoint action and reductive action, where almost no-one has dared to tread.

%Finally our work is mainly combinatorial and it is certainly a huge challenge to find more intrinsic proofs.

\section {The Main Players.} \label {2}

Assume $\mathfrak g$ of type $A_{n-1}$ that is to say $\mathfrak {sl}(n)$.

\subsection {Columns}\label {2.1}

\subsubsection {Column Notation} \label {2.1.1}

 In type $A_{n-1}$, a (standard) parabolic subalgebra is given by the composition $(c_1,c_2,\ldots,c_k)$ of $n$ with the $c_i:i\in [1,k]$ being the successive sizes  of the Levi factors (see \cite [4.1.1]{FJ1}, which of course sum to $n$. This gives rise to a diagram $\mathscr D$ having  columns $C_1,C_2,\ldots,C_k$ of $\height C_i=c_i: i \in [1,k]$.

 Let $R_i:i=1,2,\ldots,$ denote the rows $\mathscr D$.  We view them as \textit{descending} down the diagram. Set $R^s=\cup_{i=1}^s R_i$.  Set $b_{i,j}:=C_i \cap R_j$, called the $ij^{th}$ box of $\mathscr D$.

 \subsubsection {Adjacent and Neighbouring Columns} \label {2.1.2}

  Successive columns of $\mathscr D$ are said to be adjacent. Two boxes in the same row,  in not necessarily adjacent columns of $\mathscr D$, with no boxes between them are called adjacent.

   Two columns of $\mathscr D$ are said to be neighbouring \cite [4.1.2]{FJ1} of height $s$, if they are both of height $s$ and there are no columns of height $s$ strictly between them.  After Melnikov \cite {M1}, the set of hyper-surface orbital varieties is in bijection with set of neighbouring columns.

   The set of columns of height $s \in \mathbb N^+$ in $\mathscr D$ is denoted $C_1^s,C_2^s, \ldots, C_{r_s}^s$, going from left to right as the subscript increases.

    \subsubsection {Intervals} \label {2.1.3}

    If $C,C'$ are columns of $\mathscr D$ we shall always mean $C$ to lie to the left of $C'$.  Denote by $[C,C']$ the set of columns of $\mathscr D$ between $C,C'$. If $C$, $C'$, or both are omitted, then we denote the resulting set as $]C,C'],[C,C'[,]C,C'[$.

    \subsubsection {Rectangles} \label {2.1.4}

    If $C,C'$ are neighbouring columns of height $s$, we set $R^s_{C,C'}=R^s\cap [C,C'] $ viewed as a union of boxes with entries from $[1,n]$, often referred to as the rectangle defined by the columns.

    \subsection {Tableaux}\label {2.2}

   \subsubsection {Matrix Notation} \label {2.2.1}

   Let $\textbf{M}$ denote set of all $n \times n$ matrices.  As a Lie algebra under commutation it is $\mathfrak {gl}(n)$, whilst $\mathfrak {sl}(n)$ is the (simple) Lie subalgebra of matrices of trace zero. As is customary we shall barely distinguish between the two.

   The Levi factor $\mathfrak r$ of $\mathfrak p$ is just $\oplus_{i=1}^k \textbf{B}_i$, where $\textbf{B}_i$ is the $i^{th}$ Levi block of size $c_i$ centered on the diagonal.

    For each column $C_i$, let $\textbf{C}_i$ denote the rectangular block in $\mathfrak m$ lying above $\textbf{B}_i$.  We call it the $i^{th}$ column block. Its width is $c_i$ and its height $\sum_{j < i}c_j$.  In particular $\textbf{C}_1=0$.

    Then  $\mathfrak m:= \oplus_{i=2}^k \textbf{C}_i$ is the nilradical of $\mathfrak p$.

    Let $x_{i,j}$ denote the $ij^{th}$ entry of \textbf{M}.  They form its standard basis.

\subsubsection {Number Insertion} \label {2.2.2}

      A tableau $\mathscr T$ is obtained from $\mathscr D$ by inserting in the boxes  the entries $1,2,\ldots,n$, going down columns and then from left to right.

    %  If $i',j'$ lie in boxes of $\mathscr T$, the line $\ell_{i',j'}$ joining them may be associated to $x_{i',j'}$.If $i'\neq j'$ also the root vector $\alpha_{i',j'}$.  \textit{We shall barely distinguish between the three}.

      If $i'$ and $j'$ occur in boxes of strictly successive columns, then $x_{i',j'} \in  \mathfrak m$, and these elements form  a basis for $\mathfrak m$.  Likewise $x_{i',j'}$, with $i',j'$ entries joined in the same column $C_i$, form a basis for $\textbf{B}_i$, properly speaking in $\mathfrak {gl}(n)$.

      If $i'\in [1,n]$, then it lies in a unique column (resp. row) of $\mathscr T$  which we denote by $C_{(i')}$  (resp. $R_{(i')}$).

 \subsection {The Benlolo-Sanderson Invariants}\label {2.3}

 Fix a pair of neighbouring columns $C,C'$ of height $s$.  Through the reduction of \cite [3.2, 4.1.2] {FJ1} the invariant assigned to a pair of neighbouring columns is obtained by ignoring all the columns outside the given neighbouring pair and is then given by the construction of Benlolo and Sanderson \cite {BS}.   All this is detailed in \cite [Sect. 3, Sect. 5]{FJ1} where some independent proofs are given.  Here we note a slight refinement in the lemma below,

 \subsubsection {Construction} \label {2.3.1}

 Assume that $C=C_1,C'=C_k$, then $\textbf{B}_1$  (resp. $\textbf{B}_k$) are the first (resp. last) Levi blocks and both have size $s$.  Let $\textbf{M}^*_s$ be the $n-s \times n-s$ minor in the lower left hand corner of \textbf{M}.  It fits in \textbf{M} snugly between these first and last Levi blocks.  It is semi-invariant for the opposed subalgebra $\mathfrak p^-$ of $\mathfrak p$.

 Through the Killing form we may regard $\textbf{M}^*_s$ as a (polynomial) function on $\mathfrak g$.  Amazingly when $\textbf{M}^*_s$ is restricted to $\mathfrak m$ it becomes a $\mathfrak p$ semi-invariant.  In this one only has to check invariance under $x_{s-1,s},x_{n-s,n-s+1}$.

 However if we apply this recipe blindly, the resulting invariant is zero whenever the Levi has a block size $>s$.

  \subsubsection {The Cure for Vanishing} \label {2.3.2}

 The trick (motivated by quantum groups) is to evaluate $\textbf{M}^*_s$ instead on $\mathfrak m +a\Id$.  Then it is non-zero, a polynomial in $a$ but not a semi-invariant!

 Let $d^{R_{C,C'}^s}_\mathscr D$ ( or simply, $d_\mathscr D$) be the number of boxes in $\mathscr D$ strictly below $R_s$ between the pair $C,C'$.

With respect to \cite [3.4.10] {FJ1}, one easily checks that $d_\mathscr D=\sum_{i=1}^k \max (c_i-s,0)$.

 The restriction of $\textbf{M}^*_s$ is a polynomial of degree $\sum_{i=1}^{k-1}\min (c_i,c_1)$.  One checks that the leading term is semi-invariant and is just $a^{d_\mathscr D}$ times what is then the Benlolo-Sanderson invariant $I^s_{C,C'}$ \cite {BS}.  It is an irreducible polynomial \cite [5.3]{FJ1}.  If $C,C'$ have height $s$ but are not neighbouring, then the resulting polynomial is still semi-invariant but is a product \cite [1.10]{FJ2} of the factors coming from the pairs of neighbouring columns of height $s$, sandwiched between $C,C'$.

 \subsubsection {The Variables in a Benlolo-Sanderson Invariant} \label {2.3.3}

 Let $\mathfrak m^-$ be the transpose of $\mathfrak m$.  We may identify $\mathbb C[\mathfrak m]$ with $S(\mathfrak m^-)$ through the Killing form.  Then $I^s_{C_1,C_k} \in S(\mathfrak m^-)$.  More generally $I^s_{C,C'}$ is a polynomial in the co-ordinate functions $x^*_{i,j}:x_{i,j}\in \mathfrak m$.  Not every co-ordinate function need occur in $I^s_{C,C'}$ (as indicated just below).

 \subsubsection {Quantization} \label {2.3.4}

   Consider the composition $(1,2,1)$.  In this case $I^1_{C_1,C_3} = x^*_{1,2}x^*_{2,4} +x^*_{1,3}x^*_{3,4}$.  Its zero locus is a hypersurface orbital variety in $\mathfrak m$.  On the other hand the term which multiplies $a^2$ in the above construction, namely $x^*_{1,4}$ is \textit{not} a semi-invariant.

   Consider the $x_{j,i}:x_{i,j}\in \mathfrak m$.  They commute in $S(\mathfrak m^-)$.  Some may wish to consider them in the enveloping algebra  $U(\mathfrak m^-)$, where they do not commute. Then one may ask if there is a highest weight module, whose highest weight vector is annihilated by a suitably ordered form of this expression.  Thereby Kostant \cite {K2} first constructed a non-trivial ``quantization'' in type $A_3$ which turned out to be a unitarizable highest weight module for a suitable real form.

 The problem of quantization of hypersurface orbital varieties in type $A$ was solved in  \cite {JM} with the help of the Jantzen sum formula.  This drastically simplified and generalized the work of Kostant mentioned above \cite {K2}, though these modules are unitarizable if and only if the $\mathfrak m$ is commutative and then can be uniformly presented \cite {EJ}, a rarity in the theory of unitarizable modules.

  The work in \cite {JM} was generalized to classical type by Elena Perelman \cite {Pe}, for which one has to weaken the definition of quantization.   Elena was a doctorial student of the second named author but sadly she followed family tradition and left her superb work unpublished.  It has barely been approached let alone surpassed for now more than twenty years.

 A component of $\mathscr N$ is not in general an orbital variety closure, therefore not Lagrangian and so cannot be quantized. However it can be an orbital variety even if it is not a hypersurface.  The quantization of such orbital varieties has not been considered, other than in the general context of the quantization of all orbital varieties which is a difficult problem and is not always possible \cite {J4}.

  \subsubsection {The Distribution of the $a$ factors} \label {2.3.5}

  The formula for $d_\mathscr D$ suggests that each Levi block of size $t>s$ contributes a factor of $a^{t-s}$.  This is shown through the following relatively easy lemma in which we can assume $C=C_1,C'=C_k$ without loss of generality.

\begin {lemma} Let $\Id_{>s}$ be the diagonal matrix with only $a$ in the entries of the blocks of size $>s$.  Then the restriction of $\textbf{M}^*_s$ to $\mathfrak m +a\Id_{>s}$ has leading term $a^{d_\mathscr D} I^s_{C,C'}$.  Moreover in this the  $i^{th}$  Levi block contributes a factor of $a^{\max(0,c_i-s)}$.

\end {lemma}

\begin {proof}

Consider the sub-block $\textbf{M}_j$ of \textbf{M} enclosing the first $j$ Levi blocks.  The number of rows (resp. columns) of $\mathfrak m^-$ lying in $\textbf{M}_j$ is $\sum_{i=2}^j c_i$ (resp.  $\sum_{i=1}^{j-1}c_i$).  Recall that $c_1=s$. Suppose $c_j>s$ and choose $j$ minimal with this property.  Then the number of rows exceeds the number of columns in $\mathfrak m^-$ lying in $\textbf{M}_j$ by $c_j-s$.  Consequently the restriction of $\textbf{M}^*_s$ to the direct sum of the first $j$ column blocks $\oplus_{i=1}^j\textbf{B}_i$ is zero.   However if we add to this direct sum the diagonal matrix $a\Id$ at the place of the $j^{th}$ Levi factor, then the restriction is non-zero with leading term having a factor of $a^{c_j-s}$.  Moreover the contributions to this leading term are obtained by placing exactly $c_j-s$ copies of $a$ on the diagonal of the $j^{th}$ Levi block in every possible way.  Then in each contribution the number of columns in the $j^{th}$ column block $\textbf{B}_j$ has been reduced from $c_j$ to $s$.
Consequently the above argument may be repeated for every subsequent Levi block of size $>s$, eventually giving the assertion of the lemma.

%The number of rows (resp. columns) of $\mathfrak m^-$ meeting the first $j$ (resp. $j-1$) Levi blocks is $\sum_{i=2}^j c_i$ (resp.  $\sum_{i=1}^{j-1}c_i$).  Suppose $c_j>s$ and choose $j$ minimal with this property. Then the sum of the rows exceeds the sum of the columns by $c_j-s>0$. Consequently the restriction of $\textbf{M}^*_s$ to the direct sum of the first $j$ column blocks $\oplus_{i=1}^j\textbf{B}_i$ is zero.  However if we add to this direct sum the diagonal matrix $a\Id$ in the place of the $j^{th}$ Levi factor, then the restriction is non-zero with leading term having a factor of $a^{c_j-s}$.  Moreover in this leading term, the number of columns below the $j^{th}$ Levi factor has been reduced from $c_j$ to $s$.

%Thus the above argument may be repeated, reducing $j$ by $1$ and eventually gives the assertion of the lemma.

\end {proof}

\section {The Component Tableaux.} \label {3}

Assume $\mathfrak g=\mathfrak {sl}(n)$, that is of type $A_{n-1}$. Recall (\ref {2.2.2}) that a (standard) parabolic subalgebra $\mathfrak p$ determines a tableau $\mathscr T$.  The latter gives a basis for the nilradical $\mathfrak m$ and encodes the Benlolo-Sanderson generators of the polynomial algebra $\mathbb C[\mathfrak m]^{P'}$, whose common zero loci define the nilcone  $\mathscr N$.  Our aim is determine the components of the latter through ``component tableaux''.

In this section a component tableau $\mathscr T^\mathcal C$ will be constructed from $\mathscr T$  by a simple procedure which extends the construction of the composition tableau \cite [4.7]{FJ4}.  A main difference is that we can make several choices at each step, incorporated as numerical data $\mathcal C$ - \ref {3.1.5}.

In \ref {4.1.4} we use $\mathscr T^\mathcal C$ to define a Lie subalgebra $\mathfrak u^\mathcal C $ of  $\mathfrak n$, setting $\mathscr C:=\overline{B.\mathfrak u^\mathcal C }$.

%It is a non-trivial task \ref {4.2}, \ref {4.3} to show that $\mathscr C \subset N$ and then is a component of $\mathscr N$. The latter is achieved  by constructing a Weierstrass section (Section \ref {5}) attached to $\mathscr C$ and by a dimensionality estimate  (Section \ref {6}).
%
% This gives a map $\mathscr T^\mathcal C \mapsto \mathscr C$ of component tableaux to components, called the component map.  In \ref {7.5} the component map is shown to be injective.  Some evidence for surjectivity is given in Section \ref {8}.  Yet this is an entirely different task, wisely postponed to a subsequent paper.

\subsection {Generalities and Notation}\label {3.1}

\subsubsection {An Overview of the Sequential Construction of a Component Tableau} \label {3.1.1}

As in \cite [4.2,4.7]{FJ4} we augment the columns of $\mathscr D$ by increasing their height through adjoining empty boxes.  We denote this augmented diagram by $\mathscr D'$. Since we are just lengthening existing columns, there are no empty boxes in  the first row, that is in $R_1\cap \mathscr D'$.

  Modifying slightly the notation of \cite [4.7]{FJ4}, set $\mathscr T(1)=\mathscr T$ and for each $t \in \mathbb N^+$, let $\mathscr T(t)$ denote the tableau in which entries have been successively inserted into empty boxes of the first $t$ rows of $\mathscr D'$.

  %This successive insertion is carried out by moving entries to the right and down by rows.

  Set $\mathscr T(t,1)=\mathscr T(t-1)$ and let $\mathscr T(t,r+1):r \in [1,k[$ be obtained by putting an entry from $C_r$ into $C_{r+1}$ and down by $\geq 0$ rows, following the rules to be announced in \ref {3.2}.  (This may leave a ``gap'' in each row $R_t$ on its left (cf \cite [Fig. 2]{FJ4}), which can increase in length with $t$. It is of no particular importance.)

    The construction in \ref {3.2.2} implicitly selects elements in ``batches'' and may be carried out in \textit{several different ways} leading to \textit{several different choices} for $\mathscr T(t)$, for each $t \in \mathbb N^+$.

  Fixing for the moment just one particular selection in each batch, $\mathscr T(t)$ is a sub-tableau of $\mathscr T(t+1)$ and we denote their direct limit by $\mathscr T(\infty)$. Our rules imply that $\mathscr T(t) = \mathscr T(\infty)$, for all $t$ sufficiently large.

Given $C$ a column of $\mathscr T$.  View $C(t)$ (resp. $C(t,r), C(\infty)$), as the corresponding column of $\mathscr T(t)$ (resp.  $\mathscr T(t,r),\mathscr T(\infty)$).

   The rows $R_{t}$ of $C$ coincide with those of  $C(t)$ (resp. $C(\infty)$) for $t \leq \height C$. One has $\height C(t)\geq t$  and correspondingly $\height C(t,r) \geq t-1$ to the right of the above mentioned gap. Notice that this is a strict inequality if and only if $\height C \geq t$.

 \subsubsection {The Vav Conversive} \label {3.1.2}

  It is extremely inconvenient and not particularly informative to carry, with studious diligence, the notation $\mathscr T(t), \mathscr T(t,r)$ throughout.

  To avoid this as much as possible we shall use the principle of the Vav Conversive \cite [p.95]{Da} used in ancient texts, which arose from an intermingling of past and future as might result from sojourning in the desert when every day seems the same. Moreover here ambiguity is reduced because the past does not change.

  Thus instead of referring to $\mathscr T(t)$ being filled to give $\mathscr T(t+1)$, we shall view the rows of $\mathscr T(\infty)$ as having been successively completed. Since ``past'' rows are never changed the possible confusion is limited.

  % In some sense this principle has already been incorporated in not giving the step by step construction of $\mathscr T(t)$ from $\mathscr T(t-1)$.

   Yet of \textit{prime} importance is to distinguish the column $C$ of $\mathscr T$ with the corresponding column $C(\infty)$ of $\mathscr T(\infty)$ which may have added entries.

   \subsubsection {Surrounding Columns} \label {3.1.3}

   Recall the notation of \ref {2.1.1}. Two adjacent columns $C:=C_r,C':=C_{r+1}:r \in [1,k-1]$ of $\mathscr D$ are said to be surrounded by a pair of neighbouring columns $C(s),C(s)'$ of $\mathscr D$ of height $s$ if $C,C' \in [C(s),C(s)']$.

   This may arise for a finite set $S$ of values of $s$ simultaneously.  Here \textit {we do not } impose any condition on the relative order of the columns $C(s):s \in S$  (resp. $C(s)': s \in S)$.  Generally  we take $S$ of the form $[m_1,m_2]: m_1\leq m_2\in \mathbb N^+$.

    \subsubsection {Batches} \label {3.1.4}

    In the following (particularly in \ref {3.2.2}) we observe that \textit{several different} $\mathscr T(\infty)$ can be constructed and specified by a choice of  numerical data $\mathcal C$.   In the language of \ref {3.1.2} although the past is fixed, the future need not be, creating several different tableaux $\mathscr T^\mathcal C(\infty)$ in the limit\footnote{Just as certain theoretical physicists have been known to postulate parallel universes, even infinitely many! In our case there will only be finitely many  $\mathscr T^\mathcal C(\infty)$ choices for a given $\mathscr T$.}.  This might lead to confusion, so we postpone the use of the vav conversive up to the end of \ref {3.2.2}.   For ease of notation we shall sometimes omit the superscript.  This is compatible with the notation in \ref {3.1.1}.

  \textit{Unlike} \cite [4.7] {FJ4} we shall \textit{not} insert the elements of $[1,n]$ with the smallest first for the (linear) order relation $\preceq$ defined in \cite [4.7]{FJ4}.

  Instead we define disjoint subsets of $R_t\cap \mathscr T(t,r)$, called batches,
  used to pass from $\mathscr T(t,r)$ to $\mathscr T(t, r+1)$  (in one of several ways).  Here we regard the pair $(t,r)$ as being fixed at least up to the end of \ref {3.2.2}.

  %In this $\mathscr T(t)$ is unchanged, so there is no confusion using the principle of the Vav Conversive to regard them as lying in $\mathscr T(\infty)$.

  Let $C,C'$ be neighbouring columns of height $t$. As indicated in \ref {3.1.1}  a given row of $\mathscr T(t,r)$ may admit repeated entries having the same value.

  The batch $\mathscr B^{s}_{C,C'}:s\leq t$ is defined to consist of the rightmost entries having a given value $j$ in $R_{s}$  in $[C,C'[$ in $\mathscr T(t,r)$.

  %(Note that because $\height C'=t$ the box $R_{t'}\cap C'$ is already full, so an entry in $R_{t'}\cap [C,C'[$ cannot be pushed rightwards and horizontally into $R_{t'}\cap C'$).

  In the notation of \ref {2.1.2}, for all $u\in [1,r_s-1]$, we shall write $\mathscr B^s_{C^t_u,C^t_{u+1}}$ as $\mathscr B^s_u$.
  Set $\mathscr B^s =\cup_{u=1}^{r_s-1}\mathscr B^s_u$. The elements of this union we call the batches of $\mathscr B^s$.

  %We shall see that by \ref {3.2.1} $(iii)$ below implies that

  %$(*)$.  A given  $i\in \mathscr B^s$, lies in a unique $\mathscr B^s_u$.

      \subsubsection {Choices in Batches} \label {3.1.5}

      Through our construction every batch will provide exactly one entry.

      This provides the numerical data collated in $\mathcal C$ given by the ordered sequence of entries as one descends rows. Within a row no ordering is necessary.

      \

       \textbf{Definition}.  A pair of neighbouring columns is said to be free if it has not been used in a previous step.

\subsection {Detailing our Rules}\label {3.2}

   \subsubsection {Stopped Entries }\label {3.2.1}

   %{Stopped Entries }\label {3.2.1}
   %To avoid ambiguity we do not use the Vav Conversive in this paragraph.

   \

   \textbf {Definition}.  An entry of $C_r(\infty)$  is said to be stopped at $C_r$ if it does not appear in any column strictly to the right of $C_r$ in $\mathscr T(\infty)$.

     \subsubsection {Generalizing \cite [4.7]{FJ4}}\label {3.2.2}

     %To avoid any ambiguities we avoid the use of the vav conversive in this section.

   \

    % Set $\mathscr T(1)=\mathscr T$.

    Recall the notation of \ref {2.1.1}  and take $r+1\in [1,k]$.

    \

    Rule $(1)$.

     $(a)$. Take $i \in R_{t'} \cap C_r(t,r)$ for some $t'\leq  t$ and $\height C_{r+1}(t,r)=t$.  (Here $i$ is \textit{not necessarily} the lowest element of  $C_r(t,r)$.)

    $(b)$.   Assume that the adjacent pair $C_r,C_{r+1}$ is surrounded by free pairs of neighbouring columns of all the heights $s\in [t',t]$.

    Then put $i$ in $R_{t+1}\cap C_{r+1}(t,r)$, that is to say put $i$ down by $t-t'+1$ rows into  $C_{r+1}(t,r)$, to obtain $C_{r+1}(t,r+1)$.

    \

  Rule $(2)$.  Suppose $(b)$ fails in Rule $(1)$, then $i$ is stopped at $C_r$.

  \

  Rule $(3)$. If in Rule $(1)$ one has  $\height C_{r+1}(t,r)\leq t-1$, so equality holds by the last paragraph of \ref {3.1.1}. Then put $i$ in $R_{t}\cap C_{r+1}(t,r)$ to obtain $\mathscr T(t,r+1)$.

  \

   \textbf{Comments}.

   Rule $(1)$ puts $i$ in the batches $\mathscr B^s_{u_s}: s \in [t',t]$ given by the \textit{uniquely determined} pairs of neighbouring columns $C^s_{u_s},C^s_{u_s+1}$ surrounding $C_r,C_{r+1}$.

     Rule $(1)$ is the same as \cite [4.7(iii)]{FJ4}, if $\height C_{r+1}=t$ and has a left neighbour of height $t$ and that $t=t'$.  Yet if $\height C_{r+1}= t+1$, Rule $(1)$ is not \cite [4.7(ii)]{FJ4} which says that $i$ is stopped at $C_r$.

     Rule $(2)$ incorporates (part of) \cite [4.7(ii),(iv)]{FJ4}.

     Rule $(3)$ just moves $i$ horizontally by one column and rightwards in the row $R_{t+1}$.  It is rule \cite [4.7(i)]{FJ4}.

     \textit{Notice} that \textit{unlike} \cite [4.7]{FJ4} we have put this Rule $(3)$ \textit{last}.  Doing this incorporates the ordering $\preceq$ used in loc cit.

   \

   \textbf{Observations}.

   \

   $(i)$.  As $(b)$ of Rule $(1)$ requires \textit{free} pairs of neighbouring columns, it is \textit{automatic} that only one entry in each box is chosen.

   \

   $(ii)$.  There is \textit{always} one entry in each batch $\mathscr B^t_u:, t \in \mathbb N^+, u\in [1,r_s-1]$ that can be used in Rule $(1)$, namely when $i \in \mathscr B^t_u$ is farthest to the right. Recalling the last sentence of \ref {3.1.2}, $\height C^t_{u+1}=t$ \textit{in} $\mathscr T$, so $i$ can be inserted into $R_{t+1}\cap C^t_{u+1}$.

   \

   \textbf{Example $1$.}    Consider the composition $(1,2,1)$. We must first lower $3$ from $C_2$ into $R_2\cap C_3$. Then $3$ can be translated horizontally into $R_2\cap C_3$.  Yet if we do the second operation first, then $2$ is blocked \textit{and we fail to use the neighbouring columns of height $1$ and so fail to incorporate the corresponding Benlolo-Sanderson invariant}.  This example was also noted in \cite [4.7]{FJ4} and there $\preceq$ was used to prevent ``blockage''.

   \

  %Since $\height C^t_{r+1}=t$, Rule $(3)$ does not allow $i$ to be pushed into $R_t\cap C_{r+1}^t$ which already has an entry. This gives $(*)$ of \ref {3.1.4}.

  In our new set-up, one may have $\height C(\infty) > \height C$, for some columns $C$, for two different reasons. First if $t>t'$ in $(b)$ of Rule \ref {3.2.2}.

   \

   \textbf{Example $2$}.   Consider the composition $(2,1,2,1)$.  Then $3$ belongs to both $\mathscr B^1_{2,4}$ and $\mathscr B^2_{1,3}$.  In view of the available surrounding columns, $3$ may be lowered by two rows into $R_3$ below $5$. This differs from the composition tableau in which $4$ is lowered into $R_4$ below $6$ and then $2$ is moved horizontally into $C_2$ and then into $R_3$ below $5$.

   \

    Secondly  if already  $\height C_{r+1}(t,r)>\height C_{r+1}$.

    \

    \textbf{Example $3$.}  Consider the composition $(1,2,1,2)$.  The two possible component tableaux are given in Figure $2$ with (momentarily looking ahead to \ref {3.2.6}) the decoration of lines given in Figure $3$.
  \

  \

  \textbf{N.B. 1}. By Rules $(1)-(3)$, the entry $i$ of $R_{t'}\cap C_r(t,r):t'\leq t$ may only be placed in the first empty box of $C_{r+1}(t,r)$.  Thus in $\mathscr T(\infty)$ only one integer can fill a given box of $\mathscr D'$.

  Again there are no gaps in columns of $\mathscr T(\infty)$ and no gaps in rows except the gap to the left mentioned already in \ref {3.1.1}.  Finally $\mathscr T(\infty)$ is semi-standard, though not quite in the usual sense \cite [4.8]{FJ4}.

  \

  \textbf{N.B. 2.}  \textit{All} the batches of $\mathscr B^{s}:s \in [t',t]$ are only fully determined when $\mathscr T(t)$ is completed, because as entries are lowered batches change.  This makes it difficult to compute  $|\mathcal C|$, except algorithmically case by case.

  %Yet as noted in \ref {3.1.5} some elements in these batches may be determined before this step and if they are selected no further elements in these batches may be chosen, because of our condition on the pairs of neighbouring columns being free.  This makes it difficult to compute  $|\mathcal C|$, except algorithmically case by case.

  The limiting tableau $\mathscr T(\infty)$ depends on $\mathcal C$ and should be more fully denoted by $\mathscr T^{\mathcal C}(\infty)$, though we shall generally desist from doing so to avoid cumbersome notation.

  \

   $(**)$. Since a pair of neighbouring columns has to be free, it can be used only once in the above construction.

    \

All this may seem complicated to some readers.  Actually it is very natural.  Moreover we believe it well-nigh impossible to imagine these component tableaux otherwise and in the process to gain a horde of components, without some pain.

\

\textbf{Remark.} Notice that entries in $\mathscr T(\infty)$ are always translated to right adjacent columns.  As a consequence the boxes in a given column $C(\infty)$ have distinct entries, whilst this is false of rows.

  \subsubsection {Summarising the Passage from $\mathscr T(t)$ to $\mathscr T(t+1)$  } \label {3.2.3}

  \

     Fix $t \in \mathbb N^+$ .

     %For simplicity of notation omit the superscript $t$ on batches.

     Consider $\mathscr T(t)$ to have been constructed. Then complete the $(t+1)^{th}$ row of $\mathscr T(t+1)$.

     %lying   $[C^t_i,C^t_{i+1}]$, as follows.

       Recall the allowed choices of $i$ in the batches of  $\mathscr B^{s}:s \in [t',t]:t'\leq t$, given in Rule $(1)$ of \ref {3.2.2}.  In this each $i\in R_{t'}$  moves to the right across one column and down $t-t'+1$ rows, or by Rule $(2)$ may be stopped.

        After that the remaining entries in the $(t+1)^{th}$ row of $\mathscr T$ (so in $\mathscr T(t)$) are moved horizontally to the right sequentially starting from the leftmost one into the available empty boxes of $\mathscr D'$, through Rule $(3)$ of \ref {3.2.2}).  This last step is completely determined by the entries which have been lowered into $R_{t+1}$ by the steps in the paragraph above.

       \

   We shall see (Lemma \ref {3.2.5}) that this moving down by $t-t'+1$ rows  for an initial choice of $i\in \mathscr B^{s}:s \in [t',t]$  corresponds to there being $t-t'+1$ invariant generators defined by an appropriate set of $t-t'+1$ pairs of neighbouring columns surrounding the column $C$ containing $i$ and its right adjacent column $C'$ (see (ii) below).

 \subsubsection {Rules for the Construction of labelled lines in $\mathscr T^\mathcal C(\infty)$}\label {3.2.4}

  As in \cite [6.5]{FJ4} the tableau $\mathscr T^{\mathcal C}(\infty)$ is to be decorated by lines, each line $\ell_{i',j'}$ joining distinct entries of boxes.  These were obtained in loc. cit. by examining the composition map (see \cite [5.3]{FJ4}.  Here we merely \textit{impose} similar rules given below. Of course we are doing this for each choice of $\mathscr T^\mathcal C(\infty)$.  The miracle is that this procedure leads to a Weierstrass section (see Section \ref {5}) for each component of $\mathscr N$ in the image of the component map (so for every component, given that the component map is surjective).

  \

  \textbf{N.B.} To $\ell_{i',j'}$ we may associate the basis element $x_{i',j'}\in \textbf{M}$ and the root $\alpha_{i',j'}$.  In all what follows we shall barely distinguish between the three.  \textit{Yet we shall never add roots only root vectors.}

   \subsubsection {Lines with label $\ast$}\label {3.2.5}

  Recall that in the \textit{composition} tableau every ``step'' produces a vertical line which we decorate with a $\ast$.  There a step meant that if $C_r,C_{r+1}$ are adjacent tableaux with $\height C_{r+1}=t$ and $i\in R_{t}\cap C_r(t,r)$, then $i$  is moved right and downward by \textit{one row} into the empty box below the lowest box in $C_{r+1}$  through \cite [4.7(ii)]{FJ4}.

  Then a vertical line labelled with a $\ast$ is drawn from $i \in R_{t+1} \cap C_{r+1}(t,r+1)$ to the entry $j$ of the \textit{lowest } box in $C_{r+1}$.  %This produces a vertical line $r_{i,j}$ in $\mathscr T(\infty)$, with label $\ast$.

  \textit{However there is a novelty here} which arises because \cite [4.7(ii)]{FJ4} has been replaced by Rule $(1)$ of \ref {3.2.2}.  The latter allows \textit{three new} possibilities.

  First it can happen that $i \in R_{t'}\cap C_r(t,r)$ is moved right and downward by $t-t'+1$ rows into the first empty box below the lowest box $R_t\cap C_{r+1}$ in $C_{r+1}$  through Rule $(1)$ of \ref {3.2.2}.

  Then draw a vertical line labelled by a $\ast$ from $i \in R_{t+1}\cap C_{r+1}(\infty)$ to \textit{each} of the $t-t'+1$ entries $j_k$ in $R_{k}\cap C_{r+1}:k \in [t',t]$.  Recalling that $\height C_{r+1}=t$, one observes that these are the lowest $t-t'+1$ entries of $C_{r+1}$.

  \

  \textbf{Example 4.} Consider the composition $(2,1,1,2,1)$.  Then $\mathscr B^1_{2,3}=\{3\}, \mathscr B^1_{3,5}=\{4,5\}$.  Then in the first step $3,5$ are lowered into $R_2$ as shown in Figure $4$.  Then $\mathscr B^2_{1,4}=\{2,3,4\}$.  This gives three tableaux as illustrated in Figure $4$. In this, the second is the canonical tableau and in the third, $4$ is lowered by $t-t'+1=2$ rows from $R_1$ into $R_3$.

  \

  %For example in the third component of Example $4$ one has $t'=1,t=2$.

   Secondly this lowering may be repeated, but only in a more modest fashion.  Indeed $\height C_{r+1} = t$ and (to avoid cumbersome terminology) we shall simply say its height has been increased by one by placing $i$ in the first empty box below $C_{r+1}$, so in row $R_{t+1}$.

   In this the neighbouring columns of heights $s\in [t',t]$ surrounding $C_r,C_{r+1}$ have been used. As a consequence the entries of $C_r(\infty)$ in rows $R_u:u \in [t'+1,t]$ are stopped at $C_r$. On the other hand the new $C_{r+1}$ obtained from the step described in the previous paragraph has height $t$.  This means the entry of $R_{t+1}\cap C_r(\infty)$ need not be stopped at $C_r$  and then by Rule $(1)$ of \ref {3.2.2} may have gone (down by just one row) into $R_{t+2}\cap C_{r+1}(\infty)$ using the pair of neighbouring columns of height $t+1$ surrounding $C_r,C_{r+1}$.

   Thirdly this last step may be repeated several times.  This was already illustrated in Example $3$. Eventually this ``strange'' behaviour will appear to be natural and necessary

   \

   \textbf{Example 5.}   Consider the composition $(3,2,1,2,2,1,3)$.  Then $7 \in \mathscr B^1_{3,6} \cap \mathscr B^2_{4,5}$ and so can be lowered by two rows below $10$. Then $\mathscr B^2_{2,4}=\{5\}, \mathscr B^2_{4,5}=\{8\}$.  In this we have already used the pair $C_4,C_5$ of neighbouring columns of height $2$ and so $8$ is stopped at $C_4$.

   Yet we can move $5$ by one row under $8$ in $C_4$. Then $\mathscr B^3=\{3,5,7\}$, so $5$ can be lowered one row below $7$ in $C_5$. This exemplifies the second part with $t'=1,t=2$, as illustrated in Figure $6$.

   \begin {lemma}

    The number of lines in $\mathscr T^\mathscr C(\infty)$ with label $\ast$ is the number of pairs of neighbouring columns.
  \end {lemma}

  \begin {proof} This is essentially obvious except for  one subtle point.  First every pair of neighbouring columns is used at most once by $(**)$ of \ref {3.2.2} and each produce a vertical line in $\mathscr T(\infty)$ labelled by a $\ast$.

  Conversely we must show for all $s \in \mathbb N^+$ that every pair of neighbouring columns $C,C'$ is used in constructing a given component tableau.  Indeed by Observation (ii) of \ref {3.2.2} there is always one entry of $\mathscr B^s_{C,C'}$  which is a possible choice in \ref {3.2.2}, Rule $(1)$.  Such a choice is not  excluded by \ref {3.2.2}, Rule $(3)$ which is performed after \ref {3.2.2}, Rule $(1)$.  (See also Example $1$.)

  \end {proof}

 \subsubsection {The lines with label $1$}\label {3.2.6}

  Consider $C_r(\infty),C_{r+1}(\infty)$ as adjacent columns in $\mathscr T(\infty)$.

 The lines with label $1$ are deemed to arise when a box in $C_r(\infty)$ has been stopped at some $C_r$. Then the entry $i$ of that box is defined to be the left end-point of the resulting line $\ell$.  It is the rightmost appearance of $i \in \mathscr T(\infty)$.  The right end-point of $\ell$ is \textit{chosen} to be an element of $C_{r+1}$ defined inductively by descending its rows.  It is the highest available entry of the box in $C_{r+1}$ whose entry is not already a right end-point of a line starting higher up in $C_r(\infty)$, whose left-end point is stopped at $C_r$.  One easily checks that these lines are right and up-going.

 We also obtain lines joining boxes in adjacent columns of $\mathscr T(\infty)$ with the \textit{same} entry $i$.   We call these lines neutral, or having a neutral label.

 \

 \textbf{Trailer}. A preview of the labelled lines we obtain between $C_r(\infty)$ and $C_{r+1}(\infty)$ may be obtained from Figure $1$.

 \

 \textbf{Notation}. Given a column $C$ of $\mathscr T$ and $p \in \mathbb N^+$, let $C^{>p}$, be the partial column formed from the rows $R_{p'}:p'>p$ of $C$.  Let $C^{\leq p}$ be the truncated column obtained from $C$ by removing $C^{>p}$.  A similar notation is used for a column $C(\infty)$ of $\mathscr T(\infty)$.

 \

 Let us make four observations concerning the lines joining boxes in $C_r(\infty),C_{r+1}(\infty)$ with label $1$, which result from the above construction.

   \

   $(i)$.  Suppose for all $t',t \in\mathbb N^+:t\geq t'$ that $C_r(\infty)\cap R_{t'}$ and $C_{r+1}(\infty)\cap R_{t}$ do not have a common entry.

    \

    By the rules in \ref {3.2.2}, it follows that $\height C_{r+1}\geq t$ and an entry $i_{t'}$ of $C_r(\infty)\cap R_{t'}:t'\leq t$ is stopped at $C_r$ giving a horizontal line with label $1$ in $\mathscr T(\infty)$  to the entry $j_{t'}$ of $C_{r+1}\cap R_{t'}$.

    In other words we just obtain a set of horizontal lines with label $1$ joining $R_{t'}\cap C_r(\infty)$ to $R_{t'}\cap C_{r+1}$, for all $t' \in [1,t]$.

    \

    Otherwise

    \

 $(ii)$. let $t'$ be minimal such that $R_{t'}\cap C_r(\infty)$ and $ R_{t+1}\cap C_{r+1}(\infty)$ for some $t\geq t'$, share a common entry $i_{t'}$, so are joined by a neutral line.

  \

  In this recall \ref {3.2.2}, Rule $(1)$ that we must have had $\height C_{r+1}=t$. Then for all $m \in [t'+1,t]$ the entry $i_m$ of $R_m\cap C_r(\infty)$ is stopped at $C_r$, because as noted in \ref {3.2.5} the pairs of neighbouring columns of heights $s \in [t',t]$ are no longer free.

   If $t>t'$, this gives for all $m \in [t'+1,t]$, a line with label $1$ from $i_m$ to the entry $j_{m-1}$ of $R_{m-1}\cap C_{r+1}$.

  On the other hand as noted in \ref {3.2.5}, the entry $i_{t+1}$ of $R_{t+1}\cap C_r(\infty)$ need not be stopped at $C_r$ and  may  go down by one row entering $R_{t+2} \cap C_{r+1}(\infty)$.  Again these boxes are joined by neutral lines.  Moreover ``this going down by one row'' may be repeated several times, or not at all.  In the latter case the entry $i_{t+1}$ (\textit{if it exists}) of $R_{t+1}\cap C_r(\infty)$ is joined by a line with label $1$ to the entry $j_t$ of $R_t\cap C_{r+1}$.  Otherwise

  \

  $(iii)$.  Let $t''$ be maximal such that for all $s \in [t,t'']$ the boxes $ R_s\cap C_r(\infty)$ and $R_{s+1}\cap C_{r+1}(\infty)$ share a common entry $i_s \in [1,n]$.

  \

  The entry $i_{t''+1}$ of $R_{t''+1}\cap C_r(\infty)$ (\textit{if it exists}) is stopped at $C_r$ (because the box $R_{t''+1}\cap C_{r+1}(\infty)$ is filled by $i_{t''}$ and by hypothesis $i_{t''+1}$ does not enter $C_{r+1}(\infty)\cap R_{t''+2}$).  Let $j_{t}$ be the entry of $R_{t}\cap C_{r+1}$. This produces (at most) one right up-going line $\ell_{i_{t''+1},j_{t}}$ labelled by $1$ from $R_{t''+1}\cap C_r(\infty)$ to $R_{t}\cap C_{r+1}$.

  \

$(iv)$. For $t''>t+1$, the entries of $R_{t''}\cap C_{r}(\infty)$ and of $ R_{t''}\cap C_{r+1}(\infty)$  coincide (so are joined by neutral lines) or are both empty.  Thus there are no right-going lines with label $1$ from $C_{r}(\infty)$ to $C_{r+1}(\infty)$ other than those described in $(i)-(iii) $ above.

\

In the above notation we obtain

     \begin {cor}  In $\mathscr T^\mathcal C(\infty)$ right going lines with label $1$ join adjacent columns \newline $C_r(\infty),C_{r+1}(\infty)$.  They are horizontal through $(i),(ii)$ with the exception of the lines $\ell_{i_m,j_{m-1}}$ which are right and up-going by one row and the line $\ell_{i_{t''+1},j_{t}}$ which is right and up-going by $t''-t+1$ rows.   Moreover

     \

     $(i)$. The  lines with label $1$ and the composition of a neutral line with a line having label  $\ast$ joining boxes in $C_r(\infty),C_{r+1}(\infty)$ cannot have both the same start and end points.  In other words the label on a line (if it has one) is determined by its end-points.

     \

     $(ii)$. Distinct lines with label $1$ can neither have the same start-point in $C_r(\infty)$ nor the same end-point in $C_{r+1}$.

     \

     % Recall \ref {3.2.4} and set $\height C_{r+1}=t$. Consider a neutral line from $b:=R_{t'+j}\cap C_r(\infty)$ to $R_{t+1+j}\cap C_{r+1}(\infty)$ with either $m:=t+1- t'>1$ and $j=0$, or $m=1$ and $j \in \mathbb N$, composed with a line labelled by a $\ast$ from the latter box to a box in $C_{r+1}^{>t-m}$.
%
%      \
%
%      $(iii)$. The entries of $C_r(\infty)$ lying strictly above $b$ are either joined by a line to a box in $C_{r+1}$ with label $\ast$ to an entry of $C_{r+1}^{>t-m}$, or are stopped at $C_r$. In the latter case, the resulting line with label $1$ has right end-point in $C_{r+1}^{\leq(t-m)}$, so having entry strictly smaller than those of $C_{r+1}^{>t-m}$.

       Recall \ref {3.2.4} and set $\height C_{r+1}=t$. Consider a neutral line from $b:=R_{t'+j}\cap C_r(\infty)$ to $R_{t+1+j}\cap C_{r+1}(\infty)$ with either $m:=t+1- t'>1$ and $j=0$, or $m=1$ and $j \in \mathbb N$, composed with a line labelled by a $\ast$ from the latter box to a box in $C_{r+1}^{>t-m}$.

      \

      $(iii)$. The entries of $C_r(\infty)$ lying strictly above $b$ are either this composed line to $C_{r+1}^{>t-m}$, or are stopped at $C_r$. In the latter case, the resulting line with label $1$ has right end-point in $C_{r+1}^{\leq(t-m)}$, so having entry strictly smaller than those of $C_{r+1}^{>t-m}$.

      %\
%
%
%      Added to the hypothesis of $(iii)$, one may also have $r'\leq r$.  Again we may also have $r'<r$ and a horizontal neutral line from $C_{r'}$ to $C_{r'+1}$.
%      % having common entry $i$.
%
%      \
%
%      $(iv)$.  The entries of $C_{r'}(\infty)$ stopped at $C_{r'}$ define right going horizontal lines with label $1$ with entries strictly smaller than those of $C_{r+1}^{>t-m}$.

     \end {cor}

     \begin {proof} It remains to prove $(i)-(iii)$.

      $(i)$.  A line $\ell$ with label $1$ (resp. $(\ast)$) from $C_r(\infty)$ to $C_{r+1}$ arises when the entry an $i_m$ is stopped (resp. not stopped) at $C_r$.  Hence $(i)$.

     % \textbf {N.B.}  However one may have right going lines from a given box,  one having label $1$ and one being a neutral line composed with a line with label $\ast$.  For example take the composition $(1,1,2)$.  Then $1$ is not stopped at $C_1$ and we obtain a composed line $\ell_{1,2}$ with label $\ast$, whilst it is stopped at $C_2$ and we obtain a line $\ell_{1,4}$ with label $1$.  Again from Figure $1$ one sees that one have left going lines from a given box one with label $1$ and one with label $\ast$.  Again from Figure $1$ one sees that there may be several left (resp. right) going lines with label $\ast$ from a given box.

      $(ii)$ is by construction. Indeed if $i_m\in C_{r}(\infty)$ is stopped at $C_r$, then its right end-point is uniquely specified and has an entry not already being the right end-point of a further line with label $1$.

      $(iii)$.  Of course $(iii)$ is a mouthful. Hopefully it can be easily read off from Figure $1$ and its caption, following its notation using $m_1,m_2,m_3$.  In this $t=\height C_{r+1}=m_1+m_2$.

       Suppose $m>1$.  In terms of Figure $1$, we must take $m=m_1$ and $b:=R_{m_1+1}\cap C_r(\infty)$.  Then the entries strictly above $b$ lie in $C_r(\infty)^{\leq m_1}$.  They are stopped at $C_r$, giving horizontal lines with label $1$ having entries in $C_{r+1}^{\leq m_1}$, whilst $t-m=m_1+m_2-m_2=m_1$.  Yet the entries of $C_{r+1}$ increase down the rows and so an entry of $C_{r+1}^{\leq m_1}$ is strictly less than one of $C_{r+1}^{> t-m}$.

       Suppose $m=1$.  In terms of Figure $1$, we must either take $b \in R_j \cap C_r(\infty):j \in [m_1+2, m_1+m_2]$ and then its entries are stopped at $C_r$ and the resulting line with label $1$ has right end-point strictly above the largest entry in $C_{r+1}$, namely $j_{m_1+m_2}$, or we must take $b \in R_j \cap C_r(\infty):j \in [m_1+m_2+1, m_1+m_2+m_3-1]$ and the corresponding entries are joined to the entry $j_{m_1+m_2}$ of $C_{r+1}$ by a line labelled by a $\ast$.

       Finally the neutral horizontal lines lie in $R_{t''}:t''>m_1+m_2+m_3$. Thus their common end-points lie strictly below the right hand end-point of any line with label $1$, as can be easily read off from Figure $1$.
     \end {proof}

     \textbf{Remarks}.   Added to $ (i),(ii)$, one can have right going lines from a box,  one having label $1$ and one being a neutral line composed with a line with label $\ast$.  For example take the composition $(1,1,2)$.  Then $1$ is not stopped at $C_1$ and we obtain a composed line $\ell_{1,2}$ with label $\ast$, whilst it is stopped at $C_2$ and we obtain a line $\ell_{1,4}$ with label $1$.  Again (see Figure $1$) we can have left going lines from a box, one with label $1$, one with label $\ast$ and there may be several left (resp. right) going lines with label $\ast$ from a  box.

     $(iii)$  is used to prove Proposition \ref {6.2.4}.

     \

     \textbf{Overview.}   Let us spell out in detail the family of lines between
$C_r(\infty), C_{r+1}(\infty)$, described by the Corollary.  They are illustrated in Figure $1$.  In this the vertical lines are rounded for clarity.

     \

      First there can be $m_1$ horizontal lines with label $1$ between the first $m_1$ rows of $C_r,C_{r+1}$.

      Then there can be $m_2$ vertical lines with label $\ast$ from the \textit{first new} entry $i$ of $C_{r+1}(\infty)\setminus C_{r+1}$ to each of the last $m_2$ entries of $C_{r+1}$.  In this there can be up to $m_2$ lines with label $1$ from $C_r(\infty)$ going up by one row to $C_{r+1}$ to its last $m_2$ entries.

      After that there can be $m_3$ vertical lines with label $\ast$ from the $m_3$ new entries of $C_{r+1}(\infty)$ to the lowest entry of $C_{r+1}$.  In this there can be at most one line with label $1$ from $C_r(\infty)$ to the lowest entry of $C_{r+1}$ going up by $m_3$ rows.

      All non-negative integer values of $m_1,m_2,m_3$ are permitted.

      The result is illustrated in Figure $1$.  It has a striking up-down symmetry.

  \subsubsection {Translation to lines in $\mathscr T$}\label {3.2.7}

  A labelled line in $\mathscr T(\infty)$ joins boxes with entries in $[1,n]$.  If it is labelled by a $\ast$, its upper end point lies in $\mathscr T$ but its lower end-point does not lie in $\mathscr T$.  If it is labelled by a $1$ then its right hand end-point lies in $\mathscr T$, whilst its left-hand end point may or may not lie in $\mathscr T$. However there are unique boxes in $\mathscr T$ with any given entry and such a line determines a unique line in $\mathscr T$, with the same entries as end-points.

  Notice that $\mathscr T$ differs as a tableau from any $\mathscr T^\mathcal C$ \textit{only} through the lines and their labelling joining certain entries through the rules specified above and the choice of $\mathcal C$,

  \

  \textbf{Notation.}  The numerical data $\mathcal C$ which defines the limiting tableau $\mathscr T^\mathcal C(\infty)$ is thus transferred to a tableau $\mathscr T^{\mathcal C}$ with lines joining boxes in distinct columns and labelled by $1$ or by $\ast$.  The $\mathscr T^\mathcal C$ is called the component tableau defined by the numerical data $\mathcal C$.  We denote it by $\mathscr T^{i,j, \ldots}$ when $i,j,\ldots$ is a sequential choice of positive integers in batches so defining $\mathcal C$.   We recommend some practice in computing $\mathscr T^\mathcal C$.  In this the examples at the end of the main text may be helpful,

  \

  \textbf{N.B.}  It is of course no longer true that the lines with label $\ast$ are vertical in $\mathscr T^{\mathcal C}$.  Rather they are right and up or down going (or horizontal).  Again as in \cite [Lemma 6.2(i)]{FJ4} there can be several right going lines labelled by a $\ast$ from a given box of $\mathscr T^\mathcal C$.  However there is a novelty here.  Unlike \cite [Lemma 6.2(ii)]{FJ4}  there can be several left going lines labelled by $\ast$ from a given box of $\mathscr T^\mathcal C$.  All this is already illustrated by two possible choices of $\mathcal C$ for the composition $(2,1,1,2)$.

  On the other hand it is still true (Corollary \ref {3.2.6}(ii)) that a box in $\mathscr T^\mathcal C$ can have at most one left going line and at most one right going line labelled by a $1$.

  Again in $\mathscr T^\mathcal C$ the right-going lines with label $1$ need not only be up-going.  This is because by the rules of \ref {3.2.2} the above entry of $R_s\cap C_j(\infty)$ might have come from a column to its left in row $R_{s'}:s'\leq s$.  Consequently the right going lines with label $1$ being down-going and then up-going by $t-t'+1$ (with $t\geq t'$ defined in $(ii)-(iii)$ of \ref {3.2.4}. For the composition tableau one always has $t'=t$ (if $t$ exists) so the right going lines with label $1$ are down-going or up-going by at most one row (as we concluded in \cite [Lemma 5.3.2] {FJ4}).

  One might add that the rather complicated rules in \ref {3.2.2} introduce some symmetry in the set $\mathscr T^\mathcal C$ and ultimately some symmetry in the set of components constructed through component tableaux, whilst the canonical component tableau had a built-in asymmetry \cite [4.2.2]{FJ1}. In particular for the composition $(2,1,1,2)$ we obtained just one component of $\mathscr N$ using only the canonical component tableau; but trivially there is a second component obtained though the Dynkin diagram automorphism.

  Directly recovering this ``symmetry'' is not too easy in general.

  \

  For the composition tableau described in \cite [4.7]{FJ4}, the transition from $\mathscr T^\mathcal C(\infty)$ to $\mathscr T^\mathcal C$ is illustrated by going from Figure 2 to Figure 1 in \cite {FJ4}.

\section {The Excluded Roots and $B$ Saturation Sets} \label {4}

It is convenient to adopt the convention of \textbf{N.B.} of \ref {3.2.4} conflating roots, root vectors and lines in $\mathscr T^\mathcal C$.

 \subsection {General Remarks}\label {4.1}

 Our aim is to define for each component tableau $\mathscr T^\mathcal C$ a subalgebra $\mathfrak u^\mathcal C$ (or simply, $\mathfrak u$) spanned by root vectors of $\mathfrak m$ such that the component corresponding to $\mathscr T^\mathcal C$ is $\overline{B.\mathfrak u^\mathcal C}$.   The remaining positive roots of $\mathfrak m$ are called the excluded roots.   They will need to satisfy a number of properties set out in \ref {4.1.1}.  In this the component tableau is to specify a sum $e_\mathcal C$ (or simply, $e$) of root vectors given by the lines carrying a $1$ (which are not to be excluded roots \ref {4.1.1}(ii)) and a vector space $V_\mathcal C$ (or simply, $V$) spanned by the root vectors given by the lines carrying a $\ast$, which are to be amongst the excluded roots \ref {4.1.1}(i).  Again \ref {4.1.1}(iii) is required for $P$ stability.
  %whilst \ref {4.1.1}(iv) has no particular motivation but is ``halfway'' to the component being Lagrangian, which would result if the sum of the complementary root vectors formed a subalgebra, a well-known result sometimes attributed to the second-named author.

 However to ensure that a component tableau gives a component of $\mathscr N$ they will also have to satisfy a minimality property and \textit{this is far more delicate}.  For example we could have chosen the excluded roots to be all those not labelled by a $1$.  However this fails to give correct dimension (see \ref {4.3.1}) because $P.e$ need not be dense \cite [6.10.7]{FJ2} in a component of $\mathscr N$.

  \subsubsection {Initially Required Properties}\label {4.1.1}

  \

  \

 $(i)$. A root corresponding to a line given by a line carrying a $\ast$, is an excluded root.  This implies that $V \cap \mathfrak u=0$.

 \

 $(ii)$.  A root corresponding to a line given by a line carrying a $1$ is not an excluded root.  This implies that $e \in \mathfrak u$.

 \

 $(iii)$.  $\mathfrak u$ is stable under the action of $L^-$.  This implies that $\overline{B.\mathfrak u}$ is $P$ stable.

 \

 $(iv)$.  $\mathfrak u$ is a subalgebra of $\mathfrak m$.

 \

 Of course in $(iv)$ one can ask what class of subalgebras?

 %For orbital varieties we obtain subalgebras of $\mathfrak n$ complemented\footnote{The French would say supplemented, which is perhaps better terminology.} by a subalgebra of $\mathfrak n$.

 In exceptional cases the span of the excluded root vectors also form a subalgebra of $\mathfrak n$.  This is a sufficient condition for $\overline {B.\mathfrak u}$  to be an orbital variety closure.  It is not a necessary condition because  $\overline {B.\mathfrak u}$ does not determine $\mathfrak u$.  In general $\dim B.\mathfrak u \leq \frac{1}{2} \dim G.\mathfrak u$.  The necessary and sufficient condition for $\overline{ B.\mathfrak u}$ to be an orbital variety closure is that equality holds.

 \

Ultimately our procedure for constructing (algorithmically) a composition tableau $\mathscr T^\mathcal C$ gives rise to a subalgebra $\mathfrak u^\mathcal C$ of $\mathfrak m$ such that $\mathscr C := \overline{B.\mathfrak u^\mathcal C}$ is a component of $\mathscr N$. We have absolutely no idea of how else a component could otherwise have be written down, or even how to have guessed a suitable choice of $\mathfrak u^\mathcal C $.  At present even $\mathfrak u^\mathcal C$ is given only algorithmically and some might say by an arduous procedure.

  \subsubsection {Labelling.}\label {4.1.2}

  We use the same convention as in \cite [6.9.3]{FJ2}, namely in \textbf{M} we encircle the excluded roots in $\mathfrak m$.

  %(Here all the positive roots in Levi factor are excluded (see $(**)$ of \ref {4.1.4}) so encircling them would have been superfluous.)

  Again we recall our convention \cite [4.1]{FJ2} that the roots in $\mathfrak m$ corresponding to lines in $\mathscr T$ with a label $1$ (resp. $\ast$) are given the same labels in \textbf{M}.

  Then (cf \cite [6.9.3]{FJ2}) (i),(ii) of \ref {4.1.1} become,  in the special case of the canonical component (alias: composition) tableau, the following.

  \

  (i) Every $\ast$ is encircled.

  \

  (ii) No $1$ is encircled.

  \

We prove these assertions here for all component tableaux by a \textit{different} method than that used in \cite [6.9.3]{FJ2}, for the special case of the composition tableau.  Notably in the special case of the composition tableau this will give much less excluded roots than before.  In a weak sense (see \ref {4.3.11}) the present choice is a minimal choice.

 \subsubsection {Constructing the set of excluded roots.}\label {4.1.3}

 Our recipe for constructing the excluded roots is close to the simplest which recovers $(i)$, modified to obtain vanishing \ref {4.3}.  In this we shall follow closely but not exactly the construction and notations of \cite [Sect. 2]{FJ2}.

The excluded roots are constructed from the lines labelled by $\ast$ described in \ref {3.2.5}.  Here recall the labelling of the columns of $\mathscr D$ defined in \ref {2.1.1} and consider a column $C_{r+1}$ strictly to the right of $C_1$. Set $t=\height C_{r+1}$.  Recall that in $\mathscr T(\infty)$, a line $\ell$ labelled by a $\ast$ is vertical and may join the $m\geq 1$ entries $j^m_{r+1}:=\{j-k+1\}_{k=1}^m$ of the partial column $C_{r+1}^{>t-m}$, to the box $R_{t+1} \cap C_{r+1}(\infty)$ with entry $i_1$. Moreover this may be repeated finitely many times, yet with only $m=1$, joining the entries $i_u:u\in [1,v]$ of $R_{t+u}\cap C_{r+1}(\infty)$  to the  entry $j$ of $R_t\cap C_{r+1}$, as illustrated in Figure $1$.

 \

 \textbf{Notation}.  $j^m_{r+1}$ viewed often as the partial column $C^{>(t-m)}_{r+1}$ is abbreviated as \textbf{j}, or simply $j$ if it is a singleton.

 \

   Yet in the notation of \ref {2.2.2}, $i_u$ also appears in a unique column $C_h:=C_{(i_u)}$ of $\mathscr T$ with $h<r+1$ and in a unique row $R_f:=R_{(i_u)}$. Thus in particular $i_u<j$, for all $j \in \textbf{j}$.

   In what follows we generally drop the subscript on $i$, using also the simplified notation of the previous paragraph for the unique box $R_f\cap C_h$ of $\mathscr T$ containing $i$.

   \

   \textbf{N.B. 1}.  Since $\ell_{i,j}$ is labelled by a $\ast$, it follows by \ref {3.2.3} that $i$ appears in $C_{r+1}(\infty)$ in a row \textit{strictly} lower than $R_f$.  That is to say that $i$ goes non-trivially downwards as it progresses in $\mathscr T(\infty)$ through the rules of \ref {3.2.2}.

 \

 \textbf{Claim}. $\height C_g=f$, for some (unique) largest $g\in [1,h]$.
Moreover it admits a right neighbouring column $C'_g$ necessarily to the right of $C_h$.

 \

 %Indeed with respect to the construction {\ref {3.2.1}, \ref {3.2.2} of $\mathscr

  Recall the construction {\ref {3.2.2}, \ref {3.2.3} of $\mathscr T(t+1)$ and \textbf{N.B.1} above.

  For the entry $i\in C_h\cap R_f$ to appear in $C_{r+1}(\infty)$ on going rightwards from $C_h$ to $C_{r+1}(\infty)$ it must have gone in some first step strictly from $R_f$ to $R_{f'}$ into a column $C_{h'}(\infty)$ to the right of $C_h$ and to the left (not necessary strictly) of $C_{r+1}(\infty)$.

  By \ref {3.2.2}, Rules $(1),(2)$, this is only possible if the pair $C_h,C_{h'}$ is surrounded by a pair of neighbouring columns  of height $f$.

  \

  The left-hand column of this pair (which cannot be $C_{h}$) is the required column $C_g$ of height $f$.

  \

  This proves the claim.

  \

  %Let $C_g'$ be the right hand column of the above pair.

  \

  By construction $C_g$ lies to the left of $C_h$. It is  $C_h$ itself if perchance $\height C_h=f$.  Otherwise let  $C_{g_s},C_{g_{s-1}},\dots,C_{g_1}:h=:g_s>g_{s-1}>\dots\geq g_1$ be the columns of height $>f$ between $C_h,C_g$.  Set $g=g_0$.

\

  % Given a column $C_k$ of $\mathscr T$ and $p \in \mathbb N^+$, let $C_k^{>p}$, be the partial column formed from the rows $R_{p'}:p'>p$ of $C_k$.  Let $C_k^{\leq p}$ be the column obtained from $C_k$ by removing $C_k^{>p}$.
Recall the notation of \ref {3.2.6} for partial and truncated columns.

\

 Define \textbf{j} above. Recall that $\textbf{j}$ forms the entries of the partial column $C_{r+1}^{>(t-m)}$ and that $i$ is joined to each entry of $\textbf{j}$ by a vertical line labelled by a $\ast$.  (See Figure $1$).

  Let $\mathscr T_{(i,\textbf{j})}$ be the tableau obtained from $\mathscr T$ by moving $C_{r+1}^{>(t-m)}$ directly below $i$ in $C_h$ and if $ C^{>f}_h$ is not empty, successively making $C_{g_k}^{>f}:k=s,s-1,\ldots,1$ replace $C_{g_{k-1}}^{>f}$, in this skipping over the columns of height $<f$.

  \

   \textbf{Definition.} We refer to the former (resp. latter) operation as the \textbf{primary} (resp. \textbf{secondary}) shifting of columns.  The latter is absent if $\height C_h=f$, equivalently if $i$ is the lowest entry of $C_h$.

  \

 %\textit{Notice since the only columns of height $f$ in $[C_g,C_f[$ is $C_g$, this means skipping over the columns of height  $\leq f$}.

   $(*)$. Since  $\height C_g=t$ and is the first with this property, this process ends exactly  when $C_g$ is reached.

   \

  \textbf{N.B. 2}. Thus the columns which are altered are $C_{r+1}$ and those columns which  lie in $[C_g,C_h]$ of height $\geq f$.

\

   This construction is a refinement of that used in \cite [2.6]{FJ2}.  In particular an example of this  manipulation can be viewed in \cite [Figures 1,2]{FJ2}.

   \

\

Some readers may find this complicated, but it is deeply and mysteriously intertwined with the construction of a component of $\mathscr N$.

 \subsubsection {The Excluded Roots.}\label {4.1.4}

 %Recall the notation of \ref {4.1.3}. Consider the left movement of partial columns, described in \ref {4.1.3}.

 One has the following standard result (cf. \cite [Lemma 2.5]{FJ2}).

 \

 $(*)$.  For a given $w \in W$ the root subspaces of $\mathfrak m$ which do \textit{not} lie in $\mathfrak n \cap w \mathfrak n$ are the $x_{i,j}:i<j$ with $i$ after $j$ in the word form of $w$.

\

 Recall the notation of \ref {4.1.3}.

\

 Let $w_\mathscr T$ (resp. $w_{\mathscr T_{i,\textbf{j}}}$ or simply $w_{i,\textbf{j}}$) be the Weyl group element obtained from $\mathscr T$ (resp. $\mathscr T_{i,\textbf{j}}$) with word form given as in \cite [2.5]{FJ2},  that is to say we define its word form by reading the columns of $\mathscr T$ (resp. $\mathscr T_{(i,\textbf{j})}$) from bottom to top, starting from the leftmost column and moving rightwards.

 Let $\mathfrak l^+$ (resp. $\mathfrak l^-$) be the span of the positive (resp. negative) root vectors in the Levi factor $\mathfrak l$ of $\mathfrak p$.

 Set $\mathfrak u_{i,\textbf{j}}:= \mathfrak n \cap w_{i,\textbf{j}} \mathfrak n$. It is a subalgebra of $\mathfrak n$ complemented by a subalgebra of $\mathfrak n$ and so $\overline {B.\mathfrak u_{i,\textbf{j}}}$ is an orbital variety closure.

 \

   \textbf{Definition}.  The set of positive root vectors $X_{i,\textbf{j}}$  of $\mathfrak m$ \textit{not} appearing in $\mathfrak u_{(i,\textbf{j})}$, is called the set of excluded root vectors defined by the pair $(i,\textbf{j})$. (Being in $\mathfrak m$ means that excluded roots lying in the Levi factor are ignored.)  Let $\textbf{X}_{i,\textbf{j}}$ denote the span of the excluded root vectors defined by $X_{i,\textbf{j}}$.

 \

 \begin {lemma}

 \

 $(i)$.  $ (\mathfrak n \cap w_{\mathscr T}\mathfrak n) \cap \mathfrak l^+=0$.

 \

 $(ii)$.  $\mathfrak u_{(i,\textbf{j})} \cap \mathfrak l^+=0$.

  \

 $(iii)$.  $\mathfrak u_{(i,\textbf{j})}$ is stable under the action of $\mathfrak l^-$.

 \end {lemma}

 \begin {proof}

 If $k,l:k<l$ are in a fixed column of $\mathscr T$, then $k$ lies above $l$ and hence after $l$ in the word from of $w_\mathscr T$, through the definition of the latter.  Thus $(i)$ follows from $(*)$.

 \

  In the construction of $\mathfrak u_{(i,\textbf{j})}$ entire lower parts of columns are moved to the left. Then comparing $w_\mathscr T$ and $w_{(i,\textbf{j})}$, one immediately deduces that the property required for $(i)$ is not upset.  Hence $(ii)$.

  \

  $(iii)$ follows from $(ii)$.  Indeed fix a simple root $\alpha$ of the Levi factor $\mathfrak l$.  Then the roots of $\mathfrak m$ lie in $\alpha$-strings, that is to say form simple $\mathfrak {sl}(2)$ modules with one dimensional weight spaces with weights which are a fixed translate by integers multiples of $\alpha$. By $(ii)$, the intersection of  an $\alpha$-string with $\textbf{X}_{i,\textbf{j}}$ is stable under $x_\alpha$ and therefore its complement is stable under $x_{-\alpha}$, in other words $\mathfrak u_{(i,\textbf{j})}$ is stable under $x_{-\alpha}$.  Since these root vectors generate $\mathfrak l^-$, this proves $(iii)$.

  \end {proof}

%  \textbf{Definition}.  The set of positive roots $X_{i,j}$ \textit{not} appearing in $\mathfrak u_{(i,j)}$, is called the set of excluded roots defined by the pair $(i,j)$.

  \

  \textbf{Definition}. Recall \ref {4.1.3}. It is convenient to distinguish between the excluded roots $X_{i,\textbf {j}}$ called \textbf{primary} which arise from the primary shifting of columns and those called \textbf{secondary} which arise from the secondary shifting of columns.

  %$C_{r+1}^{>(t-m)}$ under $i$, which we call ``primary'' and the rest (which arise from the shifting of the remaining lower parts of columns to the left) which we call ``secondary''.

  \

  \textbf{Observations}.

  \

  $(1)$ The second co-ordinate $j'$ of a primary  excluded root $x_{i',j'}$ of  $X_{i,\textbf{j}}.$ lies in $C_{r+1}^{>(t-m)}$,  The first co-ordinate $i'$ is an entry in $C_{(i)}$ strictly above $i$ or lying in a column strictly to the right of $C_{(i)}$ and to the left of $C_{r+1}$, subject to the condition $i'<j'$.

  %$(1)$.The first co-ordinate $i'$ of a primary  excluded root $x_{i',j'}$ of  $\mathfrak u_{(i,\textbf{j})}$ is an entry in $C_{(i)}$ strictly above $i$ or lying in a column strictly to the right of $C_{(i)}$ and to the left of $C_{r+1}$.
%  The second co-ordinate $j'$ $.  All this subject to the condition $i'<j'$.

  \

  $(2)$. There are no secondary excluded roots in $X_{i,\textbf {j}}$, if $i$ is the unique lowest entry in $C_{(i)}$.

  \

  $(3)$.  In the notation of \ref {4.1.3}, the second co-ordinate $j'$ of a secondary excluded root $x_{i',j'}$ lies in $C^{>f}_{g_{k}}$ for some $k\in [1,s]$ and then
  its first co-ordinate $i'$ lies in $[C_{g_{k-1}}^{\leq f},C_{g_{k}}[$.
  Then those with first co-ordinate $i'\in C_{g_k}$ lie in the Levi factor and can be ignored.  Since the columns in $]C_{g_{k-1}},C_{g_{k}}[$ are of height $<f$, we obtain  $i' \in R^f\cap [C_{g_{k-1}},C_{g_{k}}[$. In particular the row containing $j'$ is strictly below that containing $i'$.

  \subsubsection {The Definition of $\mathfrak u^\mathcal C$.}\label {4.1.5}

 Fix a component tableau $\mathscr T^\mathcal C$. Let $\textbf{L}$ denote the set of all pairs $(i,\textbf{j})$ defined in \ref {4.1.3}.

Set $\mathfrak u^\mathcal C:=\cap_{i,\textbf{j}\in \textbf{L}}\mathfrak u_{(i,\textbf{j})}$.  Of course it is again a subalgebra of $\mathscr n$ \textit{but not } in general complemented by a subalgebra in $\mathfrak n$.

Thus \ref {4.1.1}(iv) holds.  Again by Lemma \ref {4.1.4}$(iii)$ above, we obtain  \ref {4.1.1}(iii).  It remains to prove \ref {4.1.1} $(i),(ii)$.  They are re-expressed as \ref {4.1.2} $(i),(ii)$ of \ref {4.1.2}, and established below.

 Thus $\overline{B.\mathfrak u^{\mathcal C}}$ is not in general an orbital variety closure, though it can be even if $\mathfrak u^\mathcal C$ is not complemented by a subalgebra in $\mathfrak n$.  Indeed it is a closed subvariety of $\overline{G.\mathfrak u_{\mathscr T^{\mathcal C}}}$ which itself admits a unique dense nilpotent orbit $\mathscr O$. By Spaltenstein \cite {S} for any nilpotent orbit $\mathscr O$ the intersection $\mathscr O\cap \mathfrak n$ is equi-dimensional of dimension $\frac{1}{2} \dim \mathscr O$, which is a union of orbital varieties. By definition all orbital varieties so obtain.

 Thus $B.\mathfrak u^{\mathcal C} \subset G.\mathfrak u^{\mathcal C} \cap \mathfrak n$ and so has dimension $\leq \frac{1}{2} \dim \mathscr O$ with equality if and only if $\overline{B.\mathfrak u^{\mathcal C}}$ is an orbital variety closure.

 \

 \textbf{Definition}. The excluded root vectors $X$ defining the complement of  $\mathfrak u^\mathcal C$ in $\mathfrak m$ is the union of the excluded root vectors $X_{i,\textbf{j}}: (i,\textbf{j}) \in \textbf{L}$.  Let \textbf{X} denote the space they span.

 \

 The excluded root vectors of $\mathfrak u^{\mathcal C}$ will be enclosed in a circle, particularly when they are presented in \textbf{M}.  See Figure $12$.

 \

 %It is highly advisable to compute some examples of the labelling of the root vectors in \textbf{M} by $1$, $\ast$ and their encirclement.  As an example we give the result in Figure $12$ for the five different component tableaux which result from the parabolic defined by the composition $(2,1,2,1,2,1)$.  The reader may also wish to show that the all the six Benlolo-Sanderson invariants vanish on converting the circles to zeros, for each of the five component tableaux - an easy computation but not so easy when $n$ is very large (see Section \ref {4.3}).
%
% \

 $(*)$. Recall  in $\mathscr T(\infty)$, that $i$ is joined to each entry of $\textbf{j}$ by a vertical line labelled by a $\ast$, and that all lines in $\mathscr T(\infty)$ labelled by a $\ast$ occur as $(i,\textbf{j})$ runs over $\textbf{L}$.

 \

 \begin {lemma} Every $\ast$ is enclosed in a circle.
 \end {lemma}

 \begin {proof} Suppose $(i,\textbf{j}) \in \textbf{L}$.  For all $j \in \textbf{j}$, one has $i<j$ and by the construction $\textbf{j}$ appears directly below $i$ in $\mathscr T_{i,\textbf{j}}$, consequently $i$ is after $j$ in the word form of $w_{i,j}$. Thus the assertion follows from $(*)$ and the definition of the excluded roots of $\mathfrak u^\mathcal C$.
 \end {proof}

%\textbf{Remark}. Of course this was practically by construction.

\subsection {Verifying $(ii)$.}\label {4.2}

This is much more difficult and we need some preliminaries.

\subsubsection {Non-crossing}\label {4.2.1}

Recall that in $\mathscr T^{\mathcal C}(\infty)$, there can be multiple entries of a given $i\in [1,n]$ and at least one. Moreover the entries occur in adjacent columns either in the same row or on passing from left to right go down by one or possibly several rows (see \ref {3.2.2}).

\

\textbf{Definition}. The set of entries of $i$ in $\mathscr T(\infty)$ is called an $i$-string and denoted by  $S(i)$.

\

 Notice by the \textbf{N.B.1} of \ref {3.2.2} an $i$-string can pass through at most one box of a column $C(\infty)$ of $\mathscr T^{\mathcal C}(\infty)$ (being then the unique box in $C(\infty)$ with entry $i$).

\begin {lemma}  (Notation \ref {2.1.1})

\

$(i)$. For each $i\in [1,n]$, there exists a unique $i$-string. It is determined by (and determines) the unique entry of $i\in \mathscr T$ and the label on the lines $\ell_{i,i'}$ with right end-point in $C_{r+1}:r \in [1,k[$.

The lines $\ell_{i,i'}$ in $\mathscr T^\mathcal C$ are labelled by a $\ast$ or by a $1$.

\

$(ii)$. For $i_1,i_2 \in [1,n]$ distinct, the strings $S(i_1),S(i_2)$ cannot cross.
\end {lemma}

\begin {proof}   %To simply matters we shall use the Vav Conversive (\ref {3.1.2}.

$(i)$.  Recall that $\mathscr T$ (defined in \ref {2.2.2}) is a sub-tableau of $\mathscr T^\mathcal C(\infty)$.  Then $S(i)$ starts at the unique box of $\mathscr T$ containing $i$.  Then by \ref {3.2.2} either $i$ moves rightwards from a column $C_r(\infty)$ to a right adjacent column of $C_{r+1}(\infty)$ of $\mathscr T^\mathcal C(\infty)$ and downwards by $m\geq 1$ rows by a neutral line, or $i$ is stopped at $C_r$.

In the first case by Corollary \ref {3.2.6} (see Figure $1$) it is joined by exactly $m$ up-going vertical lines all having a label $\ast$ to (successive and distinct) entries of $C_{r+1}$, so determining $m$.  This may be repeated for larger values of $r\in [1,k[$.

Again $i$ may be stopped at $C_r$ and then  it is joined by a line with label $1$ to an entry of $C_{r+1}$.

If $i$ is not stopped for any $r\in [1,k[$ then $S(i)$ reaches the last column of $\mathscr D$.

$(ii)$.  Let $C_r,C_{r+1}: r \in [1,k[$ be adjacent columns, where a crossing first occurs in $\mathscr T^\mathcal C(\infty)$. We can assume that on reaching $C_r(\infty)$ that $S(i_2)$ lies below $S(i_1)$. In particular $i_1\in R_{r_1}\cap C_r(\infty), i_2 \in R_{r_2}\cap C_r(\infty)$ with $r_2>r_1$. Then the assertion can be read off from Figure $1$ which shows that none of the neutral lines from $C_r(\infty)$ to $C_{r+1}(\infty)$ cross.

\end {proof}
%
%\textbf {Remark}.  Recall Example $5$ illustrated in Figure $6$. In this $7$ is not stopped at $C_4$ but goes down two rows entering $R_3$, whilst by contrast $8$ is stopped at $C_4$. On the other hand the $5$-string starts at $R_2$ and can go below $7$ into $R_4\cap C_5$.  However then the  pair $C_1,C_7$ of neighbouring columns of height $3$ are no longer free so preventing $7$ entering $R_4 \cap C_7$.

\subsubsection {Starting Places}\label {4.2.2}

Recall that $\mathscr T$ is a sub-tableau of any $\mathscr T^{\mathcal C}(\infty)$ constructed from $\mathscr T$.

The starting place of an $s$-string is defined to be the unique box of $\mathscr T^{\mathcal C}(\infty)$ in which $s$ appears in $\mathscr T$. It is the leftmost entry of $s \in \mathscr T^{\mathcal C}(\infty)$.

Suppose that $S(i_1),S(i_2)$ pass through a common column $C'(\infty)$. Then by Lemma \ref {4.2.1} it makes sense to say that $S(i_2)$ lies strictly below $S(i_1)$ in $\mathscr T^{\mathcal C}(\infty)$.

\begin {lemma} Under the above hypothesis, the starting place of $S(i_2)$ is to the left of $S(i_1)$.
\end {lemma}

\begin {proof}

In the notation of \ref {2.2.2},  $C_{(i_1)}$ (resp. $C_{(i_2)}$) is the column of $\mathscr T$ containing the starting place of  $i_1$ (resp. $i_2$). The lemma amounts to saying that $C_{(i_2)}$ lies to the left of $C_{(i_1)}$.

%Here we can assume that there is no string $S(i)$ lying strictly between $S(i_1)$ and $S(i_2)$, which we call minimality.

Let $C_{r+1}(\infty)$ be the leftmost column through which $S(i_1),S(i_2)$ both pass
and suppose  $C_{r+1}$ lies strictly to the right of $C_{(i_1)}$.

Let $C_r$ be the left adjacent column to $C_{r+1}$. Then $i_1 \in C_r(\infty)$ and since by hypothesis it enters $C_{r+1}(\infty)$, it must enter $C_{r+1}(\infty)\setminus C_{r+1}$, by \ref {3.2.2}.  Since by hypothesis $i_2$ lies in $C_{r+1}(\infty)$ strictly below $i_1$,  necessarily $i_2 \in C_{r+1}(\infty)\setminus C_{r+1}$.  Then by \ref {3.2.2} again, we must have had $i_2 \in C_r(\infty)$, contradicting the choice of $C_{r+1}$.

This contradiction forces $C_{r+1}=C_{(i_1)}$ and so $C_{(i_2)}$ lies to the left of $C_{(i_1)}$, as required.
\end {proof}

\subsubsection {Primary Excluded Roots}\label {4.2.3}

Recall the lines with label $\ast$  described in  \ref {3.2.5} and the lines with label $1$ described in Corollary \ref {3.2.6} with respect to the pair $C_r,C_{r+1}$ of $\mathscr T(\infty)$.  Below we use Figure $1$ whose notation we adopt.

\

\textbf{Observations}.

\

$(1)$. By Observation $(1)$ of \ref {4.1.4} the second co-ordinate $j'$ of a primary excluded root defined by the pair $C_r,C_{r+1}$ lies in $C_{r+1}$.

\

$(2)$. On the other hand if  $i'$ is an entry of $C_{r'}(\infty)$ stopped at $C_{r'}$, then by Corollary \ref {3.2.6}, the second co-ordinate $j'$ of the root vector $x_{i',j'}$ so defined lies in $C_{r'+1}$.

\

$(3)$. Observe that the entries of $C_{r+1}$ label the columns of the column block $\textbf{C}_{r+1}$.

\begin {lemma} Fix a component tableau $\mathscr T^\mathcal C(\infty)$. The roots vectors defined by the pairs $(i',j')$, where $i'$ is an entry of $C_{r'}(\infty)$ stopped at $C_{r'}$ is not excluded by any one of the primary excluded roots of  $\mathfrak u_{i,j}$ defined by a vertical line joining $i \in C_{r+1}(\infty)$ to $j \in \textbf {j} \subset C_{r+1}$.
\end {lemma}

\begin {proof} By $(1),(2)$ we can assume $r=r'$ without loss of generality.  Even more convincingly, if $r\neq r'$, then by $(3)$ these two sets of vectors lie in different column blocks of \textbf{M}.

The reminder of the proof may be visualized through Figure $1$ and this should make the proof transparent.

In Figure $1$, there are $m_1$ horizontal lines from $R_i\cap C_r(\infty)$ to $R_i\cap C_{r+1}$ with label $1$.  These cannot be primary excluded roots for the pair $(i,\textbf{j})$ since their right end-points are left in place in constructing $\mathscr T_{i,\textbf{j}}$.

\

In the first step following the appearance of the said $m_1$ horizontal lines, the entry $i_{m_1+1}$ of  $R_{m_1+1} \cap C_r(\infty) $ is joined by a neutral line to the entry of $R_{m_1+m_2+1} \cap C_{r+1}(\infty)$ which is then joined by vertical lines with label $\ast $ to the lowest $m_2$ entries of $C_{r+1}$, ending at its lowest point in row $R_{m_1+m_2}$.

Apart from the end-points of these horizontal lines, disposed of above, an entry $i'$ of $C_r(\infty)$ stopped at $C_r$ lies in up to $m_2$ places \textit{below} $i_{m_1+1}$ .

Set $i_{m_1+1}=i$.

On the other hand $i'$ lies in  a unique box of $\mathscr T$, which by Lemmas \ref {4.2.1}, \ref {4.2.2},  lies either, strictly below $i$ in the unique column $C:=C_{(i)}$ of $\mathscr T$ containing $i$, or in a column strictly to the left of $C$.  It is the first co-ordinate $i'$ of a root $x_{i',j'}$ labelled by a $1$.

Yet by Observation $(1)$, \ref {4.1.4}, the first co-ordinate $i'$ of an excluded root $x_{i'j'}$ of  $\mathfrak u_{(i,\textbf{j})}$ is an entry in $C$ strictly above $i$ or lying in a column strictly to the right of $C$.

Then by the last two paragraphs, a root with label $1$ coming from the first step and a primary excluded root cannot coincide.

%it cannot be a root labelled by a $1$ described in the paragraph above.

\

For the second step, set $t =m_1+m_2$ still using the notation of Figure $1$.

Then $i_k:k\in [t+1,t+m_3-1]$ is an entry of  $R_{k}\cap C_{r}(\infty)$ joined by a neutral line to $R_{k+1}\cap C_{j+1}(\infty)$ which is further joined by a vertical line with label $\ast$ to the unique lowest box in $R_t\cap C_{r+1}$. Then there is at most a line $\ell$ with label $1$ from $i'=i_{t+m_3}$ to the entry $j'$ of $R_{t}\cap C_{r+1}$.

The corresponding root defined by this possibly new line labelled by a $1$ is not excluded by the primary excluded roots from the new lines labelled by a $\ast$.  Indeed $i'=i_{t+m_3}$ lies strictly below $i=i_{t+k}:k\in [1,m_3-1]$, so we may apply Lemmas \ref {4.2.1}, \ref {4.2.2} and proceed just as in the first part.

Yet there is one more subtle point!   Why don't the roots with label $\ast$ coming from the first step exclude the root with label $1$ coming from the second step, since $i_{t+m_3}$ of the second part lies strictly below the  $i_{m_1+k}:k \in [1,m_2]$ of the first part?

Now in the first step, the resulting primary excluded roots lie in the rightmost $m_2$ columns of the column block $\textbf{C}_{r+1}$ in \textbf{M} \textit{and} in rows given by the entries below $i$ in $C_{(i)}$ or in a column to the right of $C_{(i)}$.

On the other hand in \textbf{M} the line $\ell$ coming from the second step (if it exists) lies in the last column of the column block $\textbf{C}_{r+1}$, but in a row defined by an entry of a column of $\mathscr T^\mathcal C$ strictly to the left of $C_{(i)}$.

\end {proof}

\subsubsection {Secondary Excluded Roots}\label {4.2.4}

\begin {lemma}  No element of $1$ can be encircled by a secondary excluded root.
\end {lemma}

\begin {proof}

Fix some component tableau $\mathscr T^\mathcal C$.

A secondary excluded root is created by a line $\ell_{i,j}$ labelled by a $\ast$.
Following the notation of \ref {4.1.3}, let $R_f\cap C_h$ be the unique box in $\mathscr T$
containing $i$.  As noted in Claim \ref {4.1.3} there is a column $C_{h'}$ strictly to the
right of $C_h$ in $\mathscr T$ such that $i$ enters $C_{h'}(\infty)$ strictly below $R_f$
and this forces there to be a pair of neighbouring columns $C_g,C_g'$ surrounding
$C_h,C_{h'}$ with $i \in \mathscr B^f_{C_g,C_g'}$.

By Observation $(3)$ of \ref {4.1.4}  the co-ordinates of a secondary excluded root
$x_{i',j'}$ satisfy $i' \in R^f\cap [C_g,C_h[, j' \in R_{f'}\cap C_{r+1}$ for some $f'>f$ and column $C_{r+1}\in ]C_g,C_h]$.

Now suppose $\ell_{i'',j'}$ is labelled by a $1$. It is enough to show that $i''$ cannot equal $i'$. Of course $x_{i'',j'}\in \mathfrak m$,
so $i''$ lies in a column strictly to the left of $C_{r+1}$.
Moreover $i'' \in \mathscr T^\mathcal C(\infty)$ must be stopped at $C_r$ and
must lie in $R_{f''}\cap C_r(\infty)$ with $f''\geq f'$.  This is because the line $\ell_{i'',j'}$ is up-going as it passes from $C_r(\infty)$ to $C_{r+1}$ (see Figure $1$).

We conclude that the $i''$-string descends strictly from $R^f$ to $R_{f''}$.

This means in particular that $i''$ belongs to a batch $\mathscr B^f_{C,C'}$ obviously to the left of $\mathscr B^f_{C_g,C_g'}$.  These batches must be distinct because the first (resp. second) is used to lower $i''$ (resp. $i$) whilst  \ref {5.2.2}$(1)$ allows only entry of each batch to be used for a given component tableau.

We conclude that $C'$ lies to the left of $C_{g}$, though not necessarily strictly.

Yet we saw above for $\ell_{i',j'}$ to be a secondary excluded root one must have $i' \in [C_g,C_h[$ and for $\ell_{i'',j'}$ to be labelled by a  $1$ one must have $i''\in [C,C'[$.  Yet by the previous paragraph these two sets have null intersection, so we cannot have $i''=i'$.

Hence the assertion of the lemma.

\end {proof}

\textbf{Example $6$.}   Consider composition $(1,2,2,1,3,2)$ and take $i=8,j=11$.  Then $f=2$ and $C_h=C_5$.  Consequently $C_{h'}=C_6$ and $C_g=C_3,C'_g=C_6$.

The excluded roots obtain from $\mathscr T_{8,11}$ drawn in Figure $8$.  Those which are secondary are $\alpha_{i',9}$ with $i'\in\{4,5,6\}$.  (Those which are primary are $\alpha_{i',11}$ with $i'\in \{7,8\}$).  Thus $f'=3,r=4$.

The second batch $\mathscr B^f_{C,C'}$ in the lemma is obtained by taking $C=C_2,C'=C_4$ and necessarily $f=2$.  This forces $i''=3$, so the resulting line $\ell_{3,9}$ with label $1$ cannot give a secondary excluded root coming from the line $\ell_{8,11}$ with label $\ast$.

\

\textbf{Remark.}  Some readers might complain that this is a rather ``degenerate'' example, since many columns of the lemma coincide and particularly the $i$-chain is reduced to a singleton. Yet the proof itself only depends on the box containing $i$, the subsequent peregrinations of the $i$-chain being unimportant for secondary excluded roots\footnote{Yet important for the primary excluded roots.}, and the rather simple secondary shifting of partial columns. What is important is condition \ref {4.1.3}$(*)$, that this shifting does not go beyond $C_g$, whilst the second author greedy for zeros to obtain ``vanishing'', wanted it to push further, but finally desisted.

%Consider the composition $(1,2,2,1.3,2)$ and set  $C=C_2,C'=C''=C_3,(C'')'=C_6$.
%
%Choose the component tableau in which $4$ is lowered below $6$ in the first step. Then $i''=3$ (resp. $8$) is lowered below $i'=5$ (resp. $11$) through the batch $\mathscr B^2_{2,3}$ (resp. $\mathscr B^2_{3,6}$). This means that we are taking $f=2,r=4,r'=5$ in the notation of the lemma.  The resulting tableau is illustrated in Figure $6$.
%
% Now see Figure $8$. Following \ref {4.1.3},  the sole entry $11$ of  $C_{r'+1}^{>f}$ is placed in $R_3\cap C_5$ below $8$ pushing $j'=j''=9$ leftwards and horizontally, skipping over $C_4$, which has only height $1<2=:s$ to the first available column of height $\geq f$ to the right of $C_3$.  This means that $f'=f''=3$ in the notation of the lemma.
%
%This makes $x_{i',j'}=x_{5,9}$ a secondary excluded root. Yet $9$ cannot end up below $3$ in $C_2$, by the rule in \ref {4.1.3} and so it cannot make $\alpha_{3,9}$ a secondary excluded root.  On the other hand it is $3$ rather than $5$ which lies in $\mathscr B^2_{2,3}$ giving a line $\ell_{i'',j''}=\ell_{3,9}$ with label $1$.

%Complicated but this is only a baby example!

\subsubsection {}\label {4.2.5}

Combining the above two lemmas we obtain.

\begin {prop}  In no component tableau is a $1$ encircled.
\end {prop}

\begin {proof} The construction of \ref {4.1.3} gives by definition \textit{all} the excluded roots.  These are encircled (in \textbf{M}).  Then Lemmas \ref {4.2.3}, \ref {4.2.4} can be interpreted as saying that no circle encloses a $1$.   This holds for every component tableau.  Hence the assertion.

\end {proof}

\textbf{Remark.}  This proof, though not so short, is easier (and more general) than that given in \cite [6.9.2, 6.9.3]{FJ2} for which the reader was advised to take a deep breath.

\subsection {A Vanishing Theorem.}\label {4.3}

\subsubsection {Meaning of the excluded roots.}\label {4.3.1}

Let $\mathscr T^\mathcal C$ be a component tableau.  Recall that $\mathcal C$ is numerical data specifying $\mathscr T^\mathcal C$ through selection of entries in batches.

Define $\mathfrak u^\mathcal C$ as in \ref {4.1.5}. Recall (\ref {4.1}) that $\mathfrak u^\mathcal C$ in complemented in $\mathfrak m$ by the direct sum of the excluded root vectors.

   As before set $\mathscr C:=\overline{B.\mathfrak u^\mathcal C}$.  It is an irreducible closed subvariety of $\mathfrak m$.  Let \textbf{g} be the number of generators of $S(\mathfrak m)^{p'}$.  We wish to show that $\mathscr C$ lies in $\mathscr N$, is of dimension $\dim \mathfrak m-\textbf{g}$ and  is a component of $\mathscr N$.

%As in \cite [2.6]{FJ2}, let $\mathfrak u^\mathscr C$ be the subspace of $\mathfrak n$

Recall that $e_\mathcal C$ (resp. $V_\mathcal C$) denotes the sum of the root vectors (resp. root subspaces) defined by the lines with label $1$ (resp. $\ast$). By Proposition \ref {4.2.5} one has $e_\mathcal C\in \mathfrak u^\mathcal C$.  On the other hand \ref {4.1.2}(i) holds by construction (cf \ref {4.1.3}) and so $V_\mathcal C \cap \mathfrak u^\mathcal C=0$.

 We may present the excluded roots with respect to a given $\mathfrak u^\mathcal C$ as being represented by a $0$ in $\mathfrak m$ rather than an $O$.
 %In the latter role it may encircle a $\ast$, but not a $1$ (Proposition \ref {4.2.5}).

 In the former role we may consider the generating Benlolo-Sanderson invariants being evaluated on $\mathfrak m$ by setting all the entries of the excluded roots equal to zero.  Since the Benlolo-Sanderson invariants are $B$ invariant and homogeneous their vanishing means that $\mathscr C \subset \mathscr N$.

%Our analysis will follow partly that of \cite [2.6]{FJ2}.  Yet there are two difficulties caused by our use of much less excluded roots and that the fact that the number of  component tableaux can become very large.

The proof of vanishing is given below.  We must do this for all choices of parabolics of $\mathfrak {sl}(n)$, all components $\mathscr C$ and all Benlolo-Sanderson invariants $I^s_{C,C'}$.

This is equivalent to proving that a triply infinite sequence of determinants growing exponentially in $n$, all vanish when the excluded roots are set equal to zero.  Considering this, the proof is rather easy!  Its proof is due to the first author and is also her inspired choice of excluded roots which are practically a minimal set.  A larger set of excluded roots is possible but there is a danger here that the resulting zero locus $\mathscr Z$  defined by this bigger set of excluded roots may be too small.  For example the second author had suggested that one might choose the set of all excluded roots to be those not labelled by a $1$. Through the existence of a Weierstrass section it is immediate for this choice that $\mathscr Z$ is contained in $\mathscr N$.  Yet by \cite [6.10.7]{FJ2}, the sum $e$ of the root vectors labelled by a $1$ need \textit{not} be regular, that is generate a $P$ orbit of codimension \textbf{g}). As in consequence $\mathscr Z$ has in general too small a dimension\footnote {The second author was heavily censored by his students for explaining proofs that \textit{do not} work.  Perhaps this will be true of our readers.  Nevertheless we do wish to emphasize the extreme delicacy of the present calculations.}.

We start with some preliminaries.

\subsubsection {Chains of Roots}\label{4.3.2}

Let $l_1<l_2<\ldots <l_k$ be elements of $[1,n]$ in $\mathscr D$.  It gives a chain of roots $\alpha_{l_i,l_{i+1}}:i \in [1,k-1]$ forming a system of type $A_{k-1}$. The corresponding monomial $\prod_{i=1}^{k-1}x_{l_i,l_{i+1}}$ is said to form a monomial chain (of type $A_{k-1}$).

The Benlolo-Sanderson invariants are polynomials in the $x^*_{i,j}$ which by definition are the co-ordinate functions on $x_{i,j}\in \mathfrak m$.  We call these the dual vectors.

A composite line $\mathscr L$ is a concatenation of lines $\ell_{i,j}$, with $i<j$,  joining boxes in $\mathscr T$ which we call constituents (or constituent lines) of $\mathscr L$. We say that a composite line $\mathscr L$ is directed if every constituent line passes strictly from left to right.  Then the entries of the boxes satisfy $l_1<l_2<\ldots <l_k$.  Thus a directed composite line gives a monomial chain of type $A$.

If two composite lines are disjoint (pass through no common boxes) then the roots in the different chains are orthogonal and the monomial chains are said to be disjoint.

Fix a pair of neighbouring columns $C,C'$ of height $s$.  By a slight and temporary modification of the notation of \ref {2.1.1} the columns in $[C,C']$ will be denoted by $C=:C_1,C_2,\ldots, C_k:=C'$.  Set $c_i=\height C_i$.

Define $d_{\mathscr D}$ as in \ref {2.3.2}.  It is exactly the number of boxes in $[C,C']$ strictly below $R_s$.

As noted in \cite [4.1.6]{FJ1}, after \cite {BS} the lowest degree coefficient of the restriction of $\textbf{M}^*_s$ to $\mathfrak m +a \Id$ starts as  $a^{d_{\mathscr D}}I^s_{C,C'}$.  Miraculously this is a $\mathfrak p'$ invariant.

\begin {lemma} The Benlolo-Sanderson invariant $I^s_{C,C'}$ is a sum (with appropriate signs) of products of $s$ monomial chains of dual vectors.  The products are given by all possible sets of $s$ disjoint directed composite lines starting from $C$ and ending in $C'$ such that every chain passes through all possible choices of $\min {(c_i,s)}$ boxes in $C_i:i\in ]1,k-1[$.

\end {lemma}

\begin {proof} The entries in the  columns strictly to the left of $C$ and strictly to the right of $C'$  are absent as both the first or last co-ordinate is a root vector occurring in  $I^s_{C,C'}$.  Thus we can assume without loss of generality that $C,C'$ are the first and last columns of $\mathscr T$, in keeping with our above change of notation.

In this case recall the definition of $\textbf{M}^*_s$ given in \ref {2.3.1}. It is an $n-s \times n-s$ minor fitting snugly between its first and last matrix blocks of \textbf{M} of size $s=\height C=\height C'$. View the entries of $\textbf{M}^*_s$ as functions on $\mathfrak m$ through the Killing form, that is to say as dual vectors.  Here the root vectors in $\textbf{M}^*_s$ that do not lie in $\mathfrak m^-$ are discarded\footnote{As noted in \cite [3.6.3]{FJ1} this leads to \textit{huge} economy of notation and thought since we not need to make $\textbf{M}^*_s$ depend on the sizes of the intermediate blocks.  This is likely to be even more significant \cite [3.6.2]{FJ1} outside type $A$.}.

By Lemma \ref {2.3.5} the $i^{th}$ Levi block contributes a factor of $a^{\max{(0,c_{i}-s)}}$ to the restriction of $\textbf{M}^*_s$ to  $\mathfrak m +a \Id$.

This result means that we must place $\max{(0,c_{i}-s)}$ copies of $a$ in the diagonal of the $i^{th}$ Levi factor of $\mathfrak p$ in any fashion and this for all $i \in ]1,k-1[$. Then the corresponding coefficient of $a^{d_{\mathscr D}}$, is exactly $I^s_{C,C'}$. Given that $\textbf{M}^*_s$ is a $n-s \times n-s$ minor, that is say a determinant of size $n-s$, the multiplication rules for calculating a determinant then give the required assertion.

\end {proof}

\textbf{Example $7$.} Consider the composition $(1,2,1)$.  View $\textbf{M}^*_1$ as the $3 \times 3$ determinant with entries $x^*_{i,j}:i \in [1,3],j\in [2,4]$.  Evaluate $\textbf{M}^*_1$ on $a\Id+\mathfrak m$ and just retain the terms for which $a$ appears on either the $(2,2)$ or $(3,3)$ entry.  The second gives $ax^*_{1,2}x^*_{2,4}$ and corresponds to the composite line $\ell_{1,2,4}$, the first gives $ax^*_{1,3}x^*_{3,4}$ and corresponds to the composite line $\ell_{1,3,4}$.  The lowest degree term in the expansion is $a(x^*_{1,2}x^*_{2,4}+x^*_{1,3}x^*_{3,4})=aI^1_{C_1,C_2}$.

 One might add that the second composite line does not leave the rectangle $R^1_{C_1,C_3}$.  In this special case there is just one.  In general the term used in constructing the Weierstrass section (\ref {5.2.6}) is always one composed from those that do not leave the rectangle $R^s_{C,C'}$.

\subsubsection {An Upper Bound on Degree}\label{4.3.3}

Let $C,C'$ be neighbouring columns of height $s$. Recall the notation of \ref {2.1.1} and write $C=C_l,C'=C_m$, for some $l,m \in [1,k]$.

Now let  the columns in $]C,C'[$  have arbitrary height.  Set $c_i=\height C_i$.

\begin {lemma} Let $M$ be a monomial which is a product of $s$ disjoint directed composite lines starting from $C$ and ending in $C'$.  Then
$$\deg M \leq (\sum_{i=l}^m \min {(s,c_i)})-s.\eqno{(*)}$$
\end {lemma}

\begin {proof}  Indeed the disjointness of the $s$ lines means that can pass through at most $s$ boxes of $C_i$ and of course at most $c_i$ boxes of $C_i$.  Thus the total number of boxes which these $s$ lines meet is at most $\sum_{i=l}^m \min {(s,c_i)})$.  Since there are $s$ component lines, one must subtract $s$ to obtain the number of constituent lines.
\end {proof}

Clearly $I^s_{C,C'}$ has degree $n-s-d_\mathscr D$, which in turn one checks is the right hand side $(*)$.  We define this to be the true degree $\deg I^s_{C,C'}$ of $I^s_{C,C'}$, that is to say we set
$$\deg I^s_{C,C'}:=(\sum_{i=l}^m \min {(s,c_i)})-s.\eqno{(1)}$$

\textbf{Remark}.  Notice we obtain the same result for the true degree of $I^s_{C,C'}$ if we suppress all the numbers not occurring in the rectangle $R^s_{C,C'}$ but sum over all the columns of $\mathscr T$.

\subsubsection {Virtual Degree}\label{4.3.4}

We define in \ref {4.3.6} a new tableau $\hat{\mathscr T}$ with respect to $I^s_{C,C'}$ by shifting some boxes to the left.

In addition following the Remark \ref {4.3.4} we suppress all entries not in $[C,C']$.  %(Probably it is enough to just suppress all those lying $R^s_{C,C'}$.)

It is clear that there will still be exactly $s$ disjoint monomial chains, and
that the monomials defined by the chains are unchanged. The only difference is a composite line need not be directed.  Indeed a constituent $l_{i,j}:i<j$ of this composite line may either join boxes in the same column or the composite line may zig-zag, that is say admit a constituent also going strictly from right to left. We call these constituent  lines, \textit{exceptional}.

Let $\hat{C}_i$ be the column of $\hat{\mathscr T}$ corresponding to the column $C_i$ of $\mathscr T$, by shifting some boxes to the left, as made explicit in \ref {4.3.6}. Set $\hat{c}_i:=\height \hat{C}_i$. We define the virtual degree $\deg \hat{I}^s_{C,C'}$ of $I^s_{C,C'}$ to be
$$\deg \hat{I}^s_{C,C'}=(\sum_{i=l}^m \min {(s,\hat{c}_i)})-s.\eqno{(2)}$$

\begin {lemma} Suppose $\deg \hat{I}^s_{C,C'}<\deg I^s_{C,C'}$, and that every exceptional constituent is an excluded root. Then $I^s_{C,C'}$ is zero when all the excluded roots are set equal to zero.
\end {lemma}

\begin {proof}

If there are no exceptional constituents, then the composite lines are directed and by Lemma \ref {4.3.3} the virtual degree is an upper bound on $\deg I^s_{C,C'}$.  Through the hypothesis of the lemma we obtain the required contradiction.

\end {proof}

\textbf{Remark 1.}  This is easy because we have swept under the carpet all the difficulties, which we must now tackle.

\textbf{Remark 2.}  Permit us to preview a consequence of the existence of a Weierstrass section. It defines (Theorem \ref {5.2.5}) a monomial in $I^s_{C,C'}$ of degree $\deg I^s_{C,C'}$ by Proposition \ref {5.1.1}(i) and of virtual degree exactly one less than this through Proposition \ref {5.1.1}(ii), taking account of \ref {4.1.2}.  Thus we have no option but to prove the \textit{very tight assertion} that
$$\deg \hat{I}^s_{C,C'}=\deg I^s_{C,C'}-1. \eqno {(3)}$$

\subsubsection {Attaching Excluded Roots to a Fixed Invariant}\label {4.3.5}

%A possible way to overcome the difficulties noted in \ref {4.3.2} is to prove vanishing for just one invariant, say  $I^s_{C,C'}$ by tailoring the set of excluded roots. This is of independent interest.

Fix a component tableau $\mathscr T^{\mathcal C}$.

\begin {lemma} With respect to $\mathscr T^{\mathcal C}$ and $I^s_{C,C'}$, there is a unique $i$-string starting in $ R_u \cap [C,C'[ $, for some minimal $u\leq s$ and passing to $R_{s'}\cap C''(\infty)$, with $s'>s$ and $C''\in ]C,C']$.

Moreover if $i$ lies strictly to the right of the rightmost column of height  $s'>s$ in $]C,C'[$, then $s'=s+1$ and $i$ enters $R_{s+1}\cap C'$.

\end {lemma}

\begin{proof}

Recall $\mathscr B^s_{C,C'}$ is a single batch and take $i$ in that batch as specified by the choice of $\mathscr T^\mathcal C$.

Consider the corresponding $i$ string in $\mathscr T^\mathcal C(\infty)$. By \ref {3.2.2}, Rule $1$, it passes into empty boxes of adjacent columns in $[C,C']$ as $\mathscr T^\mathcal C(\infty)$ is built from $\mathscr T$.  Then since $\height C=s$ it cannot cross $C$ and so starts in $R^s_{C,C'}$ in some row $R_u:=R_{(i)}$, where $u\leq s$.  Finally it enters $R_{s'}\cap ]C,C']$, for some $s'>s$.

Hence the first part of the lemma.

  Let $C_{\geq s}$ be the rightmost column of $[C,C'[$ of height $\geq s$ (possible $C$ itself)  The columns in $]C_{\geq s},C'[$ have height $<s$. Thus an entry in $]C_{\geq s},C'[$ must lie in $R_t:t<s$.  Then by \ref {3.2.2} such an entry cannot be moved strictly below $R_{s+1}$ until it leaves $]C_{\geq s},C'[$.  Moreover when it does enter $R_{s+1} \cap ]C_{\geq s},C'[$, it translates horizontally into $R_{s+1}\cap C'$.

   Hence the second part of the lemma.

   \end {proof}

  \textbf{Remarks}.  It can happen that $s'>s+1$ - see the bottom tableau in Figure $4$.  Here $i=4$ whilst $s=1,s'=3$.  It can also happen that $i$ is ``blocked'' from entering $R_{s'+1}\cap C'$.  In this example blockage  occurs if one introduces a column of height $3$ between $C_4$ and $C_5$.

  \

  \textbf{Definition and Notation}.  The unique $i$ - string associated by the lemma to the pair of neighbouring columns $C,C'$ (of height $s$) and  component data $\mathcal C$ is called the $i$-string of penetration, or penetrating $i$-string.  In this we set  $i=i^\mathcal C_{C,C'}$.

  \subsubsection {The Definition of $\hat{\mathscr T}_{I,\mathcal C}$}\label{4.3.6}

  Fix $s \in \mathbb N$ and a pair $C,C'$ of neighbouring columns $C,C'$ of height $s$.  Set $I=I_{C,C'}^s$.  To $I$ and some component data $\mathcal C$ we may associate an $i$-string of penetration by Lemma \ref {4.3.5}.

  To this data we associate a new tableau $\hat{\mathscr T}_{I,\mathcal C}$ (or simply, $\hat{\mathscr T}$) as follows.

  First the $i$ string of penetration starts in a unique box of $\mathscr T$ in the rectangle $R_{C,C'}^s$.  In $\mathscr T^\mathcal C$, this $i$-string moves to right, possibly just horizontally, but possibly also downwards by ``steps'', in going from a column to an adjacent column in $\mathscr T^\mathcal C(\infty)$.

Label the steps (defined by $i$ moving strictly downwards) by $u \in [1,t]$. This means that $i$ goes from a column $C_{r_u}(\infty)$ to an adjacent column $C_{r_u+1}(\infty)$ and goes down by $m_u\geq 1$ rows.

Set $t_u=\height C_{r_u+1}$, for all $u \in [1,t]$.  Let $j_{r_u+1}$ (or simply, $j$) be the unique lowest element in $C_{r_u+1}$ and set $j^{m_u}_{r_u+1}=\{j-m_u+1,j-m_u+2,\ldots, j\}$.  In the notation of \ref {4.1.3}, these are the entries of the partial column $C_{r_u+1}^{>t_u-m_u}$.  Observe that $C_{r_u+1}^{>t_u-m_u}$ is a partial column with $m_u$ rows.  Expressed briefly $\height C_{r_u+1}^{>t_u-m_u}=m_u$.

Now form the partial column $C^{\hat{m}}$ with $\hat{m}=\sum_{u=1}^t m_u$ rows by empiling the partial columns $C_{r_u+1}^{>t_u-m_u}:u \in [1,t]$, the first going on top and so downwards as $u$ increases.  As we shall see partly below and in \ref {4.3.9}, the above order is vital to preserve the set of primary excluded roots.

On the other hand, $i$ appears in $C_h:=C_{(i)}$ in some unique row $R_f=R_{(i)}$ of $\mathscr T$.

By the claim of \ref {4.1.3} applied to the case $u=1$, there exists a column $C_g$ of height $f$ for some unique largest $g \leq h$.

If $\height C_h >f$, let  $C_{g_v},C_{g_{v-1}},\dots,C_{g_1}:h=:g_v>g_{v-1}>\dots\geq g_1$ be the columns of height $>f$ between $C_h,C_g$.  Set $g=g_0$.

  In a manner similar to the construction in \ref {4.1.3}, define $\hat{\mathscr T}$ by first removing $C_{r_u+1}^{>t_u-m_u}$ from $C_{r_u+1}$ for all $u \in [1,t]$ and then forming $C^{\hat{m}}$ as in the fifth paragraph above. Then $C^{\hat{m}}$ is placed below $i\in R_f\cap C_h$, displacing the partial columns $C_{g_s}^{>f}:s=v,v-1,\ldots,1$  successively to the left skipping over the columns of height $<f$, till $C_{g_0}$ is reached.

 %Let $\hat{C}$ be the column of $\hat{\mathscr T}$ corresponding to the column $C$ of $\mathscr T$.

 For any column $C$ of $\mathscr T$, let $\hat{C}$ be the column of $\hat{\mathscr T}$ replacing $C$ in the above construction. In particular $\hat{C}_{g_v}$ is obtained by placing  $C^{\hat{m}}$ on top of $C_{g_v}^{\leq f}$ and by our choice of ordering in  $C^{\hat{m}}$ its entries decrease strictly on going up the rows.

 Thus the only difference between $\hat{\mathscr T}$ and any one of the $\mathscr T_{i,\textbf{j}}$ is that $C_{(i)}^{>f}$ has been replaced by the ``amalgamated'' partial column $C^{\hat{m}}$ instead of one of the partial columns $C_{r_u+1}^{>t_u-m_u}:u \in [1,t]$.

 Let $\hat{w}$ denote the Weyl group element defined by $\hat{\mathscr T}$ as in \ref {4.1.4} with respect to which one further defines as in \ref {4.1.4} the excluded roots associated to $\hat{\mathscr T}$.

 \

 Thus the primary excluded roots coming from  pairs $(i,j):i<j$  with both $i,j$ in $\hat{C_{g_v}}$, exactly arise when $i$ is above $j$ in $\hat{C_{g_v}}$.  (Some may define root vectors in $\mathfrak l^+$ and are ignored - see Definition \ref {4.1.4}).

 \

 $(*)$.  For all $i \in [1,v-1]$, $\hat {C}_{g_{i-1}}$ is the composite column with $C_{g_{i+1}}^{>f}$ placed below $C_{g_i}^{\leq f}$.

  %  $C=C_{g_s}; s \in [1,v-1]$, then $\hat{C}$ is the composite column with $C_{g_i}^{\leq f}$ below and $C_{g_{s+1}}^{>f}$ placed on top.  For later reference we note that
%
% \
%
% $(*)$.  One has $\height \hat{C}=\height C$, when $C=C_{g_s}; s \in [1,v-1]$.

 \

   %In \ref {4.1.3} we formed the tableau $\mathscr T_{i,\textbf{j}}$.  The present tableau is an amalgam of the tableaux $\mathscr T^{r,\textbf{j}}$ as $\textbf{j}$ runs through the partial columns $C_{r_u+1}^{>t_u-m_u}:u\in [1,t]$.

  \textbf{Remark.} The above amalgamation is \textit{mandatory} to show, by our present method, that the present  \textit{triply infinite exponentially growing family} of determinants all become zero on setting the excluded roots equal to zero.  We believe the proof is short, rather elegant and certainly much simpler than anything else imaginable.

% \textbf{Notation} We use $j \in C^{\hat{m}}$ to denote one of the partial columns $j^{m_u}_{r_u+1}:u \in [1,t]$.

  \subsubsection {Changes of column heights}\label{4.3.7}

 Recall that the $i$-string of \ref {4.3.6} is penetrating and so passes from $R^s\cap C_{r_u}(\infty)$ into $R_{s'}\cap C_{r_u+1}(\infty)$ for some  $s'\geq s+1$ at the last step.

 The following lemma describes the changes in column heights on going from  $\mathscr T$ to $\hat{\mathscr T}$.

 \begin {lemma}

 \

 $(i)$. $\height \hat {C}_{r_u+1}=\height C_{r_u+1}-m_u$, for all $u \in [1,t]$.

 \

 $(ii)$. $\height \hat {C}_{g_v}=\height C_{g_0}+\hat{m}$.

 \

 $(iii)$. $\height \hat{C}_{g_{i-1}} = \height C_{g_i}$, for all $i \in [1,v-1]$.

 \

 $(iv)$. $\height \hat{C}_{g_v}=s'$.
% In going from  $\mathscr T$ to $\hat{\mathscr T}$, for all $u\in [1,t]$, the column  $C_{r_u+1}$ has height reduced by $m_u$, whilst $\height \hat{C}_{g_v}=\height C_{g_0}+ \hat{m}$.  The remaining column heights are permuted.  Finally $\height \hat{C}_{g_v} =s'$.
 \end {lemma}

 \begin {proof} $(i)$ is an immediate consequence of removing a partial column of height $m_u$ from the bottom of $C_{r_u+1}$.

 $(ii)$ is an immediate consequence of placing a partial column of height $\hat{m}$ below $C_{g_v}^{\leq f}$, the latter having height $f= \height C_{g_0}$.

 $(iii)$ follows from \ref {4.3.6}$(*)$.

 $(iv)$. The $i$-string passes from $R_f \cap C_{g_v}$ penetrating into $R_{s'}:s'\geq s+1$ at the last step.  Thus $\hat{m} = s'-f= s'-\height C_{g_0}$, so by $(ii)$ we obtain $\height  \hat{C}_{g_v} =s'\geq s+1$.
 \end {proof}

  \textbf{Remarks}.  When we delete the entries not between $C,C'$ then this just effects the heights of the columns $C_i$ lying strictly to left of $C$ and the corresponding modified columns $\hat{C}_i$ .

  Concerning the columns strictly to the left of $C$, we claim that their contributions cancel out when entries strictly to the left are also deleted.

   Indeed with respect to \textit{just the entries lying strictly to the left of} $C$, the heights of columns of height $\geq f$ are cyclically permuted as we pass from $C_i$ to $\hat{C}_i$.  This follows as in the proof of  \ref {4.3.6}$(*)$.

    Thus in calculating the difference between the virtual and real degrees, the result is the same had we neither deleted the columns strictly to the left of $C$, nor their entries.  However \textit{a priori} this is an incorrect procedure!

 \

 \textbf{Example $8$}.  Consider the composition $(1,3,2,1,2)$ and the component in which $5$ goes down one step below $7$ and then one step below $9$.  Then $C_1,C_2$ have entries $\{1\},\{2,3,4\}$ respectively whilst $\hat{C}_1,\hat{C}_2$ have entries $\{1,3,4\},\{2,6\}$ respectively.  Retaining only the entries $\{1,2,3,4\}$ which lie strictly to the left of $C$, it is immediate that the last assertion of the Remark holds, namely that heights are cyclically permuted (in this case just swapped) as the columns acquire a hat.

 \subsubsection {Verifying Equation $(3)$}\label{4.3.8}

\begin {lemma}  Equation $(3)$ holds.
\end {lemma}

\begin {proof}  Retain the notation of \ref {4.3.7} and recall the conclusion of Lemma  \ref {4.3.7}.

The contribution of $(iii)$ of the lemma to  $\deg \hat{I}^s_{C,C'}-\deg I^s_{C,C'}$ is zero, since sums are just permuted.

Again the total decrease in column heights from (i) Lemma  \ref {4.3.7} is $\hat{m}$ which is exactly compensated by the increase of the column height coming from (ii).

Yet we still must take account of having to compare minima with respect to $s$.

Indeed $\min{(\height \hat{C}_{g_v},s)}=s$, which is a drop of $s'-s$ from $\height \hat{C}_{g_v}$.

On the other hand the penetrating $i$-string only penetrates the rectangle $R^s_{C,C'}$ at the last step, that is when $u=t$.  Thus we obtain $\height C_{r_{u}+1}< s$, for $u\in [1,t-1]$, whilst $\height C_{r_t+1}=s'-1\geq s$.

Combined these combined contributions to $\deg \hat{I}^s_{C,C'}-\deg I^s_{C,C'}$ sum to $-1$, even when $s'-1>s$.  This  recovers $(3)$.

\end {proof}
\textbf{Example $9$}. Consider the composition $(2,1,2,1)$ and the invariant which comes from the pair $C_2,C_4$ of neighbouring columns of height $1$.

There is a component in which $4$ goes down one step below $6$.   Here $s=1, s'=2$.

%Then $6$ is placed below $4$ displacing $5$ below $3$. Then $\height C_2=1$ and $\height \hat{C_2}=\height C_3$, whilst $\height \hat{C_3}=s'$.  Thus $\min{(s,\height \hat{C_3})}=s$, $\height \hat{C_4}=\height C_4 -(s'-s)=0$, which recovers Equation $(3)$.

There is a second component in which $3$ goes below $5$.  Here $s=1,s'=3$.

%Then $(4,5)$ goes below $3$.  Then $\height \hat {C_2}=s'$, so $\min{(s,\height \hat {C_2})}$ which  is a drop of $s'-s=2$. Yet $min{(s,\height C_3)}=1$ which is compensated by a drop of $1$ from $\height C_3$.

In both cases one recovers Equation $(3)$.

\subsubsection {Primary Excluded Roots}\label{4.3.9}

Recall Definition \ref {4.1.4} and the notation of \ref {4.3.6}.

\begin {lemma}  The primary excluded roots defined by $\hat{\mathscr T}$ are the same as those coming from the $\{\mathscr T_{(i,j)}: j \in C^{\hat {m}}\}$.
\end {lemma}

\begin {proof} This is because the primary excluded roots coming from $j \in C^{\hat {m}}$ being placed below $r$ in $C_{(r)}$ as a part $C^{\hat {m}}$, are the same  as those when we just place $j$ below $r$ in $C_{(r)}$.

In the notation of \ref {4.3.6}, take $u,u'\in [1,t]$ with $u<u'$.   Then with respect to the first paragraph above, in the first (resp. second)  case using
$\hat{\mathscr T}$(resp. $\hat{\mathscr T}_{i,\textbf{j}}$), $i' \in C_{r_u+1}^{>t_u-m_u}$ lies strictly above (resp. to the right) of  $j' \in C_{r_{u'}+1}^{>t_{u'}-m_{u'}}$.

Recalling Observation $(1)$ of \ref {4.1.4} in either case $\alpha_{i',j'}$ is a primary excluded root through \ref {4.1.4}$(*)$ and the way that the word form of a Weyl group element is constructed from the columns of $\hat{\mathscr T}$ (resp. $\hat{\mathscr T}_{i,\textbf{j}}$).

\end {proof}

%This would have ``Gang aft agley, An' lea'e us nought but grief an' pain'' \cite [To a Mouse] {B} had order been reversed in the empilation of columns in  $\hat{\mathscr T}$.

\textbf{Example $10$}.  Consider the composition $(2,1,1,2)$ and the canonical component tableau in which $3$ goes down twice by one row.  The construction of \ref {4.1.4} gives two tableaux $\mathscr T_{3,4}$ and $\mathscr T_{3,6}$, with excluded roots $\alpha_{3,4}$ (resp. $\alpha_{3,6},\alpha_{4,6}$), all primary in both cases. (Here $\alpha_{5,6}$ is not counted since it is a root of the Levi factor.) These are combined in $\hat{\mathscr T}$, which has a column in which $4$ is placed below $3$ and then $6$ below $4$, but this fails if the latter order is reversed, where crucially $\alpha_{4.6}$ given by a vertical line, fails to be an excluded root.

 %For this is vital to make use of the given order in which the $C_{r_u+1}^{>t_u-m_u}:u \in [1,t]$ were arranged in $C^{\hat{m}}$.

\subsubsection {Secondary Excluded Roots}\label{4.3.10}

\begin {lemma} The secondary excluded roots defined by $\hat{\mathscr T}$ are the same as those coming from the $\{\mathscr T_{i,j}: j \in C^{\hat {m}}\}$.
\end {lemma}

\begin {proof} This is because the same column $C_{(r)}$ is entered into for all cases.
\end {proof}

\subsubsection {Warning}\label{4.3.11}

Yet the set of excluded roots in $\mathscr T^\mathscr C$ can still be strictly decreased in the presence of several pairs of neighbouring columns, as noted below.

\

\textbf{Example $11$}. Consider the composition $(1,2,1,2)$.  For $\mathscr T^\mathcal C$ put $2$ below $4$ in $C_3$ and further put $3$ below this new appearance of $2$ in $C_3$.  Then $\mathscr T_{3,4}$ gives enough excluded roots so that \textit{both} $I^1_{C_1,C_3}$ and $I^2_{C_2,C_4}$ vanish.  Yet our general recipe includes the excluded roots coming from $\mathscr T_{2,4}$ which introduces an extra secondary excluded root $\alpha_{1,3}$ superfluous for vanishing but which does not give the wrong dimension for $\overline{B.\mathfrak u^\mathcal C}$, as one may check!  (This also follows in a more general context from Proposition \ref {6.2.4}, since $\ell_{1,2}$ ``covers'' this secondary excluded root.)

On the other hand had we taken the canonical component, then the excluded root $\alpha_{1,3}$ \textit{not} needed for the vanishing of $I^2_{C_2,C_4}$ \textit{is} needed for the vanishing of $I^1_{C_1,C_3}$.

The purist might ask ``is there a minimal set of excluded roots''.  This can be tricky.  Indeed consider the

%\textbf{Example 2.} Take the composition $(3,2,1,2,2,1,3)$ with $7$ being moved under $9,10$ and $5$ being moved under $8$ and then under $7$.  Then the exclude roots coming from the $5$-string ensure the vanishing of all the invariants except $I^1_{C_3,C_6}$, for which we need to move $9,10$ under $7$ displacing $8$ to the left.

\textbf{Example $12$.} Take the composition $(3,2,2,1,1,2,3)$ and $8$ being moved under $9$ and then $11$, and $5$ being moved under $7$ and then under $8$.  Then the excluded roots coming from the $5$-string ensure the vanishing of all the invariants except $I^1_{C_4,C_5}$, for which we need to move $9$ under $8$.

\subsubsection {Exclusion of Exceptional Constituents}\label {4.3.12}

Consider a  monomial chain \newline $(i_1,i_2,\ldots,i_u)$ of strictly increasing integers occurring in boxes of $\mathscr T$.  In $\hat{\mathscr T}$ some of these boxes are moved to the left. This can cause the resulting chain to admit an exceptional constituent, that is to say for some pair $(i_u,i_{u+1})$ the box containing $i_{u+1}$ lies  to the left (not necessarily strictly) of that containing $i_u$.  This can only happen if in the notation of \ref {4.3.6} one has $i_{u+1} \in \hat{C}_{g_t}$ for some $t\in [1,v]$, and $i_u$ lies in a column of $\hat{\mathscr T}$ to the right (not necessarily strictly) of $\hat{C}_{g_t}$.  However these roots are excluded roots (primary if $t=v$ and secondary if $t<v$) and so this monomial chain in $\hat{\mathscr T}$ is zero when the excluded roots are evaluated to zero.  Thus in Lemma \ref {4.3.4} a monomial in $I^s_{C,C'}$ is set to zero if it admits an exceptional constituent.

\subsection {Specific Vanishing}\label{4.4}

Let $E^\mathcal C_{C,C'}$ be the set of excluded root vectors coming from the penetrating $i$-string $i^\mathcal C_{C,C'}$ defined for the Benlolo-Sanderson invariant $I^s_{C,C'}$.

\begin {prop} $I^s_{C,C'}$ vanishes when the $E^\mathcal C_{C,C'}$ are set equal to zero.
\end {prop}

\begin {proof} Combine Lemmas \ref {4.3.4} and Lemmas \ref {4.3.8}-\ref {4.3.10} and the observation in \ref {4.3.12}.
\end {proof}

\subsection {Global Vanishing}\label{4.5}

Since Proposition \ref {4.4} applies to any Benlolo-Sanderson invariant we obtain the

\begin {thm}  For all parabolics in type $A$ and for all component tableaux, all the Benlolo-Sanderson invariants become zero when all the excluded root vectors are set equal to zero.
\end {thm}

Recall (\ref {4.3.1})  the subvariety $\mathscr C$ of $\mathfrak m$.

\subsection {Containment in $\mathscr N$}\label{4.6}

\begin {cor}  One has $\mathscr C \subset \mathscr N$.
\end {cor}

\subsection {Limits to Amalgamation}\label{4.7}

Not all the tableaux $\mathscr T_{i,\textbf{j}}:i,\textbf{j} \in \textbf{L}$ can be amalgamated into a single tableau.  Since the latter defines $w\in W$ such that the resulting component is $\overline {B.(\mathfrak n\cap w\mathfrak n)}$, this is possible if and only if $\overline {B.\mathfrak u^\mathcal C}$ is an orbital variety closure which is not always the case.  Indeed we believe that components which are orbital variety closures will be on average rare. Example $13$ describes a case where full amalgamation is not possible.  The calculation is indicated in Figures $9,10$.

  Thus it is already a miracle that we only need the roots coming from the penetrating $i$-string to prove vanishing because had we needed to adjoin further roots, we could not have reduced our analysis to a single tableau $\hat{\mathscr T}$.
  %Yet a single tableau attached to a component gives an element of the Weyl group making this component an orbital variety closure.  The latter may fail to hold  (see \ref {8.2}.  (We believe that components which are orbital variety closures will be on average rare.)

  A second even bigger miracle is that Proposition \ref {4.4} holds.  This result is crucial to proving the injectivity of the component map (\ref {7.3}).

\subsection {Example $14$.}\label{4.8}

It is highly advisable to compute some examples of the labelling of the root vectors in \textbf{M} by $1$, $\ast$ and their encirclement.  In Example $14$ we give the result in Figure $12$ for the five different component tableaux which result from the parabolic defined by the composition $(2,1,2,1,2,1)$.  The reader may also wish to show that the all the six Benlolo-Sanderson invariants vanish on converting the circles to zeros, for each of the five component tableaux - an easy computation but not so easy when $n$ is large (see Section \ref {4.3}).

 \subsection {Precursor}\label {4.9}

The origins of \ref {4.4} lay in the following observation.  Fix a pair of neighbouring columns $C,C'$ of \textit{maximal} height $s$ and a component tableau $\mathscr T^\mathcal C$.  Let $F^\mathcal C_{C,C'}$ be the subset of \textit{primary} excluded roots given by the lines with label $\ast$ coming from the penetrating $i=i^\mathcal C_{C,C'}$ string in $\mathscr T^\mathcal C(\infty)$. It is a subset of $E^\mathcal C_{C,C'}$. Let $\tau$ be the involution on columns of $\mathscr T$ defined by reversing the order of entries in each column.  Let $C_{(i)}$ be the unique column of $\mathscr T$ containing $i$ and set $j:=\tau(i), m=n-s$. Let $\textbf{M}_s$ be the $m\times m$ sub-matrix of  \textbf{M} in the upper right hand corner.  It contains $\mathfrak m$.  Evaluate to zero all the entries of $\textbf{M}_s$ not in $\mathfrak m$ as well as those in $F^\mathcal C_{C,C'}$.  A simple exercise shows that the resulting matrix has $j$ columns with $\leq j-1$ rows with non-zero entries.

If $j=1$, then $\det\textbf{M}_s=0$ trivially. Thus assume $j>1$.

Now set $j'=m-j$.  By the above $\textbf{M}_s$ has a $(j'+1)\times (j-1)$ sub-matrix in its lower left hand corner all of whose entries are $0$. Thus $\det \textbf{M}_s$ factors as the product of the $(j-1) \times (j-1)$ minor in its upper left hand corner and a $(j'+1) \times (j'+1)$ minor in its lower right hand corner, whose first column has only $0$ entries. Thus $\det \textbf{M}_s=0$, in this case too.

This gives the conclusion of Proposition \ref {4.4} with respect to the Benlolo-Sanderson invariant corresponding to a pair of neighbouring columns of maximal height.

As an example consider the left hand matrix in Figure $19$. One has $i=7,\tau(7)=8$ Ignore the zeros lying in column $8$ of \textbf{M}, because these do not come from the penetrating $i$-string. Then columns $1-7,11$ of $\textbf{M}_s$ have only zeros in rows $7,9-14$.  Notice that in this $x_{8,14}$ is \textit{not} an excluded root and indeed putting $14$ under $7$ pushes $8$ strictly to the left.  All this underlies the very tight nature of our vanishing result.

We attempted unsuccessfully to reduce to the above case by using the secondary excluded roots.

\section {Weierstrass Sections} \label {5}

 \subsection {General Remarks}\label {5.1}

 Let $\mathscr T^\mathcal C$ be a component tableau.  Let $e_\mathcal C$ (resp. $V_{\mathcal C}$) be the sum of the root vectors (resp. the sum of the root subspaces) labelled by a $1$ (resp. by a $\ast$) in $\mathscr T^\mathcal C$.

We shall often drop the subscript $\mathcal C$.  Again $\mathscr T^\mathcal C$ is constructed from a limiting tableau $\mathscr T^{\mathcal C}(\infty)$. Here and in the next two sections we may often drop the superscript writing simply $\mathscr T(\infty)$.

We aim to show that to $\mathscr T^\mathcal C$ we may associate a component $\mathscr C$ of $\mathscr N$ containing $e_\mathcal C$ with $e_\mathcal C+V_\mathcal C$ being a Weierstrass section.

For the latter we need the two conditions described below.

 \subsubsection {Composite lines}\label {5.1.1}

 \

 \textbf{Definition.}

  A directed composite line is a concatenation of lines  joining boxes in $\mathscr T$ all going strictly from left to right.  Composite lines are said to be disjoint if they do not go through the same boxes. The individual lines which make up the composite lines are called its constituents.  (They are directed by the concatenation going from left to right.)

 \

 \textbf{N.B.}.  In contrast to section \ref {4.3}, we \textit{only} consider throughout Section \ref {5} composite lines in $\mathscr T(\infty)$ (resp. in $\mathscr T$) whose constituents have neutral labels or are labelled by $1$ or $\ast$ and are directed.  We shall see from \ref {5.2.4} that this corresponds to selecting just one monomial for each Benlolo-Sanderson invariant.

 \

  To obtain a Weierstrass section (for the action of a parabolic on its nilradical in type $A$) we will prove the following.

 \begin {prop}  Between any two neighbouring columns $C,C'$  of height $s$ there is

 \

 (i). A unique disjoint union of  composite lines passing through all the boxes in $R^s_{C,C'}$.

 \

 (ii). The constituent lines which make up the disjoint union, exactly one $\ell^\ast_{C,C'}$ is labelled by a $\ast$, the rest are labelled by a $1$.
 \end {prop}

 \textbf{Remarks}. Evaluate the Benlolo-Sanderson invariant $I_{C,C'}$ on $\textbf{1}+e+V$ by replacing $x_{i,j}$ by $1$, if $\ell_{i,j}$ is labelled by a $1$, by $v_{i,j}$, if $\ell_{i,j}$ is labelled by a $\ast$ and if $\ell_{i,j}$ has no label $0$.  Then by (i),(ii) its evaluation is $v_{i,j}$ \cite [Lemma 4.2.5]{FJ1}.

% When there are no boxes strictly below $R_s$ between $C,C'$, this is particularly easy and only requires a knowledge of how a determinant is multiplied out.
%
% %\footnote{As late as 1928, Schur devoted an entire semester to a course on determinants and a course on invariants \cite [p. xxiii, personal remininiscences of W. Ledermann] {JMR}, so had we lived in an earlier age, Schur might have honoured us with his imprimatur.}
%  Otherwise a little more knowledge of how $I_{C,C'}$ is defined is needed.

 Condition (i) on its own is easy to satisfy. It results for example by just joining all adjacent boxes in $\mathscr D$ by a horizontal line and taking the sum $e$ of the root vectors (which all lie $\mathfrak m$).  The nilpotency class for this choice of $e$ is easy to compute from which one obtains $\dim G.e = 2
\dim \mathfrak m$.  Yet (borrowing a trick from the proof of Richardson's theorem) one has $\dim G/P = \dim \mathfrak m$ and so $\dim P.e \geq \dim G.e -\dim \mathfrak m=\dim \mathfrak m$.  This means that $P.e$ is dense in $\mathfrak m$, that is to say $e$ is a Richardson element.  This construction (of an explicit Richardson element in $\mathfrak m$) in Type $A$ is due Ringel et al \cite {BHRR}.  It caused quite a stir at the time and was generalised to biparabolics in all types \cite {JSY}, except curiously it fails for at least one biparabolic in type $E_8$.

In more generality, if $e\in\mathfrak m$, then $G.e\cap \mathfrak m \supset P.e$ and so after Spaltenstein (cf \cite [Eq. $(1)$]{FJ1} one has $\dim P.e \leq \frac{1}{2} \dim G.e$, yet equality can fail \cite [6.10.9. Example 4]{FJ2}.

In \cite [Sect. 5]{FJ2} we gave a horrifically complicated procedure to modify the Ringel et al lines so that (ii) also held.  Then in \cite [Sect. 5]{FJ4} we found a \textit{much} simpler (but totally mysterious) way to recover this result using the ``composition tableau''.  Now by modifying the construction we include the composition tableau in a rich new family of ``component tableaux'', with magical, mysterious and miraculous properties.

 \subsubsection {Separation}\label {5.1.2}

 Proposition \ref {5.1.1} is not quite enough on its own to obtain a Weierstrass section because it could happen that the same co-ordinate of $V$ is obtained twice. What we need is:

 \

 $(*)$.   The $\ell^\ast_{C,C'}$ are pairwise distinct.

 \

 \textbf{Remark} This can be  satisfied  \cite [$(P_2)$] {FJ1} by requiring that exactly one of the lines joining boxes in $R^s_{C,C'}$ labelled by a $\ast$ meets $R_s$ and we used this property to define $\ell^\ast_{C,C'}$ and thereby obtain a Weierstrass section for the canonical component \cite {FJ4}.  Yet this property can fail for an arbitrary component tableau.  We shall replace disjoint composite lines in $\mathscr T$ by disjoint composite lines in $\mathscr T(\infty)$ and then miraculously the whole analysis in the proof of the Proposition \ref {5.1.1} drastically simplifies and in this new context we also obtain $(*)$.

 \

 \textbf{Example $15$.}  Consider the composition $(2,1,1,1,2)$.  In this $\mathscr B^2_{C_1,C_5}=\{2,3,4\}$ giving three tableaux and corresponding to Benlolo-Sanderson invariant $I^2_{C_1,C_5}$.  The first and last have indeed a line with a $\ast$ which meets $R_2$, but for the middle one the lines labelled by a $\ast$ lie entirely in $R_1$.

 \subsection {Collapsing}\label {5.2}

\subsubsection {Composite and Special Composite lines in $\mathscr T(\infty)$}\label {5.2.1}

Recall that in \ref {3.2.6} we defined a neutral line between adjacent columns to be one joining boxes with the same entries. These are exactly the lines which form an $i$-string in $\mathscr T(\infty)$.

%Consider an $i$ - string in $\mathscr T(\infty)$.

By \ref {3.2.4} - \ref {3.2.6} the right end-point of an $i$ string may end with a vertical line in $\mathscr T(\infty)$ with label $\ast$ at a box with an entry $\psi_{\ast}(i)$ or with a line to a left adjacent column in $\mathscr T(\infty)$  with label $1$ ending in a box (different by Corollary \ref {3.2.6}(i)), with an entry $\psi_1(i)$, different to $\psi_{\ast}(i)$ (though beware that we shall often drop the subscript on $\psi$).  The constituents of such a  composite line $\mathscr L$ have all but one a neutral label.

A line with a neutral label in such a composite line $\mathscr L$  may be collapsed, that is to say its end-points identified.  Then $\mathscr L$ becomes a line $\ell$  in $\mathscr T$ joining (distinct) boxes with end-points $i, \psi_\ast(i)$ (resp. $i, \psi_1(i)$) and labelled by a $\ast$ (resp.  by a $1$). These boxes will never be in the same column and indeed may lie in columns which are quite far apart.

Then collapsing neutral lines of such a composite line in $\mathscr T(\infty)$ joining  boxes with entries $i,\psi(i)$ in $\mathscr T(\infty)$ gives a composite line $\mathscr L$ in $\mathscr T$ joining (unique) boxes with entries $i,\psi(i)$.

To avoid the resulting zippimg across possibly far-distant columns, it is highly advantageous to stay in $\mathscr T(\infty)$ for part of the calculation.

We may further concatenate such lines $\mathscr L$ to form composite lines in $\mathscr T(\infty)$ still (cf \cite [4.2.3]{FJ1}) requiring that concatenated lines pass from left to right or down to up.

%Again we may speak of a disjoint union of composite lines in $\mathscr T(\infty)$.

%Then collapsing neutral lines of such a composite line in $\mathscr T(\infty)$ joining  boxes with entries $i,\psi(i)$ in $\mathscr T(\infty)$ gives a composite line $\mathscr L$ in $\mathscr T$ joining (unique) boxes with entries $i,\psi(i)$.

Finally the constituents of a composite lines may only have neutral labels (joining boxes with entries $i$) and such a line collapses to a single box with entry $i$.

\

\textbf{Definition.} A special composite line in $\mathscr T(\infty)$ is one consisting of only neutral lines or lines labelled by a $1$.  Examples are easy to find.

\subsubsection {Lines Between Neighbouring Columns}\label {5.2.2}

Consider now a pair of neighbouring columns $C,C'$ of height $s$ and the rectangle $R_{C,C'}^s$ they define.

For all $i \in [1,s]$, let $b_i$ be the entry of $R_i\cap C$.

Label the columns between $C,C'$ as $\{C_i:i \in [1,k]\}$.

Fix numerical data $\mathcal C$ to define $\mathscr T^{\mathcal C}(\infty)$.  (Here we shall mainly drop the superscript.)

Recall that $\mathscr T$ is a sub-tableau of $\mathscr T(\infty)$.

\begin {prop}

\

(i).  From every box in $C$ except exactly one, to be denoted by $b_m$, there is a unique special composite line $\mathscr L_i:i \in [1,s]\setminus \{m\}$ in $R^s\cap \mathscr T(\infty)$ to a box in $C'$.

%\
%
%Let the box in conclusion of (i) be denoted by $b_m$.

\

(ii). There is a special composite line from $b_m$ which first goes into $R_{t'}\cap C_r(\infty):t'\leq s, r \in [1,k-1]$, and then by a neutral line into  $R_{t+1}\cap C_{r+1}(\infty):t\geq s$ and some column $C_{r}\in [C,C'[$. It can then be joined by a line $\ell^\ast_{C,C'}$ with label $\ast$ to an entry in $b'':=R_{t''}\cap C_{r+1}$, with $t'' \leq s$.

\

(iii). One may choose $\ell^\ast_{C,C'}$ so that the $\mathscr L_i:i\neq m$ do not pass through $b''$.

\

(iv).  One may form the composite line $\mathscr L_m$ by concatenating the composite line from $b_m$ to $b''$ and then by a special composite line to $C'$.

\

(v).  The composite lines in (i),(iv) collapse to disjoint composite lines passing through every box of $R^s_{C,C'}$ in $\mathscr T$.

\

(vi). Either (Case one), $\ell^\ast_{C,C'}$ has right end-point in $R_s\cap C_{r+1}$, or (Case two) $\height C_{r+1}\leq s$ and $\ell^\ast_{C,C'}$ joins $R_{s+1}\cap C_{r+1}(\infty)$ to the lowest entry of $C_{r+1}$.
\end {prop}

\begin {proof}

$(i)$.  For all $t\in [1,s],l\in [1,k-1]$, the unique $u$-string defined by the entry of $R_t\cap C_l$, either enters $R^s\cap C_{l+1}(\infty)$, or is stopped at $C_l$ defining a unique line with label $1$ to $R^s\cap C_{l+1}$ having right end-point $u':=\psi_1(u)$, or by Lemma \ref {4.3.5}, for some unique value $r$ of $l$ (specified by $\mathcal C$) $u$ enters strictly below $R_s$.  The uniqueness of this evolution proves (i) and specifies $b_m$ as the unique starting box of $C$ for which a $u$-string descends strictly below $R_s$.

%$(i)$. Recall Remark \ref {3.2.2}, that every column $C''(\infty)$ of $\mathscr T(\infty)$  has a most one box with a given entry.
%
%  Consider a column $C_{r'}$ in $[C,C'[$ and a box in $R^s \cap C_{r'}(\infty)$ with entry $i$. By the rules in \ref {3.2.2}, $i$ may enter a unique box of $C_{r'+1}(\infty)$ or be stopped at $C_{r'}$.  In the latter case, by \ref {3.2.6}, there is a line $\ell_{i,j}$ with label $1$ to the entry $j$ of a unique box $b$ in $C_{r'+1}$, so being the right end-point of $\ell_{i,j}$. Moreover by the last sentence of the Remark \ref {3.2.6} there can be no other line with label $1$ with right end-point $b$.
%
%  On the other hand by Lemma \ref {4.3.5} there is a unique column $C_{r}$ in $[C,C'[$ and a unique entry $u$ of $R^s \cap C_r(\infty)$ such that $u$ enters strictly below $R_s$ - simply $u$ is the unique choice an entry of $\mathscr B^s_{C,C'}$ specified by $\mathscr C$.
%
%  Applying these conclusions to successive columns starting from $R^s\cap C$ gives (i).

$(ii)$.  Note that there exist $t'\leq s, t\geq s$, a unique special composite line from $b_m$ to $R_{t'}\cap C_r(\infty)$ and a neutral line to $R_{t+1} \cap C_{r+1}(\infty)$ (with end-points $u$) going down by $t+1-t'$ rows.

 % $(ii)$.   By the penultimate paragraph above there exist $t'\leq s$ and a unique special composite line from $b_m$ first to $R_{t'}\cap C_r(\infty)$ and then there exist $t>s$ and a neutral line to $R_{t+1} \cap C_{r+1}(\infty)$. Let $u$ be their common entry, that is to say the end-points of this neutral line.

  Then there are two cases to consider.

  \

  \textbf {Case one.} This is when $t-t'>0$, that is  the common entry $u$ goes down by $t+1-t'>1$ rows to $R_{t+1}$.  In this by \ref {3.2.2} we must have $\height C_{r+1}=t$.

%  First when $t-t'>0$, that is  the common entry $u$ goes down by $t+1-t'>1$ rows to $R_{t+1}$.  In this by \ref {3.2.2} we must have $\height C_{r+1}=t$.

  Then (see \ref {3.2.4}) there is a vertical line labelled by a $\ast$ from $R_{t+1} \cap C_{r+1}(\infty)$ to each of the boxes
  $ R_{t''}\cap C_{r+1}: t'' \in [t,t']$.  Choose that with $t''=s$ and let $v$ be the entry of $R_s\cap C_{r+1}$.

  \

   $(*)$. By the last sentence of \ref {3.2.6}(ii) there \textit{might} also be a right-going line\footnote{Namely, the upper dashed/dotted blue line in Figure $1$.} with label $1$ having right end-point in $R_s\cap C_{r+1}$, but then its left end-point must lie in $R_{s'}\cap C_r(\infty):s'> s$, so by $(i)$, it is not joined by any special composite line to $C\setminus \{b_m\}$.

   \

   \textbf{Case two.}  This is when $t'=t=s$, that is the entry goes down by exactly one row.  Then (see \ref {3.2.4}) there is a unique vertical line from $R_{s+1} \cap C_{r+1}(\infty)$ to $R_{t''}\cap C_{r+1}$, where $t'' = \height C_{r+1}$, so whose entry $v$ is the lowest entry of $C_{r+1}$.

   \

   $(**)$. Again by the last sentence of \ref {3.2.6}(iii), there \textit{might} also be a right-going line\footnote{Namely, the lower dashed/dotted blue line in Figure $1$.} with label $1$ having end-point in $R_{t''}\cap C_{r+1}$, but then its left end-point must lie in $R_{s'}\cap C_r(\infty):s'> s$, so by $(i)$, it is not joined by the special composite line to $C\setminus \{b_m\}$.

   \

   %\textbf{Remark}.  This ``might'' in $(*)$ and $(**)$ is illustrated by the blue lines in Figure $1$.  As indicated in the caption to Figure $1$ the existence of these lines is conditional.

   \

   This gives $(ii)$ with $(*),(**)$ giving $(iii)$.  Again Case one (resp. Case two) of the proof of $(ii)$ gives Case one (resp. Case two) of $(vi)$.

   $(iv)$ follows as in $(i)$ noting that the special composite lines cannot pass (a second time) strictly below  $R_s$.

   $(v)$. By Corollary \ref {3.2.6}(ii), no two distinct lines with label $1$ may have the same right end-point.  Then $(v)$ follows from $(i),(iv)$.

\end {proof}
\textbf{Remark.} After collapsing $\ell^\ast_{C,C'}$ joins uniquely determined boxes in the rectangle $R^s_{C,C'}$.

\subsubsection {Uniqueness}\label {5.2.3}

Again recall the \textbf{N.B.} of \ref {5.1.1}

\begin {cor} There is a unique disjoint union of composite lines in $R_{C,C'}^s$ passing though all the elements of $R_{C,C'}^s$ given that the only constituent line labelled by a $\ast$ is $\ell^\ast_{C,C'}$.
\end {cor}

\begin {proof}  Let us first show that one cannot find a disjoint union of composite lines going though all the boxes in $R^s_{C,C'}$, with no constituent line labelled by a $\ast$.  Such lines must be obtained by collapsing special composite lines $\mathscr L_i$ in $R^s_{C,C'}$.  As in  the proof of Proposition \ref {5.2.2}(i) these are uniquely determined by their starting points in $R^s\cap C$ and ending points in $R^s\cap C'$.   Disjointness would then force there to be a line with label $1$ having the \textit{same} starting and end points in $\mathscr T(\infty)$ as $\ell^\ast_{C,C'}$. This is excluded by Corollary  \ref {3.2.6}(i), or by inspecting Figure $1$.

Finally admitting that $\ell_{C,C'}^\ast$ does belong to the composite lines in the disjoint union and is the only one, forces the $\mathscr L_i:i\neq m$ to be the special composite lines of Proposition \ref {5.2.2} and then the assertion of the Corollary obtains from the conclusion of Proposition \ref {5.2.2}.
%Collapse all the properly disjoint composite lines $\mathscr T_i:i\in [1,s]$ given in Lemma \ref {5.2.2}.
%%The lines not starting from $b_m$ give composite lines in $\mathscr T$ from $C$ to $C'$ made up of lines labelled by a $1$.
%
%They are uniquely determined by their starting point in $C$ and properly disjoint.  Only one, namely $\mathscr L_m$ admits contains a line (namely, $\ell_{r,r'}$) with label $\ast$.
%
%Hence the assertion of the lemma.

%Recall the definition of $r$ and $j_r$ in Lemma \ref {3.2.2} and the definition of $t,r'$ in \ref {5.2.1}. Then $r$ is connected in $\mathscr T(\infty)$ by a vertical line labelled by a $\ast$ to the entry $r'$ of $R_t \cap C_{j_r+1}$ via \ref {3.2.3}.  Then by \ref {3.2.4} this becomes the unique line $\ell_{r,r'}$  joining $r,r'$ in $\mathscr T$ and is labelled by a $\ast$.  Then there is a composite line starting from $b|ast$ which is the concatenation of lines from $b_\ast$ to $r$ and from $r'$ to an entry of $C'$ all with label $1$ and the line $\ell_{r,r'}$ with label $\ast$.
%
%This composite line is disjoint from those considered in the first part by Lemma \ref {5.2.1}(iv).
%
%Finally by Lemma \ref {5.2.1}(iv) these composite lines from a disjoint union of $R_{C,C'}^s$ and there uniqueness implies in turn uniqueness of the composite lines in the conclusion of the lemma with the proviso that there is at most one line with label $\ast$.
\end {proof}

\subsubsection {Disjointness}\label {5.2.4}

We want to show that the above lemma remains valid if we omit the last condition, that is to say.

\begin {prop}  There is a unique disjoint union of composite lines in $R_{C,C'}^s$ passing though all the boxes of $R_{C,C'}^s$.
\end {prop}

\begin {proof} Retain the notation of \ref {5.2.2}.  The special composite lines $\mathscr L_i:i\in [1,s]$ are uniquely determined in $\mathscr T(\infty)$ up to $C_r$ and must be disjoint.  Hence they pass through all the boxes in $R^s \cap \mathscr T(\infty)$ in  $[C,C_r]$.  Now the lines in $\mathscr T(\infty)$ labelled by a $\ast$ are vertical.  Then disjointness successively precludes a box joining lines in $R^s \cap \mathscr T(\infty)\cap [C,C_r]$ to be labelled by a $\ast$.

If $\ell_{C,C'}^\ast$ \textit{is} one of the constituent lines of our unique disjoint union of composite lines, then in $\mathscr T(\infty)$, the subset $R^s \cap [C,C_{r+1}]$ can only be formed by the union of the  $\mathscr L_i:i\in [1,s]$ and then by the first paragraph this also holds in $R^s \cap [C,C']$.  Thus the only way the conclusion of the proposition can fail is for $\ell_{C,C'}^\ast$ to be absent.

Now any continuation of the composite lines $\mathscr L_i:i \neq m$ from $R^s \cap C_r(\infty)$ to $R^s \cap C_{r+1}(\infty)$ must leave empty one entry of $R^s \cap C_{r+1}(\infty)$.  Of course this must be filled with the end point of $\mathscr L_m$ which enters $C_r(\infty)$ in $R_{s'}:s'>s$. Let $\ell$ be the required additional constituent line from $R_{s'}\cap C_r(\infty)$.  It cannot be a line with label $1$ for such a line would then form the continuation of some $\mathscr L_i:i\neq m$.  Thus it can only be a line with a neutral label to $C_{r+1}(\infty)$ concatenated with a vertical line $\ell^\ast$ labelled by a $\ast$.

Recover the notation of the proof of (ii) of Proposition \ref {5.2.2}

 Suppose that $t-t'>0$, so the common entry of $\mathscr L_m$ goes down by $t+1-t'>1$ rows.  Then $\height C_{r+1}\geq s$ and by \ref {3.2.5} this precludes any vertical lines in the continuation of the $\mathscr L_i:i\neq m$ at $C_{r+1}$. Thus $\ell^\ast$ and $\ell^\ast_{C,C'}$  must have the same upper end points so coincide.

 Suppose that $t=t'=s$.  Then there is just one vertical line from $R_{s+1}\cap C_{r+1}(\infty)$ and so again $\ell^\ast$ and $\ell^\ast_{C,C'}$  must coincide.

 This completes the proof of the proposition.

\end {proof}

\subsubsection {Separation}\label {5.2.5}

Let $\mathscr P$ denote the set of pairs of neighbouring columns in $\mathscr D$.

Fix a component tableau $\mathscr T^\mathcal C$ and retain the notation of \ref {5.2.2}.

\begin {lemma}  The $\ell^\ast_{C,C'}: (C,C') \in \mathscr P$ are pairwise distinct.
\end {lemma}

\begin {proof}  Recall Remark \ref {5.2.2} and view  $\ell^\ast_{C,C'}$ as joining entries boxes of the rectangle $R^s_{C,C'}$ viewed as lying in $\mathscr T^\mathcal C$.

 %The proof is by induction on common height of a pair in $\mathscr P$.
 Take $i \in [1,2]$ and $(C_i,C'_i)\in  \mathscr P$ of height $s_i$.  If $s_1=s_2$, the assertion is trivial since  $[C_1,C_1'] \cap[C_2,C_2']$  is at most one column so there are no common lines at all joining boxes in $[C_1,C_1']$ and joining boxes in  $[C_2,C_2']$.

 Thus we can suppose $s_1>s_2$.

 Now view $\ell^\ast_{C,C'}$ as constructed in Proposition \ref {5.2.2}(ii).

Recall $(vi)$ of Proposition \ref {5.2.2} and cases one and two of its conclusion.

In the first case the right hand end-point of $\ell^\ast_{C_1,C_1'}$ meets $R_{s_1}$ (in $\mathscr T$), so it cannot equal $\ell^\ast_{C_2,C_2'}$ whose right hand end-point lies in $R^{s_2}$ (in $\mathscr T$).

In the second case the composite line in $\mathscr T(\infty)$ joining boxes of $C_1$ to those of $C_1'$ with constituent $\ell^\ast_{C_1,C'_1}$ stays in $R^{s_1}$ with the one exception that it may enter $R_{s_1+1}\cap C_{r+1}(\infty)$ and then re-enter $C_{r+1}$ (necessarily of height $<s$) at it unique lowest element by a vertical line labelled by a $\ast$.  Thus $s_1$ is determined by $\ell^\ast_{C_1,C'_1}$ and so again $\ell^\ast_{C_i,C_i'}:i=1,2$ are distinct.

%In the second case $s$ is determined by $\ell^\ast_{C,C'}$ and so again the $\ell^\ast_{C_i,C_i'}$ for $i=1,2$ are distinct.
\end {proof}

\textbf{Example $16$}.  Consider the composition $(2,1,1,1,2)$.  When $3$ is chosen from $\mathscr B^2_{C_1,C_5}$ of Example $15$ all the lines in $\mathscr T$ labelled by $\ast$ join boxes in $R_1$, so it is not so obvious that they are distinct, except when one draws the lines!  On the other hand in $\mathscr T(\infty)$ two of the lines with label $\ast$ join boxes in distinct columns of height $1$ in $R^2$ (which is the case $s_1=s_2=1$ in the proof) whilst the third meets $R_{s+1}$ in $\mathscr T(\infty)$, with $s$ being the common height of the pair $C_1,C_5$, namely $2$.

%In this case when one draws the lines in
%
% Yet $\mathscr T(\infty)$ two are in columns of the same height and join boxes in $R^2\cap \mathscr T(\infty)$ whilst the third meets $R^3\cap \mathscr T(\infty)$, so all three are distinct.  Of course when they are drawn it is perfectly obvious that they are distinct but one cannot draw  \textit{all} cases.

%\subsubsection {Separation}\label {5.2.5}
%
%Let $\mathscr P$ denote the set of pairs of neighbouring columns in $\mathscr D$.
%
%Fix a component tableau $\mathscr T^\mathcal C$ and retain the notation of \ref {5.2.2}.
%
%\begin {lemma}  The $\ell^\ast_{C,C'}: (C,C') \in \mathscr P$ are pairwise distinct.
%\end {lemma}
%
%\begin {proof}  The proof is by induction on common height. Take $i \in [1,2]$ and $(C_i,C'_i)\in  \mathscr P$ of height $s_i$.  If $s_1=s_2$, the assertion is trivial since there are no common lines at all joining boxes in $[C_1,C_1']$ and joining boxes in  $[C_2,C_2']$.  Thus we can suppose $s_1>s_2$.
%
%Recall $(vi)$ of Proposition \ref {5.2.2}.
%
%In the first case, $\ell^\ast_{C_1,C_1'}$ meets $R_{s_1}$ in $\mathscr T$, so it cannot equal $\ell^\ast_{C_2,C_2'}$ which belongs to $R^{s_2}$ in $\mathscr T$.
%
%In the second case $s$ is determined by $\ell^\ast_{C,C'}$ and so again the $\ell^\ast_{C_i,C_i'}$ for $i=1,2$ are distinct.
%\end {proof}

\subsubsection {Weierstrass Sections}\label {5.2.6}

It is clear that Proposition \ref {5.2.2} implies Proposition \ref {5.1.1}, whilst  \ref {5.1.2}$(*)$ is just the conclusion of Lemma \ref {5.2.5}.

Hence

\begin {thm}  Let $\mathscr T^\mathcal C$ be a component tableau.  Then $e_\mathcal C+V_\mathcal C$  is a Weierstrass section.  In particular
$\dim V_\mathcal C=\textbf{g}$.
\end {thm}

\section {Components} \label {6}

Throughout this section fix a component tableau $\mathscr T^{\mathcal C}$.

In \ref {4.1.5}, we described a subalgebra $\mathfrak u^\mathcal C$ spanned by positive root vectors defined by $\mathcal C$. It is complemented in $\mathfrak m$ by a set of positive root vectors whose roots are called the excluded roots $X$ (from $\mathfrak u^\mathcal C$).  Presented in $\mathfrak m$ (or in \textbf{M}) the excluded roots are encircled.

Recall that in \ref {4.3.1} we defined  the irreducible closed subvariety $\mathscr C =\overline{B.\mathfrak u^\mathcal C}$ of $\mathfrak m$. Our present goal is to show that it a component of $\mathscr N$.  By Corollary \ref {4.6} we already have the inclusion $\mathscr C \subset \mathscr N$.

%%The first step is to show that $\dim \mathscr C \geq \dim \mathfrak m - \textbf{g}$.
%%The method of proof is to  extend \cite [Lemma 6.9.6]{FJ2} to any component tableaux.
%\subsection {Inclusion}\label {6.1}
%
% \begin {lemma}  $\mathscr C \subset \mathscr N$.
% \end {lemma}
%
% \begin {proof}
%  Let $I$ be any Benlolo-Sanderson invariant.   Since it is $B$ invariant, it is enough to show that $I(\mathfrak u_{\min}^\mathscr C)=0$.  Yet by Theorem \ref {4.4}
%  it follows that $I$ vanishes on
%   $ \overline{B.\mathfrak  u^{\mathscr C}}$ Let $\mathscr T^\mathscr C$ be a component tableau.

%Given a conical  affine algebraic variety $\mathscr V$, let $[\mathscr V]$ denote its projectivisation.

\subsection {Using Intersection Theory in Projective Space}\label {6.1}

Assume that
$$ \dim \mathscr C \geq \dim \mathfrak m - \textbf{g}. \eqno{(*)}$$

Given a conical  affine algebraic variety $\mathscr V$, let $[\mathscr V]$ denote its projectivisation, that is to say suppress $0\in V$ and identify all multiplies of a given non-zero vector.

\begin {thm}  $\mathscr C$ is a component of $\mathscr N$ of dimension $\dim\mathfrak m -\textbf{g}$.
\end {thm}

\begin {proof}

We use the argument of \cite [6.9.8]{FJ2} which we sketch for completeness.

Let $\mathscr V^e$ denote the set of all non-zero scalar multiples of $e_\mathcal C+v:v \in V_\mathcal C\setminus \{0\}$.  It is a conical.

Let $\hat{\mathscr C}$ be an irreducible component of $\mathscr N$ containing $\mathscr C$. As noted \cite [6.9.8]{FJ2} $\hat{\mathscr C}$ is conical. By $(*)$, it suffices to show that $\dim \hat{\mathscr C}= \dim \mathfrak m-\textbf{g}$.

Assume that $\dim \hat{\mathscr C} > \dim \mathfrak m-\textbf{g}$, so then $\dim [\hat{\mathscr C}] >\dim [\mathfrak m]-\textbf{g}$.

By Theorem \ref {5.2.6} one has $\dim V_\mathcal C=\textbf{g}$ and so
$\dim [\mathscr V^e]=\textbf{g}-1$.

Then by the intersection theory of closed projective subvarieties in projective space \cite [Chap. 1, Section 5, Thm. 6] {Sh}, we conclude that $[\hat{\mathscr C}] \cap [\mathscr V^e]$ is non-empty.  Translated back to affine space this means that $\hat{\mathscr C}$ contains an element of  the form $e_\mathcal C+v:v \in V_\mathcal C \setminus \{0\}$.

Yet $e_\mathcal C +V_\mathcal C$ is a  Weierstrass section by Theorem \ref {5.2.6}, so there exists a Benlolo -Sanderson generator $I$ such that $I(e_\mathcal C+v)$ is a non-zero scalar. On the other hand $e_\mathcal C+v \in \hat{\mathscr C} \subset  \mathscr N$ and so $I(e_\mathcal C+v)=0$. this contradiction proves the assertion.

\end {proof}

 \subsection {Covering}\label {6.2}

 We establish \ref {6.1}$(*)$ by adapting the method used in \cite [6.9.6]{FJ2} to all component tableaux.

 \subsubsection {Excluded Roots not Labelled by a $\ast$}\label {6.2.1}

%Let $X$ be the set of excluded roots.  In \textbf{M} they are encircled.

Recall \ref {3.2.5}, \ref {3.2.6}.  Recall that $X$ is the set of excluded root vectors.  These are encircled.

Let $Y$ (resp. $S$) be the set of root vectors which carry a $\ast$ (resp. $1$).  Recall that $Y \subset X$  by Lemma \ref {4.1.3}.  Set $Z=X\setminus Y$.

%Let $S$ be the set of roots which are labelled by a $1$. By Lemma \ref {4.2.5}, they are never encircled.

Recall \ref {3.2.6} that these roots are constructed from adjacent columns $C_r,C_{r+1}$ with $\height C_{r+1}=t$ and an element $i\in R_{t'}\cap C_r(\infty)$ which has been lowered by $m:=t+1-t'\geq 1$ rows into $R_{t+1}\cap C_{r+1}(\infty)$, illustrated in Figure $1$ by a neutral line.

Recall the partial column $C_{r+1}^{>(t-m)}$ and $\textbf{j}$ its set of entries.  Denote $\textbf{j}$ simply by $j$ if the latter is a singleton. In the latter case $j$  is just  the lowest entry of $C_{r+1}$.

This gives a subset of $Y_{i,\textbf{j}}$ of roots labelled by a $\ast$, joining $i$ to $j \in \textbf{j}$.  (They are described in \ref {3.2.5} and recalled in the proof of the proposition below).  In \textbf{M} the corresponding root vectors appear in the last $m=|\textbf{j}|$ columns on the right  of the column block $\textbf{C}_{r+1}$, in its $i^{th}$ row.

The subset $X_{i,\textbf{j}}$ of excluded roots resulting from $Y_{i,\textbf{j}}$ (defined in \ref {4.1.4} and recalled briefly below) are called the excluded roots defined by the pair  $(i,\textbf{j})$.

 \subsubsection {Covering}\label {6.2.2}

Consider a root vector $x_{i',j'}$.  Here we are just are concerned with those in $\mathfrak m$ so we require $i'<j'$ and that $i',j'$ do not lie in the same column of $\mathscr T$.

Given root vectors $x_{i',j'},x_{k',l'}$ in $\mathfrak m$, we say that $x_{i',j'}$ covers $x_{k',l'}$ if $i'=k'$ and $j'<l'$, that is they lie in the same row of \textbf{M} with the former to the left.

The set of root vectors which cover $Z_{i,\textbf{j}}:=X_{i,\textbf{j}}\setminus Y_{i,\textbf{j}}$ is denoted by $S_{i,\textbf{j}}$. The $Z_{i,\textbf{j}}$ (resp. $S_{i,\textbf{j}}$) are not necessarily disjoint for distinct pairs.

\subsubsection {Recalling the Construction of $X_{i,\textbf{j}}, Y_{i,\textbf{j}}$}\label {6.2.3}

\

Case $1$.  Here $m= t+1-t'>1$ and $\height C_{r+1}=t$.  The adjacent pair must be surrounded by $m$ (free) pairs of neighbouring columns of heights $k\in [t',t]$.  This gives $m$ vertical lines labelled by a $\ast$ from $R_{t+1}\cap C_{r+1}(\infty)$ with entry $i$, to each of the boxes of $C_{r+1}^{>(t-m)}$, that is the $R_{t-k+1}\cap C_{r+1}$, with entry $j_k=j-k+1 : k \in [1,m]$.

The entries joined by these vertical lines define $Y_{(i,\textbf{j})}$.

Following this recall \ref {4.1.3} that $C_{r+1}^{>(t-m)}$ is placed below $i$ in the unique column $C_h:=C_{(i)}$ in which $i$ appears in $\mathscr T$.  Let $R_f:=R_{(i)}$ be the unique row in which $i$ appears in $\mathscr T$.  Then $\height C_h \geq f$.  If this is a strict inequality the partial column $C_h^{>f}$ is moved leftwards displacing similar partial columns until the first column $C_g$ of height $f$ is reached, as detailed in \ref {4.1.3}.

\

Case $2$.  Here $m=1$ and $\height C_{r+1}\leq t$.  It is brought to height $t$ in $\mathscr T(t)$ - see \ref {3.1.1}, \ref {3.2.2}.  The adjacent pair must be surrounded by a (free) pair of neighbouring columns of height $t$.  This gives a  vertical line labelled by a $\ast$ from $R_{t+1}\cap C_{r+1}(\infty)$ with entry $i$, to the lowest box of $C_{r+1}$ with entry $j$.

In this case $Y_{i,j}$ is the singleton defined by $\ell_{i,j}$.

%However unlike the previous case (taking $m=1$), this box may lie (even strictly) above $R_t\cap C_{r+1}$.  Let its entry be $j$.

Following this, as in case $1$, the entry  $j$ is placed below $i$ in $C_h:=C_{(i)}$ and similarly partial columns are moved leftwards.

\subsubsection {The Covering Proposition-Primary Excluded Roots}\label {6.2.4}

\

Recall the notation of \ref {4.1.3}. Consider the left movement of partial columns, described in \ref {4.1.3} and referred to in \ref {6.2.3}.

Given a column $C$ between $C_g$ and $C_h$, let $\tilde{C}^{>f}$ be the partial column that replaces $C^{>f}$ in the said leftward movement of partial columns.  We call $\tilde{C}$ the displacing column.

Both cases described in \ref {6.2.3} are the same except for the specification of the first displacing column $\tilde{C_h}$, which in the first case is $C_{r+1}$.  In both cases  $x_{k',l'}\in X_{i,\textbf{j}}$ is obtained by taking $l'\in \tilde{C}^{>f}$ and  $k'\in [C^{\leq f},\tilde{C}]$ such that $k'<l'$.

When $C=C_h=C_{(i)}$ (resp. when $C$ lies strictly to the left of $C_{(i)}=C_h$) these roots were referred to \ref {4.1.4} as the primary $Z^P_{i,\textbf{j}}$ (resp. secondary $Z^S_{i,\textbf{j}}$ ) excluded roots in $Z_{i,\textbf{j}}$.

In the notation of \ref {4.1.4}, let $Z^P=\cup_{i,\textbf{j}\in \textbf{L}} Z^P_{i,\textbf{j}}$ (resp. $Z^S=\cup_{i,\textbf{j}\in \textbf{L}} Z^S_{i,\textbf{j}}$) denote the subset of $Z$ of primary (resp. secondary) excluded roots.

\begin {prop}  The roots of $S$ cover  $Z^P$.
\end {prop}

\begin {proof}

\

It is enough to prove this assertion for each choice of $i$ defining a $Z^P_{i,\textbf{j}}$, in the notation and hypotheses of \ref {6.2.3}.
Recall (\ref {3.2.5})  that in this $i$ moves down from $C_r(\infty)$ to $C_{r+1}(\infty)$ by $m\geq 1$ steps.
Fix $x_{i',l'}\in Z^P_{i,\textbf{j}}$.

%Retain the hypothesis and notation of \ref {6.2.3}
%.
%Consider a root $x_{i',j'}$ with $i'=k'$.  These are the roots given by a line whose left end-point lies in $[C,\tilde{C}]$.
%
%% $C_{(k')}\setminus C_{(k')}^{>(i)}$  or in a column lying $]C_{(k')},C_{(l')}[$.
%
%%$C\setminus C^{>}$ or in a columns strictly to the right of $C$.
%

%Recall \ref {4.2.1} the notion of an $s$-string and that distinct $s$-strings do not cross, Lemma \ref {4.2.1}.
%
%Consider the $i'$-string defined by $i'$ in $\mathscr T(\infty)$.  By \ref {3.2.2} it moves to the right horizontally or to the right descending by $\geq 1$ rows at a pair of adjacent columns.  In the latter case for some $j'>i'$, there exists an up-going vertical line $\ell_{i',j'}$ in the right hand column of the adjacent pair, labelled by a $\ast$. This means that $x_{i',j'} \in Y$.
%
% Otherwise the $i'$-string does not descend.
%
% Consider covering for the primary excluded roots.
%
% The roots that we do not need to cover are those of $Y$ - not necessarily just those of $Y_{(i,\textbf{j})}$.

 Consider $i\in R_f$, in the notation of \ref {6.2.1}.  Take $C=C_{(i)}=C_h$. The $i$-string must descend strictly (see proof of Claim in \ref {4.1.3}) and perhaps in  several steps to reach $C_{r+1}(\infty)$.  Moreover at each downward step it always enters an \textit{empty} box in $[C,\tilde{C}]$ - as specified in \ref {3.2.2}.

   By Observation $1$, \ref {4.1.5},  $i'$ (as an element of $\mathscr T$) lies in $C$ strictly above $i$ or lies in a column strictly to the right of $C$ and to the left of $\tilde{C}=C_{r+1}$.

View $\mathscr T$ as a sub-tableau of $\mathscr T(\infty)$.
Then by the two paragraphs above, $i'\in \mathscr T$  must be distinct from $i$ and lie
   strictly above $i\in \mathscr T(\infty)$ in the column of $[C_{(i)},C_{r+1}]$ in which they both appear. Then by Lemma \ref {4.2.1}(ii) the $i'$-string lies strictly above the $i$-string, in any column of $\mathscr T(\infty)$ through which they both pass.

    This condition on strings puts $i$ strictly below $i'$ in $C_r(\infty)$.  Then if
    $i'$ is stopped at $C_r$, there is a line $\ell_{i',j'}$ labelled by a $1$ with $j' <l'$ by Corollary \ref {3.2.6}(iii).  In this case  $x_{i',j'}\in S$ and covers $x_{i',l'}\in Z$. If not by Figure $1$, \ref {3.2.6}, this can only happen if $m=1$ (with $m$ as defined in the first paragraph of the present proof) and then $x_{i',j'}\in Y$.

    \end {proof}

    \subsubsection {The Covering of Secondary Excluded Roots}\label {6.2.5}

     None of the secondary excluded roots lie in $Y$, so the argument is less convoluted.   Again we need only consider the covering of some $Z^S_{i,\textbf{j}}$, and we drop the subscript.

 For a secondary excluded root $x_{k',l'}$  the displacing column $\tilde{C}$ is the rightmost column in $C_g,C_h$ of height $>f$.  Thus  $l'\in \tilde{C}^{>f}$, and  $k'\in  [C,\tilde{C}]$, with of course $k'<l'$.

 %Then one easily checks that $k'<l'$ if and only if $k' \in R^f$.

 % One may remark that in both cases we may compute $Z_{(i,r)}$ and $S_{(i,r)}$, as shown in the examples.

 Recall the notation of \ref {4.1.3}.  Recall that a subset $Z^S$ of secondary excluded roots occur when partial columns are displaced to the left starting from $C_{(i)}$ and ending in the left hand column $C_g$ of a pair of columns $C_g,C'_g$ of height $f$ surrounding $C_{(i)}$.  In this $\height C_{(i)} \geq f$ and moreover a strict inequality is needed otherwise $Z^S$ is empty.

  Recall that in the notation of \ref {4.1.3} the co-ordinates of the secondary excluded roots lie in $[C_{g_{k-1}},C_{g_{k}}]:k \in [1,s]$ with $g_0=g$. Moreover the columns in $]C_{g_{k-1}},C_{g_k}[$ all have height $<f$, by choice of notation.

     \subsubsection {A Key Lemma}\label {6.2.6}

   In above fix $i\in [1,s]$ and  set $C_{k_1}=C_{g_{i-1}},C_{k_2}=C_{g_i}$.  They are columns of height $\geq f$ with at most $C_{g_0}$ of height $s$. Strictly between them there are only columns of height $<f$.

   Here the entries of the $C_l:l\in [k_1,k_2]$ are those given by the numbering of $\mathscr T$ and the lines between boxes in $[C_{k_1},C_{k_2}]$ are those coming from the given component tableau $\mathscr T^\mathscr C$ chosen.

   Observe that

  \

  $(*)$.  A neutral line from $R_{t_1}:t_1\leq f$ to $R_{t_2}:t_2>f$ must have at least one entry outside $[C_{k_1},C_{k_2}]$, so a neutral line cannot penetrate $R_s$ in $[C_{k_1},C_{k_2}]$.

  \

  This is because the pair of neighbouring columns $C_g,C'_g$ of height $f$ is used to lower $i\in C{(i)}=C_{g_s}$ strictly below $R_f$, as noted in the paragraph following the claim in \ref {4.1.3} and so this pair is not free to allow a neutral line to penetrate $R_s$ in $[C_{k_1},C_{k_2}]$.

 \begin {lemma} There is a disjoint union of special composite lines going through all the boxes of $R^s\cap [C_{k_1},C_{k_2}]$ in $\mathscr T(\infty)$.
 \end {lemma}

 \begin {proof}

 The proof is similar to that of Proposition \ref{5.2.2}(i).

 %Set $C=C_{k_1},C'=C_{k_2}$.

 For all $t \in[1,s],l\in [k_1,k_2-1]$, the unique $u$-string in $\mathscr T(\infty)$ defined by the entry of $R_t\cap C_l(\infty)$, either enters $R^s\cap C_{l+1}(\infty)$, or is stopped at $C_l$ defining a unique line with label $1$ to $R^s\cap C_{l+1}$ having right end-point $u'$.  Indeed the third possibility arising in Proposition \ref  {5.2.2}(i) of a neutral line penetrating $R_s$ in $[C_{k_1},C_{k_2}]$ is excluded by $(*)$.

  As in Proposition \ref {5.2.2}, the existence and uniqueness of the evolution of the special composite lines starting from  boxes in $R^f\cap C_{k_1}$ and ending in  $R^f\cap C_{k_2}$, proves the lemma.

 \end {proof}

 \subsubsection {The Covering Roots}\label {6.2.7}

 For all $k \in [1,s]$, let $Z^S(k)$ denote the subset $Z^S$ coming from the pair $C_{g_{k-1}},C_{g_{k}}$.
 By \ref {4.1.4}, Observation $3$, the first co-ordinate $k'$ of $x_{k',l'}$ lies in $R^f \cap [C_{g_{k-1}}^{\leq f},C_{g_{k}}[$.

  %Precisely the first co-ordinate $k'$  of $x_{k',l'}$ lies in $[C_{g_{k-1}}^{\leq f},C_{g_{k}}]$ and the second co-ordinate $l'$ in $C_{g_k}^{>f}$.  Those with first co-ordinate in $C_{g_k}$ belong to the Levi factor and can be ignored.  Since the columns strictly in $]C_{k_1},C_{k_2}[$ are of height $<f$, we obtain $k'\in R^f \cap [C_{g_{k-1}}^{\leq f},C_{g_{k}}[$.

 Let $S(k)$ be the subset of $S$ obtained by collapsing the disjoint composite lines of Lemma \ref {6.2.6}.

 %Through the lemma we conclude that for each box of $b'_v:=R_v\cap C_{g_k}:v \in [1,f]$ there exists a left  going line with label $1$ from $b'_v$ a box $b_u$ in $R^f\cap [C_{g_{k-1}},C_{g_k}[$. Let $k'$ (resp. $l'$) be the entry of $b_u$ (resp. $b'_v$) and denote this line by $\ell_{k',l'}$.  It defines an element $x_{k',l'}\in S'_k$.

 \begin {lemma}  The elements of $S$ cover those of $Z^S$.
 \end {lemma}

 \begin {proof} It is enough to show that the elements of $S(k)$ cover those of $Z^S(k)$ for all $k \in [1,s]$. Through Lemma \ref {6.2.6} we conclude that to a box in $R^f \cap [C_{g_{k-1}},C_{g_k}[$, with label  $i'$ there exists a strictly right going line $\ell_{i',j'}\in S'_k$ with label $1$ to a box in $R^f \cap [C_{g_{k-1}},C_{g_k}]$ with label $j'$, with necessarily $j'>i'$.

 On the other hand the elements of $Z^S(k)$ take the form $\ell_{k',l'}$ with first co-ordinate in $[C_{g_{k-1}},C_{g_{k}}[\cap R^f$ and second co-ordinate in $C_{g_k}^{>f}$. Thus $l'>j'$, with $l'$ defined in the first paragraph. Since $i'$ in the first paragraph is an arbitrary element of $[C_{g_{k-1}},C_{g_{k}}[\cap R^f$ we can take $i'=k'$ and then $x_{i',j'}\in S'_k$ covers $x_{k',l'}\in Z'_k$, as required.
 \end {proof}

 \subsection {Tangent Spaces}\label {6.3}

 We now generalize the claim in \cite [6.9.7]{FJ2} to any component tableau $\mathscr T^\mathcal C$, to obtain a proof of \ref {6.1}$(*)$.

  Let \textbf{S},\textbf{X},\textbf{Y},\textbf{Z} be the linear span of the $x_{k,l} \in S,X,Y,Z$, respectively.  Then by definition of $\mathfrak u,e_\mathcal C$ and the latter subspaces we have
  $$e_\mathcal C=\sum_{\alpha \in S}x_\alpha \in\mathfrak u, \quad \mathfrak m =\mathfrak u \oplus \textbf{X}, \quad \textbf{X}=\textbf{Z} \oplus \textbf{Y}.  \eqno{(1)}$$

   In addition  $V_\mathcal C=\textbf{Y}$.  Then by Theorem \ref {5.2.6}, $e_\mathcal C+V_\mathcal C$ is a Weierstrass section in $\mathfrak m$. In particular $\dim \textbf{Y}$ equals the number \textbf{g} of Benlolo-Sanderson generators of the polynomial ring $\mathbb C[\mathfrak m]^{P'}$.

 \begin {prop}

 \

 $(i)$. $(\mathfrak u+\mathfrak n.e)+ \textbf{Y}=\mathfrak m$.

 \

 $(ii)$. $\mathfrak n.e \cap \textbf{Y}=\{0\}$.

 \

 $(iii)$. $\dim B.\mathfrak u= \dim (\mathfrak u +\mathfrak n.e) = \dim \mathfrak m -\textbf{g}$.

 \

  $(iv)$. $(\mathfrak u +\mathfrak n.e)\oplus \textbf{Y}=\mathfrak m$.

  \

  $(v)$.  $\mathscr C =\overline {B.\mathfrak u}$ is an irreducible component of $\mathscr N$.
 \end {prop}

 \begin {proof} $(i)$.  Take $(k,l) \in Z$.  By the Proposition \ref {6.2.4} there exists $(i,j)\in S$ with $i=k,j<l$. Moreover by Corollary {3.2.6}(ii) this element of $S$ is uniquely determined by $i$.  Again since $j<l$ one has $x_{j,l}\in \mathfrak n$.  Then $x_{i,l}\in \mathfrak n.x_{i,j}$.

 Just as in \cite [6.9.7]{FJ2} this is not quite the end of the story.  Indeed there can be a second pair $(l,m) \in S$, which if it exists is uniquely determined by $l$.  Then $[x_{j,l},x_{l,m}]= x_{j,m}$, whilst the asserted uniqueness gives
 $$x_{j,l}.e= x_{j,m}-x_{i,l}.\eqno{(2)}$$

 Moreover generally for \textit{any} root vector $x \in \mathfrak n$ the expression $x.e$ is at most the sum of two root vectors\footnote{Indeed this is an advantage of using $e$ whilst for a generic element $\tilde{e}\in \mathfrak u$, there are many terms in $\mathfrak n.\tilde{e}$.}.

 In the above we note that $l<m$, so in \textbf{M}  and the two terms in the right hand side of $(2)$ must belong to distinct columns blocks of \textbf{M} with the first term in a column block strictly to the right of that containing the second term. In particular if the  first term lies in  the last column block, then the second term does not appear.

 Now $\mathfrak u,\textbf{Y},\mathfrak m$ are all direct sums of their respective intersections with the column blocks, it follows from $(1)$ that the left hand side of $(1)$ contains the entire last column block of \textbf{M}.

 Then correspondingly $(i)$ obtains by right to left induction on column blocks.

 $(ii)$ is essentially \cite [Lemma 3.1(ii)]{FJ2} but due to its possible interest we give a more detailed proof.

 First given a finite dimensional vector space $\mathfrak m$ and $c$ an indeterminate, use $O(c^2)$ to denote a sum of terms which multiply $c^m:m\geq 2$.

 Now take $p \in \mathbb C[\mathfrak m]$, that is to say a polynomial function on $\mathfrak m$, which more particularly is  a sum of terms of positive degree. Then given $e,v \in \mathfrak m$, one has
 $$p(e+cv +O(c^2))=p(e+cv) +O(c^2). \eqno {(3)}$$

 It is enough to prove $(3)$ with $p$ homogeneous of degree $m \geq 1$, that is a product of $m$ linear functions $f_i:i\in [1,m]$.  Then the left hand side above equals $\prod_{i=1}^m f_i(e+cv)+O(c^2)=p(e+cv)+O(c^2)$.  Hence $(3)$.

 Now recall that $\mathfrak n$ acts locally nilpotently on $\mathfrak m$, so $(\exp cx)e$ is  a polynomial in $c$ with coefficients in $\mathfrak m$.

 Now take $e=e_\mathcal C,V=V_\mathcal C$.  Suppose $(ii)$ fails.  Then there exists $x \in \mathfrak n$ such that $x.e:=v \in V\setminus \{0\}$.  Choose a linear function $\xi$ on $V$ such that $\xi(v)\neq 0$.  By definition of a Weierstrass section, there exists a sum of Benlolo-Sanderson invariants $p$ (hence a polynomial $p$ which is a sum of terms of positive degree) such that $p(e+cv)=c\xi(v)\neq 0$ for all $c\neq 0$.

 The invariance of $p$ implies that $p(\exp(cx).e)=p(e)=0$. Yet the left hand side equals $p(e+cx.e +O(c^2)=c\xi(v)+O(c^2))$.  The latter has derivative at $c=0$ equal to $\xi(v)\neq 0$ contradicting the previous assertion.  Hence $(ii)$.

 Let $N$ be the closed irreducible subgroup of $B$ with Lie algebra $\mathfrak n$.  One has $B=NH$ and since $\mathfrak u$ is $H$ stable we obtain $B.\mathfrak u=N.\mathfrak u$.

 %We may obtain a lower bound on $\dim N.\mathfrak u$ by considering the action of $N$ on a (generic)

 Take a (generic) point $\tilde{e} \in \mathfrak u$. The action $N.\mathfrak u$ of $N$ on $\mathfrak u$ augments its dimension by the dimension of the tangent space through $\tilde{e}$ modulo $\mathfrak u$, in brief from $\dim \mathfrak u$ to $\dim (\mathfrak u +\mathfrak n.\tilde{e})$.  Just taking $\tilde{e}=e$ the latter dimension is already $\geq \dim \mathfrak m - \dim \textbf{Y}=\dim \mathfrak m - \textbf{g}$, by $(i)$.

  Thus $\dim B.\mathfrak u \geq \dim \mathfrak u +\mathfrak n.e \geq  \dim \mathfrak m - \textbf{g}$ and so by Theorem \ref {6.1} equality holds throughout.  Hence $(iii)$.

% Now those $X.e:X \in \mathfrak u$ which are linearly independent modulo $\mathfrak u$ form one-parameter subgroups of $N$ making  $N.\mathfrak u$ a fibre bundle over $\mathfrak u$ of dimension $\dim (\mathfrak u +\mathfrak n.e)\geq \dim \mathfrak m -g$, by (i).  (\mathfrak u +\mathfrak n.e) equality holds, which in view of (i) gives (iii).

 %$(iii)$ follows from $(i)$ and $(ii)$.  Indeed the vector spaces in $(1)$ are spanned by root vectors, whilst the right hand side of $(2)$ is a sum of just two root vectors one already in \textbf{Z}.

 $(iv)$ follows from $(iii)$ and $(i)$.  $(v)$ follows from $(iii)$ and Theorem \ref {6.1}.
 \end {proof}

 \textbf{Remark}.  it is possible to show \textit{combinatorially} that $(ii)$ implies $(iii)$. However perhaps the reader prefers our present use of intersection theory in projective space, which besides also shows (Theorem \ref {6.1} and Proposition \ref {6.3}) that $\mathscr C$ is a component of $\mathscr N$.

\subsection {The Component Map}\label {6.4}

By Proposition  \ref {6.3} there is map $\mathscr T^\mathcal C \mapsto \mathscr C$ of the set of component tableaux to the set of irreducible components of $\mathscr N$.  We call it the component map.

\section {Injectivity of the Component Map} \label {7}

Fix a composition hence a tableau $\mathscr T$,  and a pair of \textit{distinct} component tableaux $\mathscr T^{\mathcal C},\mathscr T^{\mathcal C'}$ associated to this composition.

\subsection {Choices in Batches.} \label {7.1}  Recall the notion of a batch \ref {3.1.4}, and the notation given there.

\begin {lemma}  There exists $s$ minimal and then $u$ minimal with the following property.  There exist $i,i' \in \mathscr B^s_u$ \textit{distinct} with $i\in \mathscr B^s_u$ (resp. $i' \in \mathscr B^s_u$) defining the construction of $\mathscr T^\mathcal C$  (resp. $\mathscr T^{\mathcal C'}$).

\end {lemma}

\begin {proof} Since entries in batches are determined with $s$ increasing and that the order of appearance in a given $\mathscr B^s_u:u\in [1,r_{s-1}]$ is immaterial, it would otherwise be true that $\mathscr T^\mathcal C=\mathscr T^{\mathcal C'}$.
\end {proof}

\subsection {An Exchange Lemma.} \label {7.2}

Let $C,C'$ be the neighbouring columns of height $s$ defining $\mathscr B^s_u$ and set $I=I^s_{C,C'}$.

Then the rectangle $R^s_{C,C'}$ admits $i:=i^\mathcal C_{C,C'}$ (resp. $i':=i^{\mathcal C'}_{C,C'}$) as the penetrating $i$-string (resp. $i'$-string) for $\mathcal C$ (resp. $\mathcal C'$).

%Fix a composition hence a tableau $\mathscr T$,  and a pair of distinct component tableaux $\mathscr T^{\mathcal C},\mathscr T^{\mathcal C'}$ associated to this composition.

  \begin {lemma} There is a line $\ell$ in $\mathscr T^{\mathcal C}$ labelled by a $\ast$, which in $\mathscr T^{\mathcal C'}$ is labelled by a $1$ and vice-versa.
  \end {lemma}

  \begin {proof}

  By interchanging labels on tableaux, we can assume, so by convention, that $i'>i$.

  Since $i\in \mathscr B^s_u$ arises for $\mathscr T^\mathcal C$, there exists $r \in [1,k]$ so that by Rule $(1)$ of \ref {3.2.2} applied to  $\mathscr T^{\mathcal C}(\infty)$, $i$ goes from $R_m\cap C_r(\infty)$ down by $m'\geq 1$ rows into $R_{m+m'}\cap C_{r+1}(\infty)$, by a neutral line and then by possibly by several lines labelled by a $\ast$ each to a different entry of $C_{r+1}$.

  This means in particular that $\height C_{r+1}=m+m'-1$.  Thus by Rule $(2)$ of \ref {3.2.2} applied to  $\mathscr T^{\mathcal C'}(\infty)$, the entry $i$ is stopped at $C_r$, giving a line $\ell$ with label $1$.

  We claim that the entry of right end point one of the of the lines, say $\ell'$,  with label $\ast$ in $\mathscr T^{\mathcal C}(\infty)$ can be chosen to be the same as the right end point of the line $\ell$  with label $1$ in $\mathscr T^{\mathcal C'}(\infty)$.

  The assertion is trivial if $m'=1$ because both lines necessarily go to the unique lowest entry of $C_{r+1}$.

  Otherwise take $m=m_1+1, m'=m_2$ following assiduously Figure $1$.  The line $\ell$ with label $1$ in  $\mathscr T^{\mathcal C'}(\infty)$ is horizontal and has right end-point in $R_{m_1+1}\cap C_{r+1}$, whilst in $\mathscr T^{\mathcal C}(\infty)$ there is a unique (vertical) line with label $\ast$ to $R_{m_1+1}\cap C_{r+1}$.  This gives the required assertion.

   Yet this is not quite the end of the story.  Indeed the proof makes a choice of component tableau being defined by the choice of $i \in \mathscr B^m_u$.

  Yet the argument of the first part equally well applies to the second part when  $i'>i$, the only difference is that $i'$ goes down by $\geq 1$ rows from $C_{r'}(\infty)$ into $C_{r'+1}(\infty)$ with $r'>r$.

\end {proof}

\textbf{Example $17$}.  Consider the composition $(2,1,1,2)$.  This admits two component tableaux. In both cases $3$ enters $\mathscr T(\infty)$ below $4$.  Then the batch $\mathscr B_{C_1,C_4}^2$ has two elements $\{2,3\}$.  Let $\mathcal C'$ (resp.  $\mathcal C$) be that in which $i'=3$ is chosen (and placed below $6$) (resp. $i=2$) is chosen and placed below $3$ in $C_3$. Then in $\mathscr T^{\mathcal C'}$ (resp. $\mathscr T^{\mathcal C}$) the line $\ell_{2,4}$ has label $1$ (resp. $\ast$) and comes from the first part of the proof, whilst for $\ell_{3,6}$, the labelling is reversed, as comes from the second part of the proof.
Further examples arises from Figures $18$-$20$ where the two root vectors $l,l'$ have been marked in red assuming the printer courteous enough to provide colour!  Colour in Figure $1$ is not so vital since the different lines are otherwise identified (see Caption to Figure $1$).  The colour scheme in the remaining figures, indicate repeated entries according to the rules of \ref {3.2.2}, is more aesthetic than essential.

\subsection {The Rightmost Line.} \label {7.3}

The lines $\ell:=\ell_{i,j},\ell':=\ell_{i',j'}$ in the conclusion of the exchange lemma are to be denoted (as usual) by the entries they join.  Notice that in this $i\in C_r(\infty),i'\in C_{r'}(\infty)$ whilst $j\in C_{r+1},j'\in C_{r'+1}$. In particular the right end point of $\ell'$ lies strictly to the right of that of $\ell$.   We call the former the rightmost line (for the given pair of component tableaux).  In Example $17$ above $r=2,r'=3$.

\

\textbf{N.B.} By our convention above, the rightmost line will be labelled by a $1$ in $\mathscr T^\mathcal C$.

\

Let $x_\ell$ (resp. $x_{\ell'}$) be the co-ordinate vector in \textbf{M} defined by $\ell$ (resp. $\ell'$).

Let $X(i^\mathcal C_{C,C'})$, be the set of excluded roots vectors constructed as in \ref {4.1.4} \textit{but only} through the lines with label $\ast$ defined by the steps in the penetrating $i=i^\mathcal C_{C,C'}$-string (recall \ref {4.3.5}) until just after this $i$-string leaves the rectangle $R^s_{C,C'}$,
at which point we say \textit {by convention} that the $i$-string is halted, though it may continue! (see Example $19$).

Notice that $X(i^\mathcal C_{C,C'})\subset X$.   This inclusion may be strict and is key in proving the lemma below, through the following observation.

\

$(*)$. The second co-ordinate $l$ of the root vector $x_{k,l}\in X(i^\mathcal C_{C,C'})$  lies to the left of $C_{r+1}$ (but not necessarily strictly).

\

\textbf{Definition}. Let $\mathfrak u^\mathcal C_{C,C'}$ be the complement in $\mathfrak m$ of the space spanned by the set $X(i^\mathcal C_{C,C'})$.

\

%Let $\mathfrak u^\mathcal C_{C,C'}$ be the complement in $\mathfrak m$ of the excluded roots created by the steps in the penetrating  $i=i^\mathcal C_{C,C'}$-string until just after it leaves $R^s_{C,C'}$.

 One has $\mathfrak u^\mathcal C \subset \mathfrak u^\mathcal C_{C,C'}$,  because the right hand side has less excluded vectors. By Proposition \ref {4.4} we have $I(\mathfrak u^\mathcal C_{C,C'})=0$. (A key fact in what follows!)

%This also holds with $\mathcal C$ replaced by $\mathcal C'$, but it will not be needed.

Recall that $x_{\ell'}=x_{i',j'}$.

\begin {lemma} $Bx_{\ell'}\subset \mathfrak u^\mathcal C_{C,C'}$.
\end {lemma}

\begin {proof}  %We must show that there is no excluded root $x_{k,l}$ coming from any line labelled by a $\ast$ from penetrating $i$-string in the upper right-hand quadrant \textbf{Q} of \textbf{M} above $x_{i',j'}$.

Let \textbf{Q} be the set of co-ordinate vectors in the upper right-hand quadrant of \textbf{M} with lower left hand corner at $x_{\ell'}=x_{i',j'}$.  Obviously $B.x_{\ell'}$ is contained in the linear span of $\textbf{Q}$ (but note the former is \textit{not itself} a vector subspace).

 Thus we must show that no root vector $x \in \textbf{Q}$ lies in $X(i^\mathcal C_{C,C'})$.
 %is an excluded root vector coming from any line labelled by a $\ast$ arising from the steps of the penetrating $i$-string, until it stops \textit{in the sense of \ref {7.3}}.

 \

 $(**)$. One has $j' \in C_{r'+1}$ which lies strictly to the right of $C_{r+1}$.

 \

   One has $x_{i',j'}\in \textbf{Q}$, but since it is labelled by a $1$ in $\mathscr T^\mathcal C$ so is not an excluded vector (Prop. \ref {4.2.5}) and in particular does not lie in $X(i^\mathcal C_{C,C'})$.

   Take $x_{k,l}\in \textbf{Q}$.

  % Then $l \geq j'$, or $k \leq i'$.

  Suppose $k\leq i',l > j'$.  Then $x_{k,l}\notin X(i^\mathcal C_{C,C'})$, since by $(**)$ this would force a line labelled by a $\ast$ coming from the penetrating $i$-string to have its second co-ordinate $l$ lying to the right of $C_{r'+1}$, so strictly to right of $C_{r+1}$. This is excluded by $(*)$.  It is illustrated in Example $19$, Figure $14$.

   Suppose $k<i', l=j'$.   Then $x_{k,l}\notin X(i^\mathcal C_{C,C'})$, as this would have to be created by partial columns shifted to the left starting from a partial column to the right of $C_{r'+1}$ itself lying strictly to the right of $C_{r+1}$. This is again excluded by $(**)$. It is illustrated in Examples $18,20$, Figures $13,15$.  The latter is slightly sharper since the shifting starts already at $C_{r'+1}$, whereas in the former it starts at $C_{r'+2}$.

\end {proof}

\subsection {Partial Linearity for the Rightmost Line.} \label {7.4}

Recall the notation of \ref {7.2}.

\begin {cor} $ I(x_{\ell'}+\overline{B.\mathfrak u^{\mathcal C}})=0$.
\end {cor}

\begin {proof} Since $x_{\ell'}$ is a point and $I$ is a regular function, we can remove closure in the statement of the lemma. Consider $x_{\ell'}+b.u'$ with $b \in B.u' \in \mathfrak u^\mathcal C$.  Since $b^{-1}x_{\ell'}=u'' \in \mathfrak u^\mathcal C_{C,C'}$, by Lemma \ref {7.3}, we obtain $x_{\ell'}+u'=b(b^{-1}x_{\ell'}+u')=b(u''+u')\in B.\mathfrak u^\mathcal C_{C,C'}$, by the linearity of $\mathfrak u^\mathcal C_{C,C'}$.  Then $I(x_{\ell'}+B.\mathfrak u^{\mathcal C})\subset I(B.\mathfrak u^{\mathcal C}_{C,C'})=0$, by Theorem \ref {4.4}.
\end {proof}
\textbf{Remark.} Even under the hypothesis of the lemma, it can be false that  $B.x_{\ell'} \subset \mathfrak u^\mathcal C$.  Thus we may need to restrict co-ordinates following $I$.  For this see Examples $18-20$.

\subsection {Injectivity}\label {7.5}

Retain the notation of \ref {7.2}.  Let $\mathcal C,\mathcal C'$ be distinct numerical data associated to a fixed composition.

 \begin {thm} Let $\mathcal C,\mathcal C'$ be distinct numerical data associated to a fixed composition.  Then $\overline {B.\mathfrak u^\mathcal C} \neq \overline {B.\mathfrak u^{\mathcal C'}} $.
 \end {thm}

 \begin {proof}  Define $\ell^*_{C,C'}$, introduced in Proposition \ref {5.1.1}(ii).

 We claim that $\ell'=\ell^*_{C,C'}$.  The choice of $\ell^*_{C,C'}$ is made in Proposition \ref {5.2.2}(ii).  Here its left end-point is $i=i^\mathcal C_{C,C'}$.  If \textbf{Case two} of the proof applies its right end-point is uniquely determined making the claim trivial.  Suppose \textbf{Case one} holds, then its right-end point is uniquely determined by the condition that its right end point lies in $R_s$. Since its left end-point $i\in \mathscr B^s_{C,C'}\subset R_s$.  This forces it to coincide with the horizontal line $\ell$ labelled by a $1$ in  $\mathscr T^{\mathcal C'}$ and hence with $\ell'$ labelled by a $\ast$ in $\mathscr T^{\mathcal C}$.

 By definition of a Weierstrass section, this  means that the restriction of $I=I^s_{C,C'}$ to $e_{\mathcal C'}+V_{\mathcal C'}$ is just the linear function on the root vector labelled by a $\ast$, namely $x_{\ell'}$,

 Thus $I(e_{\mathcal C'}+ x_{\ell'})\neq 0$, so a fortiori $I(\overline{B.\mathfrak u^{\mathcal C'}}+ x_{\ell'})\neq 0$.

 Suppose that  $\overline {B.\mathfrak u^\mathcal C} = \overline {B.\mathfrak u^{\mathcal C'}} $. Then $I(\overline{B.\mathfrak u^{\mathcal C}}+ x_{\ell'})\neq 0$.  This contradicts the conclusion of Corollary \ref {7.4}, so proving the required assertion.

 \end {proof}

%\textbf{Example $18$.}  Consider the composition $(3,2,1,3,1,2)$.  For $\mathcal C$ put $5$ below $9$ which describes the penetrating $i$-string for the rectangle $R^2_{C,C'}$, where $C=C_2,C'=C_6$, which stops at $C_4$.  (This does not include putting $7$ below $10$ which makes $x_{7,10}$ an excluded root vector.) For $\mathcal C'$ put $7$ below $10$, $8$ below the new appearance of $7$ and $3$ below $9$.  Then the pair $\mathcal C',\mathcal C$ makes $\ell_{8,10}$ the rightmost line.
%
%This means that $x_{7,10}$ lies in $\mathfrak u^{\mathcal C'}_{C,C'}$ but not in $\mathfrak u^{\mathcal C'}$. Thus we still have $B.x_{8.11} \in \mathfrak u^{\mathcal C'}_{C,C'}$.  Note that here $r=3,r'=4$.

 \textbf{Example $18$.}  Consider composition $(3,2,1,3,2,1)$.
 For $\mathcal C$ put $5$ below $9$ which describes the penetrating $i$ string for the rectangle $R^3_{C,C'}$, $C=C_1,C'=C_4$.  It is halted at $C_4$.  (Notice it does not include putting $10$ down below $12$, which make $x_{7,11}$ a secondary excluded root vector.) For $\mathcal C'$ put $8$ below $11$. The
  pair $\mathcal C',\mathcal C$ makes $\ell_{8,11}$ the rightmost line.  Now, by the observation in parentheses $x_{7,11}$ lies in $\mathfrak u^{\mathcal C'}_{C,C'}$ but not in $\mathfrak u^{\mathcal C'}$. Thus we still have $B.x_{8,11} \in \mathfrak u^{\mathcal C'}_{C,C'}$.
  Here $r=3,r'=4$.

\

 \textbf{Example $19$.}  Consider composition $(3,2,1,2,2,1,3)$.  For $\mathcal C$ put $7$ down below $10$ which describes the penetrating $i$ string for the rectangle $R^2_{C,C'}$, where $C=C_4,C'=C_5$.  It is halted at $C_5$.   Now for $\mathcal C'$ put $9$ below $11$ and $8$ below $10$.  The pair $\mathcal C',\mathcal C$ makes $\ell_{9,11}$ the rightmost line.

  Yet the $i$-string in $\mathscr T^{\mathcal C'}$ may be continued by putting $7$ into $C_7$ below $14$.  This makes $x_{9,14}$  a primary  excluded root vector. This means that $x_{9,14}$ lies in $\mathfrak u^{\mathcal C'}_{C,C'}$ but not in $\mathfrak u^{\mathcal C'}$. Thus we still have $B.x_{9.11} \in \mathfrak u^{\mathcal C'}_{C,C'}$.  Here $r=4,r'=5$.

 \

 \textbf{Example $20$.}  Consider the composition $(3,2,1,3,1,2)$.  For $\mathcal C$ put $5$ below $9$ which describes the penetrating $i$-string for the rectangle $R^2_{C,C'}$, where $C=C_2,C'=C_6$, which is halted at $C_4$.  (This does not include putting $7$ below $10$ which makes $x_{7,10}$ a primary excluded root vector.) For $\mathcal C'$ put $7$ below $10$, $8$ below the new appearance of $7$ and $3$ below $9$.  Then the pair $\mathcal C',\mathcal C$ makes $\ell_{8,10}$ the rightmost line.

 This means that $x_{7,10}$ lies in $\mathfrak u^{\mathcal C'}_{C,C'}$, but not in $\mathfrak u^{\mathcal C'}$. Thus we still have $B.x_{8,11} \subset \mathfrak u^{\mathcal C'}_{C,C'}$. Here $r=3,r'=4$.

\section {Towards Surjectivity} \label {8}

\subsection {Surjectivity for Partitions} \label {8.1}

Our most general result for subjectivity of the component map is when $(c_1,c_2,\ldots,c_k)$ is a partition of $n$.  Ironically it is rather boring and indeed paradoxical. Indeed to paraphrase Hilbert, every schoolboy in Cambridge (England!) knows the number of compositions of $n$ to be $2^{n-1}$, whilst only the combined genius of Hardy and Ramanujan (working together at Trinity College, Cambridge during the Great War) could give a formula for the number of partitions of $n$ and even that was only asymptotic.  Of course the number of partitions of $n$ is just the number of symmetric groups $S_n$ orbits in the set of compositions of $n$.  On the other hand $S_n$ also permutes the weights of the generators of the ``semi-invariant'' algebra $I=\mathbb C[\mathfrak m]^{P'}$, \cite [Prop. 3.4.8, Lemma 4.5]{FJ1}. Yet we hope we have convinced the reader that describing the components of the zero locus $\mathscr N$ of $I_+$ is a hard problem, whilst by contrast we have

\begin {lemma}  For a partition of $n$, the nilcone $\mathscr N$ is irreducible.
\end {lemma}

\begin {proof} In the case of a partition, neighbouring columns are adjacent.  Thus the Benlolo-Sanderson invariant generators are polynomials in the co-ordinates described by lines joining a pair of adjacent columns and are therefore pairwise distinct for the different pairs.  Since every generator is itself an irreducible polynomial \cite [Cor, 5.3]{FJ1}, it follows that the quotient algebra $Q:=\mathbb C[\mathfrak m]/I_+$ is a tensor product of (finitely generated) commutative domains. This does not quite say that $Q$ is a domain (it could admit nilpotent elements) but it does say that $\mathscr N$ is a finite Cartesian product of irreducible closed affine varieties and hence is irreducible \cite [Thm. 3, p.35]{Sh}.
\end {proof}

\subsection {Surjectivity in the Presence of Linear Generators.} \label {8.2}

\

Our second less general but more telling result is when $I$ admits ``enough'' linear generators. This allows ``factorisation'' by which the components of $\mathscr N$ are clearly seen.

Examples comes from the composition $(2,1^k,2)$ for $k\geq 2$.  All cases are similar but the case $k=3$ is already interesting.  In general there are $k-1$ linear generators and a generator of degree $2+k$ which factorises modulo the linear generators as a product of two $2 \times 2$ minors (having degree $2$) and $k-2$ further linear factors, whose zero loci give altogether $k$ irreducible components.  On the other hand $\mathscr B^2_1$ has cardinality $k$ leading to $k$ component tableaux.  Then injectivity implies surjectivity.  Even without banal counting one easily links up a component given by factorisation with a component tableau.

The composition $(2,1,1,1,2)$ exhibits a further surprising feature.  In this case $\mathscr B^2_1=\{2,3,4\}$, giving components labelled by $2,3,4$.  (See Example $15$). The component labelled by $3$, say $\mathscr T^{3}$, has an extra linear generator namely $x_{3,5}$.  It is the commutator $[x_{3,4},x_{4,5}]$.  Thus $\mathfrak u^{3}$ is complemented by a (Lie) subalgebra and so the component it defines is an orbital variety closure.

 Yet this is false for the other two components!  These are interchanged by a Dynkin diagram involution, so it suffices to consider just one.

The image of the (canonical) component tableau $\mathscr T^{4}$ under the component map is $\mathscr V:=\overline{B.\mathfrak u^4}$. By \cite [Lemma 6.9.9]{FJ2} it is enough to show that $\frac{1}{2}\dim G.\mathscr V > \dim \mathscr V$.

 Now the right hand term has dimension $\dim \mathfrak m -\textbf{g}=\dim \mathfrak n-5$, by Theorems \ref {6.3}, \ref {6.1}.

On the other hand $G.e_4 \subset G.\mathscr V$. Here the lines with label by $1$ defining $e_4$  form monomial chains $(1,3,5,6),(2,4),7$.  It follows that the
nilpotency class of $e_4$ is $(4,2,1)$ whose conjugate partition is $(3,2,1,1)$.  Thus $\dim G.\mathscr V \geq\dim G.e = 2\dim \mathfrak n - 3.2-2.1 = 2(\dim \mathfrak n-4)$, which by the previous paragraph results in  the required strict inequality.

By contrast the nilpotency class of $e_3$ is $(3,2,2)$ giving $\dim G.e_3 = \dim \mathscr N - 3.2-3.2 = 2(\dim \mathfrak n-6)$.  All this means is that $P.e_3$ is not dense in $\overline{B.\mathfrak u^3}$.  This does not mean that we have found a component with no dense $P$ orbit, since $P.(e_3+x_{3,6})$ is a dense $P$ orbit in $\overline{B.\mathfrak u^3}$. Matters of this nature were studied in \cite [6.10]{FJ2} for the canonical component.  We have not extended this theory to all components.

\subsection {Further Factorisation via Linear Generators.} \label {8.3}

\

\textbf{Example $21$.} Consider the composition $(3,2,1,1,1,2,3)$. In this case there are $6$ component tableaux.  First we lower $6$ (resp. $7$) below $7$ (resp. $8$).  Then $\mathscr B^2_{2,6}={i}:i=5,6,7$.  When $i$ is lowered into $R_3$ we obtain $\mathscr B^3_{3,i}$.  Then $3$ or $i$ can be lowered into $R_4$ giving $3\times 2=6$ tableaux.

As in the Example of \ref {8.3} with $k=3$, the Benlolo-Sanderson invariant of degree $5$, namely $I^2_{C_2,C_6}$ factorizes modulo the linear generators $x_{6,7},x_{7,8}$ as $x_{6,8}AB$ where $A$ (resp. $B$) is the $2\times 2$ minor defined by the rows $4,5$ and columns $6,7$ (resp. rows $7,8$ and columns $9,10$). Thus the resulting zero set has three components.

  The component defined by $x_{6,8}=0$ corresponds to taking $i=6$ in the first paragraph above.   Modulo $x_{6,7},x_{7,8},x_{6.8}$, the Benlolo-Sanderson invariant of degree $10$, namely $I^3_{C_1,C_7}$ factorizes via \cite [Lemma 1.10]{FJ2} as a product of two irreducible degree five polynomials.  Their zero sets respectively correspond to taking $3$ (resp. $6$) from the batch $\mathscr B^3_{C_1,C_7}=\{3,6\}$.

  The component defined by $A=0$ (resp. $B=0$) corresponds to taking $i=5$ (resp. $i=7$) in the first paragraph above. Adjoining $I^3_{C_1,C_7}$, the resulting set has components corresponding to taking either elements of the batch $\mathscr B^7_{C_1,C_7}$, which is $\{3,5\}$ (resp. $\{3,7\}$).  However unlike the previous case $I^3_{C_1,C_7}=0$ does not factorize neither modulo $A$, nor modulo $B$.  Thus we cannot be sure to have so obtained \textit{all} the components of $\mathscr N$.

  This situation is rather analogous to the two component tableaux for each of the compositions $(2,1,1,2)$ and $(3,2,2,3$).  The first case corresponds to taking $k=2$ in \ref {8.2} in which $I^2_{C_1,C_4}$ factorizes modulo the linear invariant $x_{3,4}$  as two irreducible minors of degree $2$ thereby proving there are just two components of $\mathscr N$.  In the second case $I^3_{C_1,C_4}$ does \textit{not} factorize modulo $I^2_{C_2,C_3}$ (which is a $2 \times 2$  minor), so we cannot be sure that there are again just two  components of $\mathscr N$.  Yet it is plausible - just as it is plausible that in our previous case we have found all components.

 \section {Index of Notation and of Notions} \label {9}

 Notation and notions frequently used are given below in the section where they are first defined and sometimes also first mentioned.

 \

 \textbf{Notations}

 \

\ref {1}.  \ \ \ \ $\mathbb C,[m,n],]m,n],[m,n[,]m,n[$.

\ref {1.1}. \ \ $\mathscr N$.

\ref {1.3}. \ \ $G, \mathfrak p$.

\ref {1.4}. \ \ $B, \mathfrak n,\mathfrak b,P,\mathfrak m,\mathfrak p, H, \mathfrak h, W, P'$.

\ref {2.1.1}. $\mathscr D,C_i,c_i,R_i,R^i,b_{i,j}$.

\ref {2.1.2}. $C_i^s$.

\ref {2.1.3}. $[C,C'], [C,C'[, ]C,C'],]C,C'[$.

\ref {2.1.4}. $R^s_{C,C'}$.

\ref {2.2.1}. $\textbf{M},\mathfrak r, \textbf{B}_i,\textbf{C}_i,x_{i,j}$.

\ref {2.2.2}. $\mathscr T,C_{(i)},R_{(i)}$.

\ref {2.3.1}. $\textbf{M}^*_s$.

\

\ref {2.3.2}. $d_\mathscr D^{R_{C,C'}^s},d_\mathscr D^{R^s_{C,C'}}$.

\ref {2.3.3}. $\mathfrak m^-, x^*_{i,j}$.

\ref {2.3.5}. $\Id_s$.

\ref {3.1.1}. $\mathscr T(t),\mathscr T(\infty),\mathscr T(t,\ell),C(t),C(\infty),C(t,\ell)$.

\ref {3.1.4}. $\mathscr T^\mathcal C(\infty)\mathscr B^t_{C,C'},\mathscr B^s$.

\ref {3.2.4}. $\ell_{i',j'}, \alpha_{i',j'}$.

\ref {3.2.6}. $C^{\leq p},C^{>p}$.

\ref {3.2.7}.   $\mathscr T^\mathcal C, \mathscr T^{i,j,\ldots}$.

\ref {4.1.3}. $C_g,C'_g, \mathscr T_{i,j}, \mathscr T_{i,\textbf{j}}$.

\ref {4.1.4}. $w_\mathscr T,w_{i,j},\mathfrak u_{i,j},\mathfrak l^+$.

\ref {4.1.5}. $\mathfrak u^\mathcal C, X, \textbf{X}$.

\ref {4.2.1}. $S(s)$.

\ref {4.3.1}. $\mathscr C, \textbf{g}$.

\ref {4.3.5}.  $i_{C,C'}^\mathcal C$.

\ref {4.9}. \ \ \ $\textbf{M}_s$.

\ref {5.1}. \ \ \ $e_\mathcal C,V_\mathcal C$.

\ref {5.1.1}. $\ell_{C,C'}^\ast$.

\ref {5.2.2}. $\mathscr L_i$.

\ref {6.2.1}.  $Y,S,Z, Y_{i,\textbf{j}}, X_{i,\textbf{j}}$.

\ref {6.2.2}.  $Y_{i,\textbf{j}}, X_{i,\textbf{j}}$.

\ref {6.3}. \ \ \ $\textbf{S},\textbf{X},\textbf{Y},\textbf{Z}$.

\

\textbf{Notions.}

\

\ref {1.2}.  \ \ Weierstrass Section.

\ref {1.5}.  \ \ Composition Tableau.

\ref {1.7}.  \ \ Component Tableau.

\ref {1.8}.  \ \ Orbital Variety.

\ref {1.9}.  \ \ Hypersurface Orbital Variety.

%\ref {2.1.1}.  Diagram.

\ref {2.1.2}.  Adjacent and Neighbouring Columns.

\ref {2.1.4}.  Rectangles.

\ref {2.2.2}.  Tableau.

\ref {2.3}. \ \ \ Benlolo-Sanderson invariants.

\ref {2.3.4}.  Quantization, Unitarizable Highest Weight Modules.

\ref {3}. \ \ \ \ \  Component Tableau.

\ref {3.1.2}. The Vav Conversive.

\ref {3.1.3}. Surrounding Columns.

\ref {3.1.4}.  Batches.

\ref {3.2.5}.  Lines with label $\ast$.

\ref {3.2.6}.  Lines with label $1$, lines with a neutral label.

\ref {4}. \ \ \ \ \ Excluded Roots.

\ref {4.1.3}. Primary and Secondary shifting of roots.

\ref {4.1.4}. Primary and Secondary excluded roots.

\ref {4.1.5}. Encircling.

\ref {4.2.1}. An $s$-string.

\ref {4.2.2}. Starting Places.

\ref {4.2.3}. Primary Excluded Roots.

\ref {4.2.4}. Secondary Excluded Roots.

\ref {4.3.5}. The $i$ - string of penetration.

\ref {5}. \ \ \ \ \ Weierstrass Sections.

\ref {5.1.1}. Composite lines.

\ref {5.2}. \ \  Collapsing.

\ref {5.2.1}. Special Composite Lines.

\ref {6.4}. \ \  The Component Map.

\ref {7.1}. \ \  The Exchange Lemma.

\ref {7.3}. \ \ The Rightmost Line.

%\subsection*{Tableaux}
\begin{figure}[H]
\begin{center}
\begin{tikzcd}[row sep=0.5em,%tiny,
column sep = 1.8em]
  &C_r&&&C_{r+1}\\
R_1&\boxed{i_1} \arrow[-,line width = 1 pt]{rrr}{1}&&&\boxed{j_1}\\
R_2&\boxed{i_2}\arrow[-,line width = 1 pt]{rrr}{1}&&&\boxed{j_2}\\
 \vdots&\vdots&&\vdots\\
R_{m_1}&{\boxed{i_{m_1}}}\arrow[-,line width = 1 pt]{rrr}{1}&&&\boxed{j_{m_1}}\\
R_{m_1+1}&\red{\boxed{i_{m_1+1}}}\arrow[-,dashed, line width = 1 pt,  green,rrrddddd]&&&\boxed{j_{m_1+1}}\\
R_{m_1+2}&{\boxed{i_{m_1+2}}}\arrow[-, black , line width = 1 pt]{rrru}{1}&&&\boxed{j_{m_1+2}}\\
R_{m_1+3}&{\boxed{i_{m_1+3}}}\arrow[-, black , line width = 1 pt]{rrru}{1}&&&\boxed{j_{m_1+3}}\\
 \vdots&\vdots\arrow[-, black , line width = 1 pt]{rrru}{1}&&&\vdots\\
  \vdots& &&&\boxed{j_{m_1+m_2}}\\
R_{m_1+m_2+1}&\red{\boxed{i_{m_1+m_2+1}}}\arrow[-,dash dot, line width = 1 pt,  blue ]{rrru}{1}\arrow[-,dashed, green, line width = 1 pt,rrrd]&&&\red{\boxed{i_{m_1+1}}} \arrow[-,swap,  line width = 1 pt,red]{u}{*} \arrow[-,swap, bend right=40,swap, line width = 1 pt,red]{uuu}{*}\arrow[-,swap, bend right=45,swap, line width = 1 pt,red]{uuuu}{*}\arrow[-,swap, bend right=50,swap, line width = 1 pt,red]{uuuuu}{*}\\
R_{m_1+m_2+2}&\red{\boxed{i_{m_1+m_2+2}}}\arrow[-,dashed, line width = 1 pt, green, rrrd]&&&\red{\boxed{i_{m_1+m_2+1}}}\arrow[-,swap,bend right=40,  line width = 1 pt,red]{uu}{*}\\
R_{m_1+m_2+3}&\red{\boxed{i_{m_1+m_2+3}}}\arrow[-,dashed, line width = 1 pt, green, rrrd]&&&\red{\boxed{i_{m_1+m_2+2}}}\arrow[-,swap,bend right=45,  line width = 1 pt,red]{uuu}{*}\\
 \vdots&\vdots&&&\vdots\\
R_{m_1+m_2+m_3-1}&\red{\boxed{i_{m_1+m_2+m_3-1}}}\arrow[-,dashed, line width = 1 pt,  green, rrrd]&&&\vdots&\\
R_{m_1+m_2+m_3}&\red{\boxed{i_{m_1+m_2+m_3}}}\arrow[-,dash dot, swap, bend right=8,  line width = 1 pt, blue]{uuuuuurrr}{1}&&&\red{\boxed{i_{m_1+m_2+m_3-1}}}\arrow[-,swap,bend right=50,  line width = 1 pt,red]{uuuuuu}{*}\\
R_{m_1+m_2+m_3+1}&\red{\boxed{i_{m_1+m_2+m_3+1}}}\arrow[-,dashed, swap,  line width = 1 pt, green]{rrr}&&&\red{\boxed{i_{m_1+m_2+m_3+1}}}\\
\end{tikzcd}
\caption{ This illustrates Overview,  \ref {3.2.6}.  The dashed green lines are neutral lines from $C_r(\infty)$ to $C_{r+1}(\infty)$. The solid black lines with label $1$ join an entry in $C_r(\infty)$ to an entry in $C_{r+1}$.  This is also true of the dashed/dotted blue lines, but the existence of the latter is conditional on their being an entry in the left hand factor and from which there is no outgoing neutral line.  (In particular the uppermost dashed/dotted line would be absent if there were no entry $i_{m_1+m+2+1}$ in the left hand column or if it were joined to an entry in the right hand column by a neutral line, as in top diagram rather the bottom diagram in Fig. $3$.) The vertical red lines join an entry of $C_{r+1}(\infty)\setminus C_{r+1}$ to $C_{r+1}$.
 } \label{fig1}
\end{center}
\end{figure}
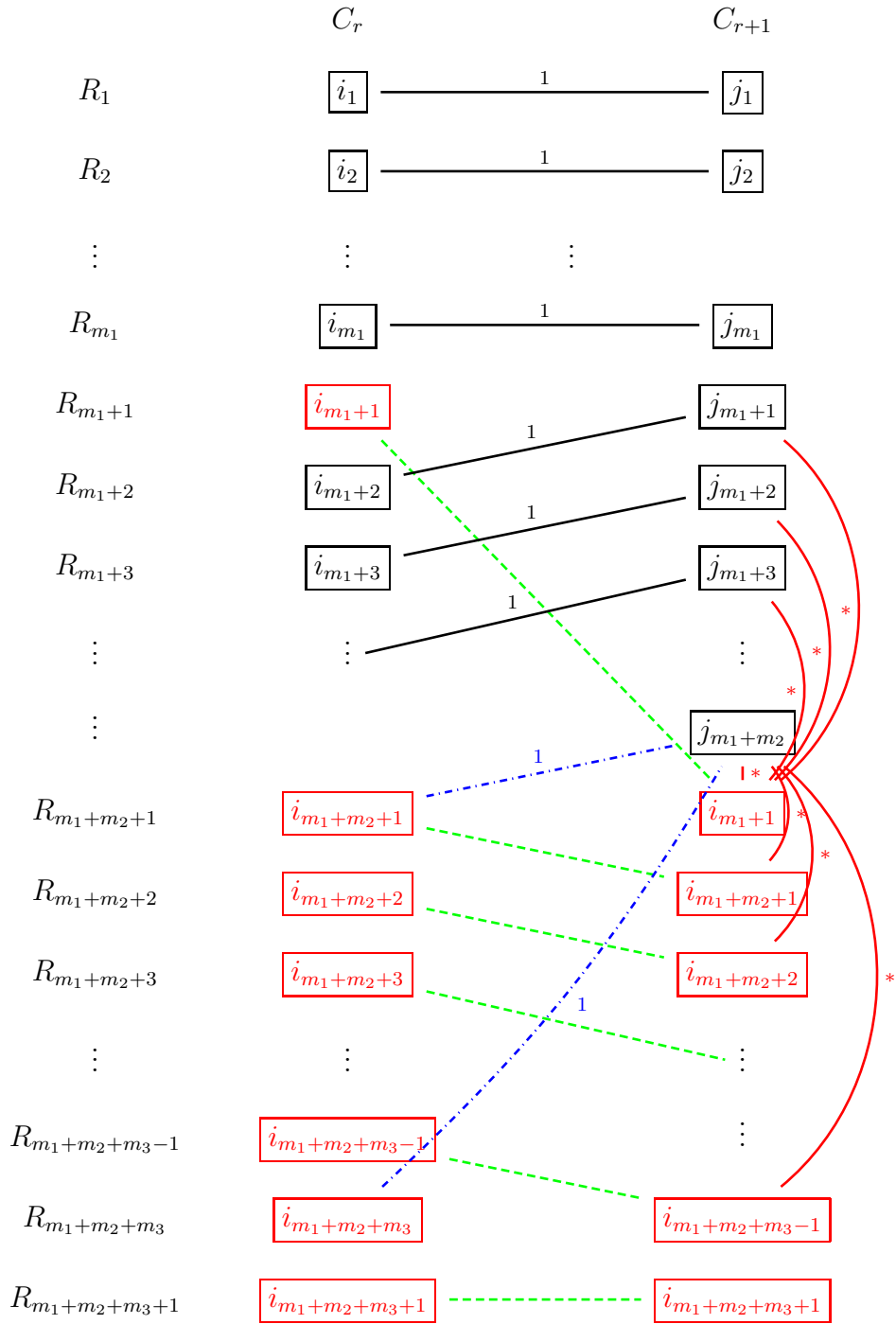
%%%%%%%%%%%%%%%%%%%%%%%%%%%%%%%%%%%%%%%%%%%%%%%%%%

\begin {center} \textbf{Examples.}

\

\end {center}
Some of the examples of the text are illustrated, with the numbering retained.
In this a composition defines the standard parabolic being considered. Less and less details are given as it is presumed that the reader will become joyfully adept at using the rules set out in the main text, particularly \ref {3.2.2}, \ref {3.2.6}, \ref {4.1.4}, \ref {4.3.6}, \ref {6.2}, \ref {7.3}.

\

\textbf{Example 3}. %Let $P$ a parabolic define by the composition $(1,2,1,2)$.

\begin{figure}[H]
\begin{center}
\begin{tikzcd}[row sep=0.2 em,%tiny,
column sep = 0.3em]
& 1&2&4&5\\
& & 3&&6 \\
\end{tikzcd} \, 	$\rightarrow$  \,
\begin{tikzcd}[row sep=0.2 em,%tiny,
column sep = 0.2em]
 1&  2&{4}&5\\
&3& \red2 & 6 \\
\end{tikzcd} \,  $ \begin{array}{cc}\nearrow&\begin{tikzcd}[row sep=0.2 em,%tiny,
column sep = 0.3em]
 1&    2 &{4}&5\\
&   3&    {   {\red 2}}& 6 \\
&&&\red2\\
\end{tikzcd}  \\

\searrow&\begin{tikzcd}[row sep=0.2 em,%tiny,
column sep = 0.3em]
 1&  2&{4}&5\\
&   3&  \red 2& 6 \\
&&\red3&3\\
\end{tikzcd}

 \end{array}$
\\

\caption{ Consider the composition $(1,2,1,2)$.  There are two component tableaux.  Starting from $\mathscr T$, one first moves $2$ below $4$.  Then $\mathscr B^2_{2,4} = \{3,2\}$.  Thus either $2$ can be lowered below $6$ or $3$ below $2$ and translated horizontally under $6$.}
\label{fig2}
$ \begin{array}{cc}\nearrow&\begin{tikzcd}[row sep=0.2 em,%tiny,
column sep = 0.3em]
 1&    2 &{4}&5\\
&   3&    {   {\red 2}}& 6 \\
&&&\red2\\
\end{tikzcd}  \rightarrow \begin{tikzcd}[row sep=0.5 em,%tiny,
column sep = 0.8em]
 1\arrow[- ,  r,"1"]&    2 \arrow[-,dashed, line width = 1 pt, green, rd]&{4}\arrow[- ,  r,"1"]&5\\
&   3\arrow[- ,  ru,"1"]&   {   {\red 2}} \arrow[-,dashed, line width = 1 pt, green, rd]\arrow[- ,  u,red,"*",swap]& 6 \\
&&&\red2\arrow[- ,  u,red,"*",swap]\\
\end{tikzcd} \\
\searrow&\begin{tikzcd}[row sep=0.2 em,%tiny,
column sep = 0.3em]
 1&  2&{4}&5\\
&   3&  \red2& 6 \\
&&\red3&3\\
\end{tikzcd} \rightarrow \begin{tikzcd}[row sep=0.5 em,%tiny,
column sep = 0.8em]
 1\arrow[- ,  r,"1"]& 2&{4}\arrow[- ,  r,"1"]&5\\
&   3 \arrow[-,dashed, line width = 1 pt, green, rd]&  \red 2  \arrow[- ,  u,red,"*",swap] \arrow[- ,  r,"1"]& 6 \\
&&\red3 \arrow[- ,  uu ,red, bend  left=35 ,"*"] \arrow[-,dashed, line width = 1 pt, green, r]&3\\
\end{tikzcd}
 \end{array}$

\caption{  This describes the decoration of the lines in
%Figure $2$ coming from \ref {3.2.6} as depicted in Figure $1$}
$\mathscr T(\infty)$ coming from the two component tableaux of Figure $2$ following \ref {3.2.6} and Figure $1$. Notice that in the top diagram, there is no line (with label $1$) to $6$ because there is no entry in the box $R_3\cap C_3$ (which could be remedied by adding a column of height $3$ on the extreme  left) and because $2$ is joined by a neutral line to $C_3$.}

\label{fig3}
\end{center}
\end{figure}

\

\

\

\

%\begin{landscape}
\textbf{Example 4}. %Let $P$ a parabolic define by the composition $(2, 1, 1, 2, 1)$.

\begin{figure}[H]
\begin{center}
\begin{tikzcd}[row sep=0.2 em,%tiny,
column sep = 0.2em]
& 1&   { { {\black 3}}}&   { { {\black 4}}}&   { {\black 5}}&7\\
&   { {\black 2}}& & &6& \\
\end{tikzcd}\, $\rightarrow$
\bigskip\bigskip\bigskip
 \begin{tikzcd}[row sep=0.2 em,%tiny,
column sep = 0.2em]
&  1&{3}&   { { {\black 4}}}&   { {\black 5}}&7\\
&   { {\black 2}}&  { {\black 2}} &  { {\red3}} &6& \\
\end{tikzcd}
%$\rightarrow$ %\\
%\begin{tikzcd}[row sep=0.2 em,%tiny,
%column sep = 0.5em]
%       &C_1&C_2&C_3&C_4&C_5\\
%&  1&   { { {\black 3}}}&   { { {\black 4}}}&   { {\black 5}}&7\\
%&   { \encircled 2} &  {\encircled 2 }&  { \encircled 3 }&6&5 \\
%\end{tikzcd}
$ \begin{array}{cc}
\nearrow&\begin{tikzcd}[row sep=0.2 em,%tiny,
column sep = 0.2em]

& 1& 3& 4&  {5}&7\\
& 2& { 2} & \red 3 &6&  {\red 5}\\
& & &  \red 2 & &
\end{tikzcd}   \\
 \rightarrow &\begin{tikzcd}[row sep=0.2 em,%tiny,
column sep = 0.2em]

& 1& 3& 4&  {5}&7\\
& 2& 2&{\red 3} &6& {\red 5} \\
& & & &   \red 3&
\end{tikzcd}\\
\searrow  &\begin{tikzcd}[row sep=0.2 em,%tiny,
column sep = 0.2em]

&  1&  3 &   {{\red 4}}& 5 &7\\
&  2& 2 & 3 &6&6 \\
& & & &  \red 4 &
\end{tikzcd}

\end{array}$
\\

\caption{ This describes the three component tableaux for the composition $(2,1,1,2,1)$ following assiduously the rules of section \ref {3.2.2}.  Notably in the bottom tableau $4$ is lowered by two rows below $6$ using the fact that $4 \in \mathscr B^1_{3,5} \cap  \mathscr B^2_{1,4}$ and this excludes any further use of these pairs of neighbouring columns.   }
\label{fig4}

\

$ \begin{array}{cc}
\nearrow&\begin{tikzcd}[row sep=0.5 em,%tiny,
column sep =1em]
& 1\arrow[- ,  r,"1"]& 3  \arrow[-,dashed, line width = 1 pt, green, rd] & 4\arrow[- ,  r,"1"]&  { 5} \arrow[-,dashed, line width = 1 pt, green, rd]&7\\
& 2 \arrow[-,dashed, line width = 1 pt, green, r]& { 2}  \arrow[-,dashed, line width = 1 pt, green, rd] & \red 3\arrow[- ,  u,red ,"*"] &6\arrow[- ,  ur,"1"]&  {\red 5}  \arrow[- ,  u,red ,"*"]\\
& & &  \red 2 \arrow[- ,  uu, bend right=35 ,red ,"*"] &2 &2
\end{tikzcd}   \\
 \rightarrow &\begin{tikzcd}[row sep=0.5 em,%tiny,
column sep = 1em]

& 1\arrow[- ,  r,"1"]& 3 \arrow[-,dashed, line width = 1 pt, green, rd]& 4\arrow[- ,  r,"1"]&  { 5} \arrow[-,dashed, line width = 1 pt, green, rd]&7\\
& 2 \arrow[-,dashed, line width = 1 pt, green, r]& 2\arrow[- ,  ur,"1"]&{\red 3}\arrow[- ,  u,red ,"*"]  \arrow[-,dashed, line width = 1 pt, green, rd]&6\arrow[- ,  ur,"1"]& {\red 5}\arrow[- ,  u,red ,"*"] \\
& & & &   \red 3 \arrow[- ,  u,red ,"*"]&
\end{tikzcd}\\
\searrow  &\begin{tikzcd}[row sep=0.5 em,%tiny,
column sep = 1em]

&  1\arrow[- ,  r,"1"]&  3  \arrow[-,dashed, line width = 1 pt, green, rd]&   {{4}} \arrow[-,dashed, line width = 1 pt, green, rdd]& 5\arrow[- ,  r,"1"] &7\\
&  2 \arrow[-,dashed, line width = 1 pt, green, r]& 2 \arrow[- ,  ur,"1"]& \red 3\arrow[- ,  u,"*",red ] &6&6 \\
& & & &  \red 4 \arrow[- ,  u,red ,"*",swap]\arrow[- ,  uu, bend right=35 ,red ,"*"]  &
\end{tikzcd}

\end{array}$

\

\caption{ This describes the decoration of the lines in Figure $4$.}
\label{fig5}
\end{center}
\end{figure}

\

\textbf{Example 5. } %Consider the composition $(3,2,1,2,2,1,3)$.

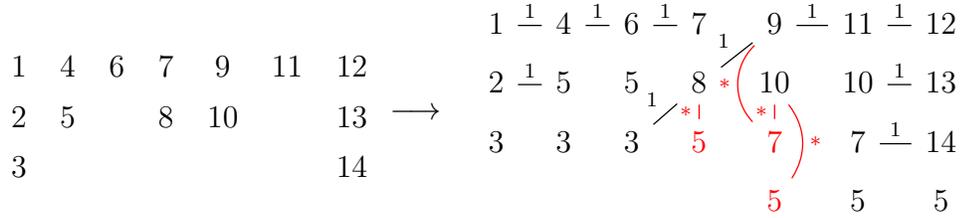
\begin{figure}[H]
\begin{center}{
$\begin{tikzcd}[row sep=0.2 em,%tiny,
column sep = 0.2em]
&  1&4& 6& 7&9&11&12\\
&  2&5&  & 8&10&&13\\
&  3& &  & & & &14\\
\end{tikzcd} \longrightarrow  \begin{tikzcd}[row sep=0.5em,%tiny,
column sep = 0.8em]
&  1\arrow[- ,  r ,"1"]&4\arrow[- ,  r ,"1"]& 6\arrow[- ,  r ,"1"]& 7&9\arrow[- ,  r ,"1"]&11\arrow[- ,  r ,"1"]&12\\
&  2\arrow[- ,  r ,"1"]&5&5  & 8\arrow[- ,  ur ,"1"]&   10 &  10 \arrow[- ,  r ,"1"]&13\\
&   3&    3 & 3 \arrow[- ,  ur ,"1"]& \red 5 \arrow[- ,  u , red,"*"]&\red 7 \arrow[- ,  u,red ,"*"]\arrow[- ,  uu, bend left=45,red ,"*"]  & 7\arrow[- ,  r ,"1"] &14\\
& && & & \red 5  \arrow[- ,  uu,red ,"*",bend right=35 ,swap] & 5&5 \\
\end{tikzcd}$

\caption{
This describes \textit{just one} component tableaux for the composition $(3,2,1,2,2,1,3)$.  In this $7$ is moved down two rows below $10$ and $5$ by one row below $8$ and then by a further one row below $7$ illustrating the paragraph just above Lemma \ref {3.2.5}.  The lines with label $1$ resulting from \ref {3.2.6} are also indicated.
 } \label{fig6}
  }\end{center}
\end{figure}

\textbf{Example 6}.% Consider the composition $(1, 2, 2, 1, 3, 2)$

\begin{figure}[H]
\begin{center}
\begin{tikzcd}[row sep=0.5em,%tiny,
column sep = 0.8em]
& 1&2&4&6&7&10\\
&  & 3&5& & 8  & 11\\
&  &  & & & 9  &
\end{tikzcd} $\longrightarrow$ \begin{tikzcd}[row sep=0.5em,%tiny,
column sep = 0.8em]
& 1\arrow[- ,  r ,"1"]&2\arrow[- ,  r ,"1"]&4&6\arrow[- ,  r ,"1"]&7\arrow[- ,  r ,"1"]&10\\
&  & 3&5\arrow[- ,  ur ,"1"]&\red4 \arrow[- ,  u , red,"*"]\arrow[- ,  r ,"1"]& 8  & 11\\
&  &  &\red 3 \arrow[- ,  u , red,"*"] &3\arrow[- ,  r ,"1"] & 9\arrow[- ,  ur ,"1"]  &  \red 8 \arrow[- ,  u , red,"*"]
\end{tikzcd}
\caption{
%The component tableau $\mathscr T^{4,8,3}(\infty)$ meaning that $4,8,3$ are sequentially chosen from the batches. The resulting lines in the component tableau obtained from Figure $1$ are drawn but omitting the neutral lines.
Choosing $4,8,3$ sequentially from the batches in the composition $(1,2,2,1,3,2)$ defines the component tableau $\mathscr T^{4,8,3}(\infty)$.
The resulting lines in the component tableau obtained from Figure $1$ are drawn but omitting the neutral lines.}
 \label{fig7}
  \end{center}
\end{figure}
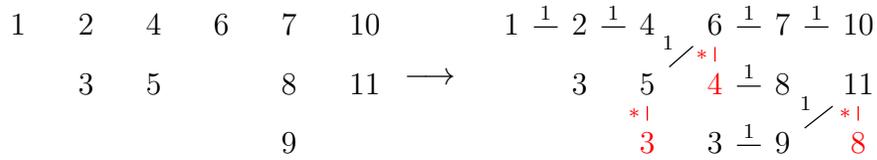
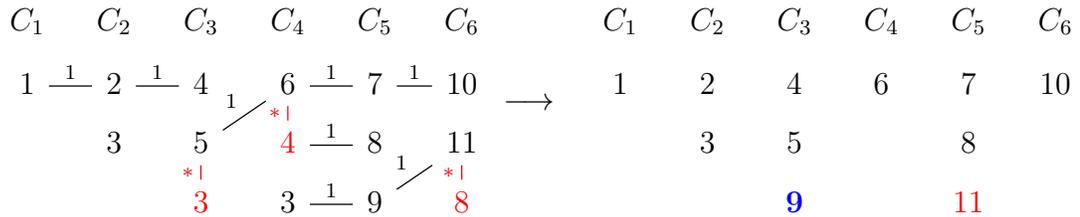
\begin{figure}[H]
\begin{center}
\begin{tikzcd}[row sep=0.5em,%tiny,
column sep = 0.8em]
 &C_1&C_2&C_3&C_4&C_5&C_6\\
& 1\arrow[- ,  r ,"1"]&2\arrow[- ,  r ,"1"]&4&6\arrow[- ,  r ,"1"]&7\arrow[- ,  r ,"1"]&10\\
&  & 3&5\arrow[- ,  ur ,"1"]&\red4 \arrow[- ,  u , red,"*"]\arrow[- ,  r ,"1"]&  8  & 11\\
&  &  &\red 3 \arrow[- ,  u , red,"*"] &3\arrow[- ,  r ,"1"] & 9\arrow[- ,  ur ,"1"]  &  \red 8 \arrow[- ,  u , red,"*"]
\end{tikzcd}   $\longrightarrow$
\begin{tikzcd}[row sep=0.5em,%tiny,
column sep = 0.8em]
 &C_1&C_2&C_3&C_4&C_5&C_6\\
& 1&2&4&6 &7&10\\
&  & 3&5 & & 8  & \\
&  &  &\blue\textbf 9 &&\red 11    &
\end{tikzcd}
\caption{
The right hand tableau is $\mathscr T_{8,11}$ of \ref {4.1.4} with respect to $\mathscr T^{4,8,3}(\infty)$.
Thus $11$ is put under $8$ and then $9$ is pushed under $6$ creating secondary excluded roots $\ell_{i',9}$: with $i'=4,5,6$ which one may observe are not labelled by a $1$, so satisfying Lemma \ref {4.2.4}. }  %Although $9$ could be further pushed under $3$ this would not follow our rule in \cite [4.1.3] {FJ5} and would result in the root $\ell_{3,9}$ labelled by a $1$ being encircled.

 \label{fig8}

\end{center}
\end{figure}

\textbf{Example 7 %19
} %Consider the composition $(2, 1, 2, 2, 1)$, as above

\

\begin{figure}[H]
\begin{center}
\begin{tikzcd}[row sep=0.5em,%tiny,
column sep = 1em]

& 1&3&4&6&8\\
&2& &5& 7& \\
\end{tikzcd}  $\longrightarrow  $
 \begin{tikzcd}[row sep=0.5em,%tiny,
column sep = 1em]
&  1\arrow[- ,  r ,"1"]&3&4\arrow[- ,  r ,"1"]&6\arrow[- ,  r ,"1"]&8\\
&   2& 2\arrow[- ,  ur ,"1"]&5& 7&7 \\
&   & &  \red {3}\arrow[- ,  u,red ,"*"]\arrow[- ,  uu, bend left=35 ,red ,"*"] \arrow[- ,  ur ,"1"]&\textbf {\violet 5}\arrow[- ,  u,red ,"*"]&  5& \\
\end{tikzcd}

\caption{This describes $\mathscr T^{3,5}$ for the composition $(2,1,2,2,1)$.
Given that  $3$ labelled in red, belongs to both batches  $\mathscr{B}^{1}$ and $\mathscr{B}^{2}$, whilst  $5$ in purple belong only  to $\mathscr{B}^{2}$, Therefore $3$ can  descend by two rows below $5$. Then $5$ descends one row below $7$.
 } \label{fig9}
 \end{center}
\end{figure}
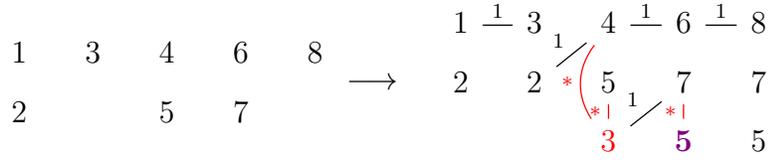

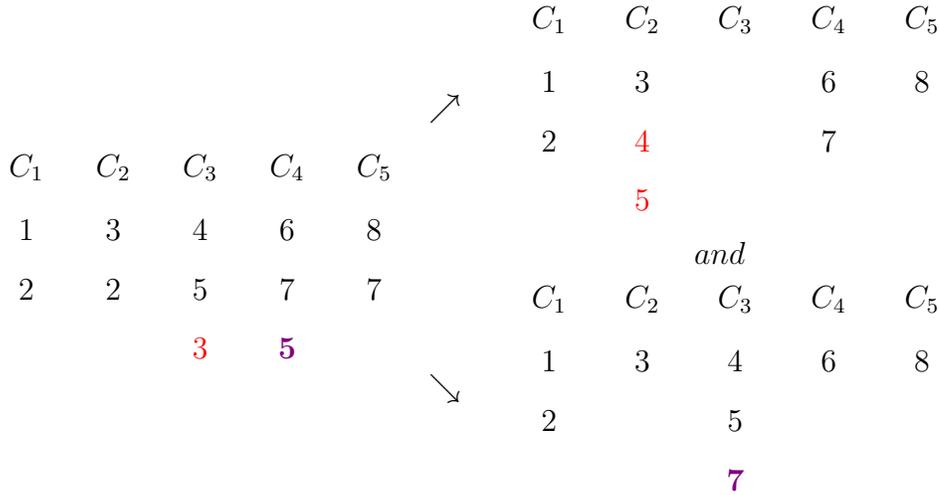
\begin{figure}[H]
\begin{center}
 \begin{tikzcd}[row sep=0.5 em,%tiny,
column sep = 0.8em]
 &C_1&C_2&C_3&C_4&C_5\\
&  1&3&4&6&8\\
&   2& 2&5& 7&7 \\
&   & &  \red {3}&\textbf {\violet 5} & \\
\end{tikzcd} $
 \begin{array}{cc}
 \nearrow &
\begin{tikzcd}[row sep=0.5em,%tiny,
column sep =1em]
 &C_1&C_2&C_3&C_4&C_5\\
&  1&3& &6&8\\
&   2& \red 4& & 7& \\
&   &\red5 & & & \\
\end{tikzcd} \\
 & and \\
\searrow &
\begin{tikzcd}[row sep=0.5em,%tiny,
column sep = 1em]
 &C_1&C_2&C_3&C_4&C_5\\
&  1&3&4&6&8\\
&   2& &  5& & \\
&   & &\textbf {\violet 7}& & \\
\end{tikzcd}
 \end{array}$
\caption{The top tableau on the right describes $\mathscr T_{3,\textbf{j}}$. Here $\textbf {j}:=(4,5)$ which is put below $3$. The bottom tableau describes $\mathscr T_{5,7}$.  In this $7$ is put below $5$.
These two tableaux cannot be amalgamated into a single tableau and so this component is not an orbital variety closure.}
%To draw the tableau    $\hat{\mathscr{T}}   $, you must position    $4   $ and    $5   $ beneath    $3   $, and    $7   $ beneath    $5   $. However, it is crucial not to place    $4   $,    $5   $, and    $7   $ beneath    $3   $ together, as this would exclude the element    $x_{3,7}   $ containing a    $1   $. Therefore,    $\hat{\mathscr{T}}   $ needs to be represented in two separate tableaux: one with    $4   $ and    $5   $ beneath    $3   $, and the other with    $7   $ beneath    $5   $, ensuring it doesn't interfere with our Weistrass section.
 \label{fig10}
\end{center}
\end{figure}

\

\

\

\

\

\

%%%%%%%%%%%%%%%%%%%%%%%%%%%%%%%%%%%%%%%%%%%%%%%%%%%%%%%%%%%%%%%%

\textbf{Example 8}. %Let the parabolic define by the sequence $ (2, 1, 2, 1, 2, %1) $
\begin{figure}[H]
\begin{center}
\begin{tikzcd}[row sep=0.5em,%tiny,
column sep = 1em]
       &C_1&C_2&C_3&C_4&C_5&C_6\\
R_1& 1&3&4&6&7&9\\
R_2& 2& &5& & 8  &
\end{tikzcd}
\caption{ This describes $\mathscr T$ for the composition $(2,1,2,1,2,1)$.  All five component tableaux are given below following the rules of \ref {3.2.2}, \ref {4.1.4} which the reader will now have no difficulty in applying.  Here $\mathscr T^{\mathcal C}(\infty)$ appears on the left,  $\mathscr T^{\mathcal C}$ appears on the right and below the corresponding matrix \textbf{M} is given with labels $1$, $\ast$ and the excluded root vectors encircled.
 } \label{fig2}
\end{center}
\end{figure}

\begin{center}
\begin{tikzcd}[row sep=0.5em,%tiny,
column sep = 1 em]
  {1}  \arrow[-,r,"1"]&3 &4 \arrow[-,r,"1"]& 6&7 \arrow[-,r,"1"]& {9} \\
  {2} & {2}   \arrow[-,ur,"1"]& {5}   &  {5}  \arrow[-,ur,"1"]& {8}& {8} \\
  &   &\red 3  \arrow[-,u,"*", red]  \arrow[-,uu,"*",bend right, red]  &  3   \arrow[-,ur,"1"]& \red 6\arrow[-,u,"*", red]  \arrow[-,uu,"*",bend right, red]& 6  \\
\end{tikzcd}  \;\;\;\;\;\;\;\; ,  \;\;\;\;\;\;\;\;
\begin{tikzcd}[row sep=0.5em,%tiny,
column sep = 1 em]
  {1}  \arrow[-,r,"1"]&    {3}  \arrow[-,r,"*", red] \arrow[-,dr,"*", red,bend right] \arrow[-,drrr,"1",bend right=75]&    {4}  \arrow[-,r,"1"] &   {6} \arrow[-,r,"*", red] \arrow[-,dr,"*", red,bend right]&  {7} \arrow[-,r,"1"]&  {9} \\
   {2} \arrow[-,rru,"1"]&  &  {5}   \arrow[-,rru,"1"] &  &  {8}&  \\
\end{tikzcd}  \\
\end{center}
%\;\;\;\;\;\;\;\;\;\;\;\;\;\;\;\;\;\;\;\;\;\;\;\;\;\;\;\;\;\;
\begin{center}
\footnotesize \begin{tikzpicture}
 \matrix [matrix of math nodes,left delimiter=(,right delimiter=)] (m)
 {
1 & 0  & {1} &   &   &   &   &   &      \\
 0& 1  &   & {1}   &   &   &   &   &      \\
  &   &  1 &  \red\cir { $\ast$} & \red\cir{ $\ast$} &   &   & 1  &      \\
  &   &   &   1& 0  &  1 &   &   &      \\
  &   &   &   0&  1 &   &  1 &   &      \\
  &   &   &   &   &  1 &   \red\cir{ $\ast$}&  \red\cir{ $\ast$} &      \\
  &   &   &   &   &   &  1 &  0 &  1    \\
  &   &   &   &   &   &   0&  1 &      \\
  &   &   &   &   &   &   &   & 1    \\
};
 %simple rectangle
 \draw (m-2-1.south west) rectangle (m-1-2.north east);
 \draw (m-3-3.south west) rectangle (m-3-3.north east);
 \draw (m-5-4.south west) rectangle (m-4-5.north east);
  \draw (m-6-6.south west) rectangle (m-6-6.north east);
   \draw (m-8-7.south west) rectangle (m-7-8.north east);
    \draw (m-9-9.south west) rectangle (m-9-9.north east);
 %\draw (m-6-5.south west) rectangle (m-5-6.north east);
 %\draw[rounded corners,ultra thick, draw=black, fill=blue, opacity=0.2] (m-2-3.south west) rectangle (m-1-4.north east);
%% \draw[rounded corners,ultra thick, draw=black, fill=red, opacity=0.2] (m-3-4.south west) rectangle (m-3-4.north east);
 %\draw[rounded corners,ultra thick, draw=black, fill=blue, opacity=0.2] (m-4-5.south west) rectangle (m-3-6.north east);
 %\draw[rounded corners,ultra thick, draw=black, fill=red, opacity=0.2] (m-4-3.south west) rectangle (m-4-3.north east);
 %\draw[rounded corners,ultra thick, draw=black, fill=blue, opacity=0.2] (m-4-5.south west) rectangle (m-3-6.north east);

 \end{tikzpicture}\\
\end{center}
 \bigskip

\begin{center}
\begin{tikzcd}[row sep=0.5em,%tiny,
column sep = 1 em]
  {1}  \arrow[-,r,"1"]&3 \arrow[-,r,"1"]&4 & 6&7 \arrow[-,r,"1"]& {9} \\
  {2} & {2}   & {5}  \arrow[-,ur,"1"] &   \arrow[-,u,"*", red] \red 4  \arrow[-,ru,"1"]& {8}& {8} \\
  &   & \red2  \arrow[-,u,"*", red]    &  {2} \arrow[-,ur,"1"]& \red 6\arrow[-,u,"*", red]  \arrow[-,uu,"*",bend right, red]& 6  \\
\end{tikzcd}  \;\;\;\;\;\;\;\; ,  \;\;\;\;\;\;\;\;
\begin{tikzcd}[row sep=0.5em,%tiny,
column sep = 1 em]
  {1}  \arrow[-,r,"1"]&    {3}  \arrow[-,r,"1"] &    {4}  \arrow[-,r,"*", red]  \arrow[-,rr,"1", bend left=35]&   {6} \arrow[-,r,"*", red] \arrow[-,dr,"*", red]&  {7} \arrow[-,r,"1"] \arrow[-,r,"1"]&  {9} \\
   {2} \arrow[-,rr,"*", red]&  &  {5}  \arrow[-,ru,"1"]  &  &  {8}&  \\
\end{tikzcd}  \\
\end{center}

\begin{center}
\begin{tikzpicture}
\footnotesize \matrix [matrix of math nodes,left delimiter=(,right delimiter=)] (m)
 {
1 & 0  & {1} &   & \cir{ }  &   &   &   &      \\
 0& 1  &   &   &\red\cir{ $\ast$}   &   &   &   &      \\
  &   &  1 & 1 & \cir{ } &   &   &   &      \\
  &   &   &   1& 0  &  \red\cir { $\ast$} & 1  &   &      \\
  &   &   &   0&  1 &   &   &   &      \\
  &   &   &   &   &  1 &   \red\cir{ $\ast$}&  \red\cir{ $\ast$} &      \\
  &   &   &   &   &   &  1 &  0 &  1    \\
  &   &   &   &   &   &   0&  1 &      \\
  &   &   &   &   &   &   &   & 1    \\
};
 %simple rectangle
 \draw (m-2-1.south west) rectangle (m-1-2.north east);
 \draw (m-3-3.south west) rectangle (m-3-3.north east);
 \draw (m-5-4.south west) rectangle (m-4-5.north east);
  \draw (m-6-6.south west) rectangle (m-6-6.north east);
   \draw (m-8-7.south west) rectangle (m-7-8.north east);
    \draw (m-9-9.south west) rectangle (m-9-9.north east);
 %\draw (m-6-5.south west) rectangle (m-5-6.north east);
 %\draw[rounded corners,ultra thick, draw=black, fill=blue, opacity=0.2] (m-2-3.south west) rectangle (m-1-4.north east);
%% \draw[rounded corners,ultra thick, draw=black, fill=red, opacity=0.2] (m-3-4.south west) rectangle (m-3-4.north east);
 %\draw[rounded corners,ultra thick, draw=black, fill=blue, opacity=0.2] (m-4-5.south west) rectangle (m-3-6.north east);
 %\draw[rounded corners,ultra thick, draw=black, fill=red, opacity=0.2] (m-4-3.south west) rectangle (m-4-3.north east);
 %\draw[rounded corners,ultra thick, draw=black, fill=blue, opacity=0.2] (m-4-5.south west) rectangle (m-3-6.north east);

 \end{tikzpicture}\\
\end{center}
  \bigskip

\begin{center}
\begin{tikzcd}[row sep=0.5em,%tiny,
column sep = 1 em]
{1}  \arrow[-,r,"1"]&3 &4 \arrow[-,r,"1"]& 6\arrow[-,r,"1"] &7 & {9} \\
{2} & {2}   \arrow[-,ur,"1"]& {5} & {5} &  {8}\arrow[-,ur,"1"]   &\red7\arrow[-,u,"*", red] \\
 &   & \red 3 \arrow[-,uu,"*",bend left, red] \arrow[-,u,"*", red]   & 3\arrow[-,ur,"1"]& \red {5} \arrow[-,u,"*", red]& 5 \\
\end{tikzcd}  \;\;\;\;\;\;\;\; ,  \;\;\;\;\;\;\;\;
\begin{tikzcd}[row sep=0.5em,%tiny,
column sep = 1 em]
  {1}  \arrow[-,r,"1"]&    {3} \arrow[-,drrr,"1",bend right=75] \arrow[-,r,"*", red]\arrow[-,dr,"*", bend right, red] &    {4}  \arrow[-,r,"1"] &   {6} \arrow[-,r,"1"] &  {7} \arrow[-,r,"*", red]&  {9} \\
   {2} \arrow[-,rru,"1"]&  &  {5}  \arrow[-,rr,"*", red]  &  &  {8}\arrow[-,ru,"1"]&  \\
\end{tikzcd}
\end{center}

\begin{center}
\footnotesize \begin{tikzpicture}
 \matrix [matrix of math nodes,left delimiter=(,right delimiter=)] (m)
 {
1 & 0  & {1} &   &   &   &   &   &      \\
 0& 1  &   &  1 &   &   &   &   &      \\
  &   &  1 & \red\cir { $\ast$} & \red\cir { $\ast$} &   &   &   & 1     \\
  &   &   &   1& 0  &  1&   &  \cir{ } &      \\
  &   &   &   0&  1 &   &   & \red\cir { $\ast$}  &  \red\cir { $\ast$}    \\
  &   &   &   &   &  1 &  1&\cir{ } &      \\
  &   &   &   &   &   &  1 &  0 &  \red\cir { $\ast$}    \\
  &   &   &   &   &   &   0&  1 &      \\
  &   &   &   &   &   &   &   & 1    \\
};
 %simple rectangle
 \draw (m-2-1.south west) rectangle (m-1-2.north east);
 \draw (m-3-3.south west) rectangle (m-3-3.north east);
 \draw (m-5-4.south west) rectangle (m-4-5.north east);
  \draw (m-6-6.south west) rectangle (m-6-6.north east);
   \draw (m-8-7.south west) rectangle (m-7-8.north east);
    \draw (m-9-9.south west) rectangle (m-9-9.north east);
 %\draw (m-6-5.south west) rectangle (m-5-6.north east);
 %\draw[rounded corners,ultra thick, draw=black, fill=blue, opacity=0.2] (m-2-3.south west) rectangle (m-1-4.north east);
%% \draw[rounded corners,ultra thick, draw=black, fill=red, opacity=0.2] (m-3-4.south west) rectangle (m-3-4.north east);
 %\draw[rounded corners,ultra thick, draw=black, fill=blue, opacity=0.2] (m-4-5.south west) rectangle (m-3-6.north east);
 %\draw[rounded corners,ultra thick, draw=black, fill=red, opacity=0.2] (m-4-3.south west) rectangle (m-4-3.north east);
 %\draw[rounded corners,ultra thick, draw=black, fill=blue, opacity=0.2] (m-4-5.south west) rectangle (m-3-6.north east);

 \end{tikzpicture}
\end{center}
  \bigskip

\begin{center}
 \begin{tikzcd}[row sep=0.5em,%tiny,
column sep = 1em]
  {1}  \arrow[-,r,"1"]&  {3} \arrow[-,r,"1"]& {4}  &   {6}\arrow[-,r,"1"]& {7} & {9} \\
  {2} &2  & 5\arrow[-,ru,"1"]&  \red 4   \arrow[-,u,"*", red]& {8}\arrow[-,ru,"1"]&  \red 7   \arrow[-,u,"*", red] \\
  &   & \red2 \arrow[-,u,"*", red]   & 2 \arrow[-,ur,"1"] & \red 4\arrow[-,u,"*", red] &  4  \\
\end{tikzcd} \;\;\;\;\;\;\;\; ,  \;\;\;\;\;\;\;\;
\begin{tikzcd}[row sep=0.5em,%tiny,
column sep = 1 em]
  {1}  \arrow[-,r,"1"]&    {3}  \arrow[-,r,"1"] &    {4}\arrow[-,drr,"*", bend right, red]  \arrow[-,r,"*"] &   {6} \arrow[-,r,"1"]&  {7} \arrow[-,r,"*", red]&  {9} \\
   {2} \arrow[-,rr,"*", red] \arrow[-,rrrr, "1", bend right]&  &  {5}  \arrow[-,ru,"1"]  &  &  {8}\arrow[-,ru,"1"]&  \\
\end{tikzcd} \end{center}

\begin{center}
\footnotesize \begin{tikzpicture}
 \matrix [matrix of math nodes,left delimiter=(,right delimiter=)] (m)
 {
1 & 0  & {1} &   &  \cir{ }  &   &   &   &      \\
 0& 1  &   &   &  \red\cir{$\ast$}&   &   & 1  &      \\
  &   &  1 &1 & \cir{ }  &   &   &   &     \\
  &   &   &   1& 0  & \red\cir{{ $\ast$} } &   &  \red\cir{ { $\ast$}} &      \\
  &   &   &   0&  1 & 1  &   &  &     \\
  &   &   &   &   &  1 &  1& \cir{ } &      \\
  &   &   &   &   &   &  1 &  0 &  \red\cir { $\ast$}    \\
  &   &   &   &   &   &   0&  1 &    1  \\
  &   &   &   &   &   &   &   & 1    \\
};
 %simple rectangle
 \draw (m-2-1.south west) rectangle (m-1-2.north east);
 \draw (m-3-3.south west) rectangle (m-3-3.north east);
 \draw (m-5-4.south west) rectangle (m-4-5.north east);
  \draw (m-6-6.south west) rectangle (m-6-6.north east);
   \draw (m-8-7.south west) rectangle (m-7-8.north east);
    \draw (m-9-9.south west) rectangle (m-9-9.north east);
 %\draw (m-6-5.south west) rectangle (m-5-6.north east);
 %\draw[rounded corners,ultra thick, draw=black, fill=blue, opacity=0.2] (m-2-3.south west) rectangle (m-1-4.north east);
%% \draw[rounded corners,ultra thick, draw=black, fill=red, opacity=0.2] (m-3-4.south west) rectangle (m-3-4.north east);
 %\draw[rounded corners,ultra thick, draw=black, fill=blue, opacity=0.2] (m-4-5.south west) rectangle (m-3-6.north east);
 %\draw[rounded corners,ultra thick, draw=black, fill=red, opacity=0.2] (m-4-3.south west) rectangle (m-4-3.north east);
 %\draw[rounded corners,ultra thick, draw=black, fill=blue, opacity=0.2] (m-4-5.south west) rectangle (m-3-6.north east);

 \end{tikzpicture}\\
\end{center}
  \bigskip

\begin{center}
\begin{tikzcd}[row sep=0.9 em,%tiny,
column sep = 1 em]
  {1}  \arrow[-,r,"1"]&  {3} \arrow[-,r,"1"] & {4}  &   {6}\arrow[-,r,"1"]  \arrow[-,d,"*", red]& {7} & {9} \\
  {2} &2  &5&  \red4 & {8}\arrow[-,ru,"1"]&   \red7   \arrow[-,u,"*", red] \\
  &   & \red2 \arrow[-,u,"*", red]  \arrow[-,uur,"1"]  & \red5 \arrow[-,uu,"*", red,bend right=40, red]& 5 &  5 \\
\end{tikzcd}  \;\;\;\;\;\;\;\; ,  \;\;\;\;\;\;\;\;
\begin{tikzcd}[row sep=0.5em,%tiny,
column sep = 1 em]
  {1}  \arrow[-,r,"1"]&    {3}  \arrow[-,r,"1"] &    {4}  \arrow[-,r,"*", red] \arrow[-,rrd,bend right=40,"1"]  &   {6} \arrow[-,r,"1"] &  {7} \arrow[-,r,"*", red]&  {9} \\
   {2} \arrow[-,rr,"*", red]&  &  {5}  \arrow[-,ru,bend right,"*", red]  &  &  {8}\arrow[-,ru,"1"]&  \\
\end{tikzcd}
 \end{center}
\begin{figure}[H]
\begin{center}
 \footnotesize \begin{tikzpicture}
 \matrix [matrix of math nodes,left delimiter=(,right delimiter=)] (m)
 {
1 & 0  & {1} &   & \cir{ }   &    &   &   &      \\
 0& 1  &   &   & \red\cir{ { $\ast$}}  & 1  &   &   &      \\
  &   &  1 &1 & \cir{ } &   &   &   &      \\
  &   &   &   1& 0  &  \red\cir{ { $\ast$}}&   & 1  &      \\
  &   &   &   0&  1 &  \red\cir{ { $\ast$}} &  &  &      \\
  &   &   &   &   &  1 &  1& \cir{ }&      \\
  &   &   &   &   &   &  1 &  0 &  \red\cir{{ $\ast$} }   \\
  &   &   &   &   &   &   0&  1 &    1  \\
  &   &   &   &   &   &   &   & 1    \\
};
 %simple rectangle
 \draw (m-2-1.south west) rectangle (m-1-2.north east);
 \draw (m-3-3.south west) rectangle (m-3-3.north east);
 \draw (m-5-4.south west) rectangle (m-4-5.north east);
  \draw (m-6-6.south west) rectangle (m-6-6.north east);
   \draw (m-8-7.south west) rectangle (m-7-8.north east);
    \draw (m-9-9.south west) rectangle (m-9-9.north east);
 %\draw (m-6-5.south west) rectangle (m-5-6.north east);
 %\draw[rounded corners,ultra thick, draw=black, fill=blue, opacity=0.2] (m-2-3.south west) rectangle (m-1-4.north east);
%% \draw[rounded corners,ultra thick, draw=black, fill=red, opacity=0.2] (m-3-4.south west) rectangle (m-3-4.north east);
 %\draw[rounded corners,ultra thick, draw=black, fill=blue, opacity=0.2] (m-4-5.south west) rectangle (m-3-6.north east);
 %\draw[rounded corners,ultra thick, draw=black, fill=red, opacity=0.2] (m-4-3.south west) rectangle (m-4-3.north east);
 %\draw[rounded corners,ultra thick, draw=black, fill=blue, opacity=0.2] (m-4-5.south west) rectangle (m-3-6.north east);

 \end{tikzpicture}\\

\caption{ From the above one may verify that the Exchange Lemma \ref {7.2} and the criterion for partial linearity Lemma \ref {7.3} hold. There are six Benlolo-Sanderson invariants and the reader may care to verify that each vanish when the encircled roots are set equal to zero as prescribed by Proposition \ref {4.4}.  }

% For the composition $(2,1,2,1,2,1)$ all five component tableaux are computed.  The resulting labelled lines in $\mathscr T$ as also in \textbf{M}.  From the latter one may verify the swapping lemma and that the criterion for partial linearity of \cite [Section 7]{FJ5} hold.  Thus the five tableaux lead to distinct components of $\mathscr N$.
  \label{fig11}
\end{center}
\end{figure}
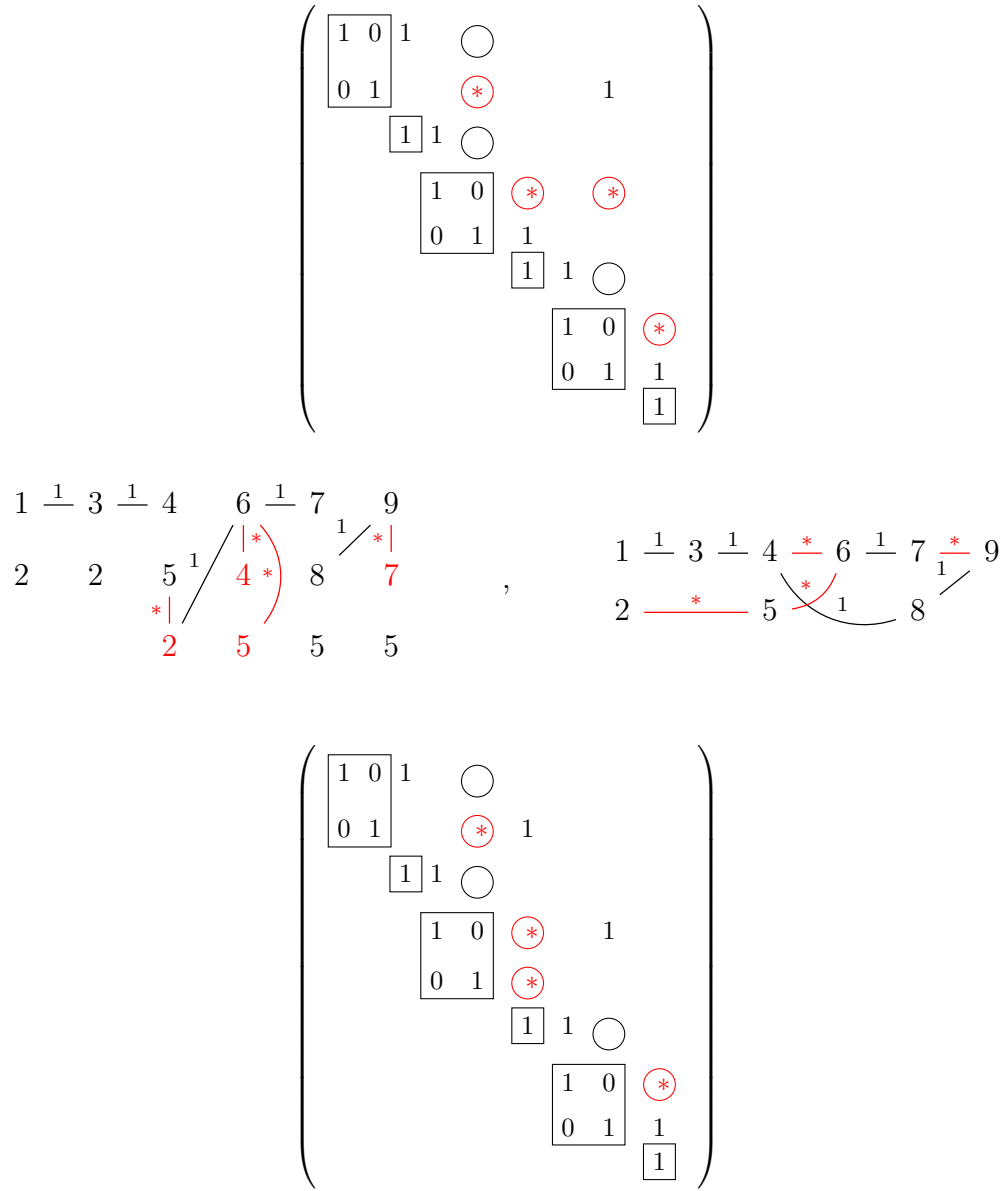

\

\small{\textbf{Example 18}
\begin{figure}[H]
\begin{center}
$\begin{array}{c|c}  \mathcal C & \mathcal C'\\
\begin{tikzcd}[row sep=0.5em,%tiny,
column sep = 1em]
       &{C_1}& {C_2}& {C_3}&{C_4}&{C_5}&{C_6}\\
R_1& 1  \arrow[-,r,"1"]&4  \arrow[-,r,"1"]&6  \arrow[-,r,"1"]&7  \arrow[-,r,"1"]&10&12\\
R_2& 2  \arrow[-,r,"1"]& 5&5&8 \arrow[-,r,"1",red]& 11 \arrow[-,ur,"1"]  &\blue10   \arrow[-,u,"*", blue] \\
R_3& 3& 3& 3 \arrow[-,ur,"1"]& 9& 9  &9\\
R_4&  & & & \red 5  \arrow[-,u,"*", green] \arrow[-,uu,"*", red,  bend right]&  5 &5\\
\end{tikzcd} & \begin{tikzcd}[row sep=0.5em,%tiny,
column sep = 1em]
       &{C_1}& {C_2}& {C_3}&{C_4}&{C_5}&{C_6}\\
R_1& 1  \arrow[-,r,"1"]&4  \arrow[-,r,"1"]&6  \arrow[-,r,"1"]&7  \arrow[-,r,"1"]&10&12  \arrow[-,d,"*", blue]\\
R_2& 2  \arrow[-,r,"1"]& 5&5  \arrow[-,r,"1",red] &8& 11  \arrow[-,ur,"1"] &\blue10 \\
R_3& 3& 3&3& 9  \arrow[-,ur,"1"]& \red8 \arrow[-,u,"*", red]  &8\\
R_4&  & & & \green 3 \arrow[-,u,"*", green]&3 & 3\\
\end{tikzcd}
\\
\begin{tikzpicture}
 \matrix [matrix of math nodes,left delimiter=(,right delimiter=),row sep=0.35em, column sep=0.1em] (m)
 {
1 & 0 &0   &1 &   &   &   &   &   &   &   &   \\
 0& 1  & 0  &   & 1&   &   &   &   &   &   &   \\
0 &  0 &1   &   &   &   &   &  1 &   &   &   &   \\
 &   &   &1   & 0  &1 &   & \cir{ { }}   &\cir{ { }}    &   &   &   \\
 &   &   & 0  &1   &   &   &  \red\cir{ { $\ast$}}   &\cir{ { $\ast$}}    &   &   & \\
 &   &   &   &   &1   & 1  &\cir{ { }}  &\cir{ }     &   &   &   \\
 &   &   &   &   &   &1 & 0  & 0  &1 & \blue\cir{ { }}  &   \\
 &   &   &   &   &   &  0 & 1  &  0 &   &\red{\bold1}   &   \\
 &   &   &   &   &   &0 &0 &1   &   & &   \\
 &   &   &   &   &   &   &   &   &1   & 0  & \blue\cir{ { $\ast$}}   \\
 &   &   &   &   &   &   &   &   &  0 &1   & \blue1  \\
 &   &   &   &   &   &   &   &   &   &   &1 \\
};
 %simple rectangle
%\draw[red,thick ] (m-1-1.north west) rectangle (m-11-11.south east);
 \draw (m-3-1.south west) rectangle (m-1-3.north east);
  \draw(m-9-7.south west) rectangle (m-7-9.north east);
 \draw(m-5-4.south west) rectangle (m-4-5.north east);
  \draw (m-11-10.south west) rectangle (m-10-11.north east);
  \draw (m-6-6.south west) rectangle (m-6-6.north east);
    \draw(m-12-12.south west) rectangle (m-12-12.north east);

\end{tikzpicture}&\begin{tikzpicture}
 \matrix [matrix of math nodes,left delimiter=(,right delimiter=), row sep=0.1em, column sep=0.1em] (m)
 {
1 & 0 &0   &1 &   &   &   &   &\cir{ { }}   &   &   &   \\
 0& 1  & 0  &   &1 &   &   &   & \cir{ { }}  &   &   &   \\
0 &  0 &1   &   &   &   &   &   & \cir{ { $\ast$}}  &   &   &   \\
 &   &   &1   & 0  &1 &   &   &  \cir{ { }} &   &   &   \\
 &   &   & 0  &1   &   &   & \red{\bold1}   & \cir{ { }}  &   &   & \\
 &   &   &   &   &1   & 1  & &   \cir{ { }} &   &   &   \\
 &   &   &   &   &   &1 & 0  & 0  &1 &\cir{ { }}  &   \\
 &   &   &   &   &   &  0 & 1  &  0 &   &   \red\cir{{ $\ast$}} &   \\
 &   &   &   &   &   &0 &0 &1   &   &1 &   \\
 &   &   &   &   &   &   &   &   &1   & 0  & \blue\cir{ { $\ast$}}  \\
 &   &   &   &   &   &   &   &   &  0 &1   &  \blue1 \\
 &   &   &   &   &   &   &   &   &   &   &1 \\
};
 %simple rectangle
%\draw[red,thick ] (m-1-1.north west) rectangle (m-11-11.south east);
 \draw (m-3-1.south west) rectangle (m-1-3.north east);
  \draw(m-9-7.south west) rectangle (m-7-9.north east);
 \draw(m-5-4.south west) rectangle (m-4-5.north east);
  \draw (m-11-10.south west) rectangle (m-10-11.north east);
  \draw (m-6-6.south west) rectangle (m-6-6.north east);
    \draw(m-12-12.south west) rectangle (m-12-12.north east);

\end{tikzpicture}\\
\end{array}$

\caption{The component tableau $\mathscr T^\mathcal C$ (resp. $\mathscr T^{\mathcal C'}$) on the left (resp. right) is given for the composition $(3,2,1,3,2,1)$, with $8$ (resp. $5$) chosen from $\mathscr B^2_{C_2,C_5}$.  This makes $\ell_{8,11}$ the rightmost line as defined in \ref {7.3}. The Benolo-Sanderson is $I:=I^2_{C_2,C_5}$. The encircled root vector $x_{7,11}$ is not required for $I$ to vanish, as a consequence of Lemma \ref {7.3}, verifiable independently in this case. Thus the condition Corollary \ref {7.4} for partial linearity is satisfied. As a consequence the resulting components $\mathscr C$ (resp. $\mathscr C'$) of $\mathscr N$ are distinct, as predicted by Proposition \ref {7.5}.}
%The f The blue $1$, asterisks, and circles within $\mathscr{C}^1$ and $\mathscr{C}^2$ originate from the neighboring columns $C_3$ and $C_6$, each with a height of $1$. One can say the pair $C_3,C_6$ is outside  the neighboring columns $C_1,C_5$ of hight $3$, therefore one can ignore the blue $1$ , $\ast$ and to study the freeness of      rightmost line. one can notice that the red lines  are interchangeable
  \label{fig13}
\end{center}
 \end{figure}
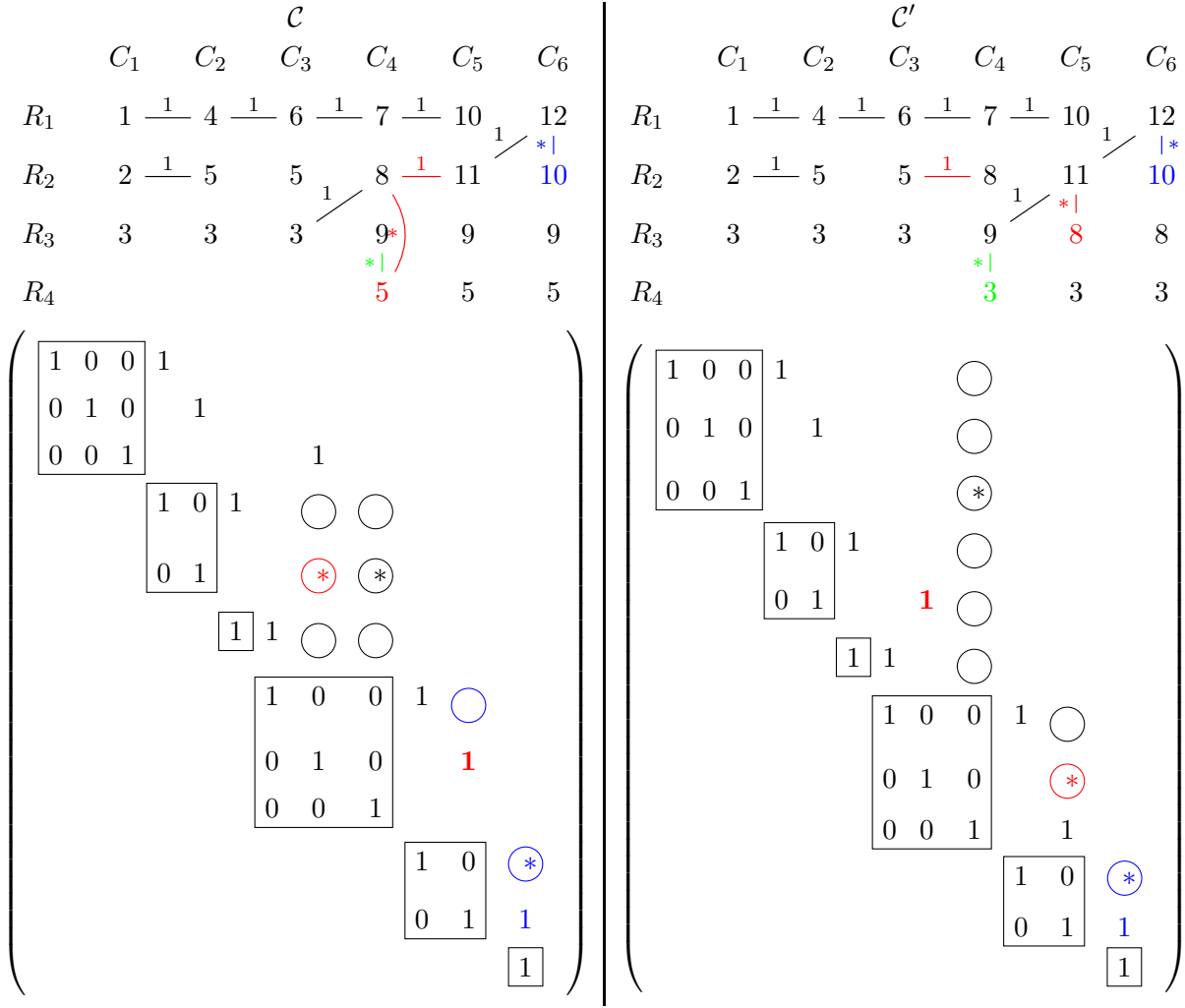

 \textbf{Example 19} %Let the parabolic defined by the sequence $(3,2,1,2,2,1,3)$.

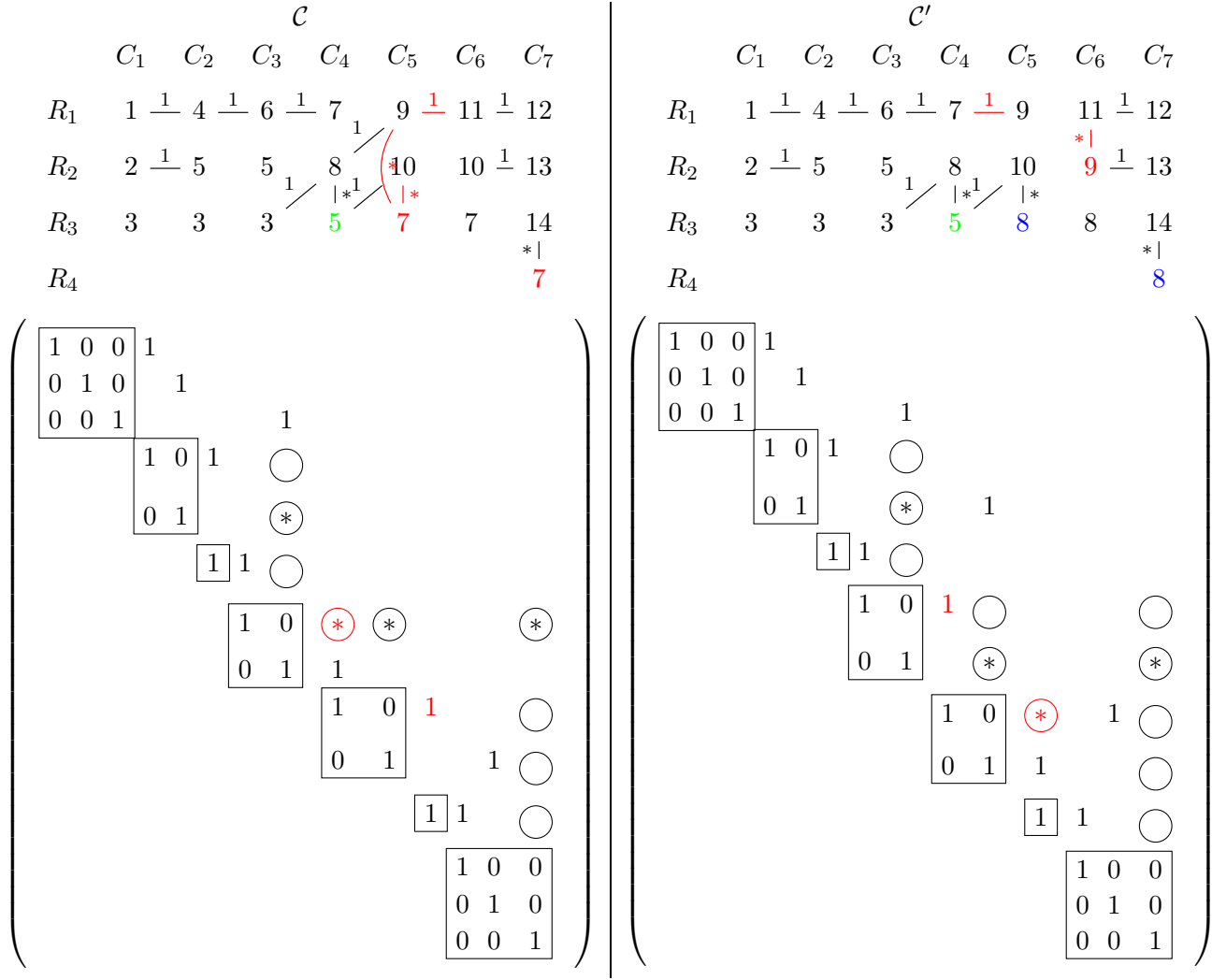
\begin{figure}[H]
\begin{center}
$\left.\begin{array}{c|c}
\mathcal  C  & \mathcal C'\\
  \begin{tikzcd}[row sep=0.5em,%tiny,
column sep = 0.5em]
       &{C_1}& {C_2}& {C_3}&{C_4}&{C_5}&{C_6}&{C_7}\\
R_1& 1  \arrow[-,r,"1"]&4  \arrow[-,r,"1"]&6  \arrow[-,r,"1"]&7  & 9\arrow[-,dd,"*", red,bend right]   \arrow[-,r,"1",red]&11  \arrow[-,r,"1"]&12\\
R_2& 2  \arrow[-,r,"1"]&5  &5  &8  \arrow[-,d,"*"] \arrow[-,ru,"1"]& 10  \arrow[-,d,"*",red]  & 10   \arrow[-,r,"1"]&13\\
R_3& 3&3  & 3  \arrow[-,ru,"1"] &\green5 \arrow[-,ru,"1"]   & \red 7   &  7  &14\\
R_4& &  &  & && & \red 7\arrow[-,u,"*"] \\
\end{tikzcd} &   \begin{tikzcd}[row sep=0.5em,%tiny,
column sep = 0.5em]
       &{C_1}& {C_2}& {C_3}&{C_4}&{C_5}&{C_6}&{C_7}\\
R_1& 1  \arrow[-,r,"1"]&4  \arrow[-,r,"1"]&6  \arrow[-,r,"1"]& 7 \arrow[-,r,"1",red]  &9  &11  \arrow[-,r,"1"]&12\\
R_2& 2  \arrow[-,r,"1"]&5  &5  & 8  \arrow[-,d,"*"]& 10  \arrow[-,d,"*"]  & \red 9  \arrow[-,u,"*",red]  \arrow[-,r,"1"]&13\\
R_3& 3&3  & 3  \arrow[-,ru,"1"] &\green5 \arrow[-,ru,"1"]   &\blue8  & 8  &14\\
R_4& &  &  & && & \blue8\arrow[-,u,"*"] \\
\end{tikzcd}\\ \small
\begin{tikzpicture}
 \matrix [matrix of math nodes,left delimiter=(,right delimiter=),row sep=0.05em, column sep=0.em] (m)
 {
1 & 0  &0   &1   &   &   &   &   &   &   &   &  &  & \\
0 & 1  &0   &   &  1 &   &   &   &   &   &   &  &  & \\
0 & 0  &1   &   &   &   &   & 1  &   &   &   &  &  & \\
 &   &   & 1  & 0  &1   &   & \cir{ }  &   &   &   &  &  & \\
  &   &   & 0  & 1  &   &   &\cir{ $\ast$}   &   &   &   &  &  & \\
  &   &   &   &   &  1 & 1  &\cir { }   &   &   &   &  &  & \\
  &   &   &   &   &   & 1  &0  & \red \cir{ $\ast$}  &\cir{ $\ast$}   &   &  &  &\cir{ $\ast$} \\
    &   &   &   &   &   & 0  &1  & 1 &   &   &  &  & \\
    &   &   &   &   &   &   &  & 1 &  0 &\red 1   &  &  &\cir { } \\
    &   &   &   &   &   &   &  & 0 &  1 &   &  &1  &\cir { } \\
&   &   &   &   &   &   &  & &   & 1  & 1 &  &\cir { } \\
&   &   &   &   &   &   &  & &   &   & 1 & 0 &0 \\
&   &   &   &   &   &   &  & &   &   & 0 & 1 &0 \\
&   &   &   &   &   &   &  & &   &   & 0 & 0 &1 \\
};
 %simple rectangle

 \draw (m-3-1.south west) rectangle (m-1-3.north east);
  \draw(m-14-12.south west) rectangle (m-12-14.north east);
 \draw(m-5-4.south west) rectangle (m-4-5.north east);
  \draw (m-8-7.south west) rectangle (m-7-8.north east);
  \draw (m-6-6.south west) rectangle (m-6-6.north east);
  \draw(m-11-11.south west) rectangle (m-11-11.north east);
  \draw(m-10-9.south west) rectangle (m-9-10.north east);

\end{tikzpicture}& \small \begin{tikzpicture}
 \matrix [matrix of math nodes,left delimiter=(,right delimiter=),row sep=0.em, column sep=0.em] (m)
 {
1 & 0  &0   &1   &   &   &   &   &   &   &   &  &  & \\
0 & 1  &0   &   &  1 &   &   &   &   &   &   &  &  & \\
0 & 0  &1   &   &   &   &   & 1  &   &   &   &  &  & \\
 &   &   & 1  & 0  &1   &   & \cir{ }  &   &   &   &  &  & \\
  &   &   & 0  & 1  &   &   &\cir{ $\ast$}   &   & 1  &   &  &  & \\
  &   &   &   &   &  1 & 1  &\cir { }   &   &   &   &  &  & \\
  &   &   &   &   &   & 1  &0  &\red 1  & \cir{ }  &   &  &  &\cir{ } \\
    &   &   &   &   &   & 0  &1  &  & \cir { $\ast$}  &   &  &  &\cir { $\ast$} \\
    &   &   &   &   &   &   &  & 1 &  0 &\red \cir { $\ast$}  &  & 1 &\cir { } \\
    &   &   &   &   &   &   &  & 0 &  1 & 1  &  &  &\cir { } \\
&   &   &   &   &   &   &  & &   & 1  & 1 &  &\cir { } \\
&   &   &   &   &   &   &  & &   &   & 1 & 0 &0 \\
&   &   &   &   &   &   &  & &   &   & 0 & 1 &0 \\
&   &   &   &   &   &   &  & &   &   & 0 & 0 &1 \\
};
 %simple rectangle

 \draw (m-3-1.south west) rectangle (m-1-3.north east);
  \draw(m-14-12.south west) rectangle (m-12-14.north east);
 \draw(m-5-4.south west) rectangle (m-4-5.north east);
  \draw (m-8-7.south west) rectangle (m-7-8.north east);
  \draw (m-6-6.south west) rectangle (m-6-6.north east);
  \draw(m-11-11.south west) rectangle (m-11-11.north east);
  \draw(m-10-9.south west) rectangle (m-9-10.north east);

\end{tikzpicture}
\end{array}\right.$
\caption{The component tableau $\mathscr T^\mathcal C$ (resp. $\mathscr T^{\mathcal C'}$) on the left (resp. right) is given for the composition $(3,2,1,2,2,1,3)$, with $7$ (resp. $8$) chosen from $\mathscr B^2_{C_2,C_5}$.  This makes $\ell_{9,11}$ the rightmost line as defined in \ref {7.3}. The Benolo-Sanderson is $I:=I^1_{C_3,C_6}$. The encircled root vectors in the last column, notably $x_{7,14},x_{9,14}$,  are not needed for $I$ to vanish, as a consequence of Lemma \ref {7.3} verifiable independently in this case. Thus the condition Corollary \ref {7.4} for partial linearity is satisfied. As a consequence the resulting components $\mathscr C$ (resp. $\mathscr C'$) of $\mathscr N$ are distinct.} \label{fig14}
\end{center}
 \end{figure}

 % \caption{The component tableau $\mathscr T^\mathcal C$ (resp. $\mathscr T^{\mathcal C'}$) on the left (resp. right) is given for the composition $(3,2,1,2,2,1,3)$, with $7$ (resp. $8$) chosen from $\mathscr B^2_{C_2,C_5}$.  This makes $\ell_{9,11}$ the rightmost line \ref {7.3}. The Benolo-Sanderson $E:=I^2_{C_2,C_5}$. The encircled root vector $x_{9,14}$ is not required for $I$ to vanish, as a consequence of Lemma \ref {7.3} verifiable independently in this case. Thus the condition Corollary \ref {7.4} for partial linearity is satisfied. As a consequence the resulting components $\mathscr C$ (resp. $\mathscr C'$ of $\mathscr N$ are distinct.}

 \textbf{Example 20.}
 \begin{figure}[H]
\begin{center}
$\begin{array}{c|c}  \mathcal  C  & \mathcal C'\\
  \begin{tikzcd}[row sep=0.5em,%tiny,
column sep = 0.5em]
       &{C_1}& {C_2}& {C_3}&{C_4}&{C_5}&{C_6}\\
R_1& 1  \arrow[-,r,"1"]&4  \arrow[-,r,"1"]&6  \arrow[-,r,"1"]& 7  &10 \arrow[-,d,"*"] \arrow[-,r,"1"]&11 \\
R_2& 2  \arrow[-,r,"1"]&5  & 5  &8\arrow[-,ur,"1",red]  &\green7 \arrow[-,r,"1"] &12 \\
R_3& 3 &3  &3 \arrow[-,r,"1"] &9\arrow[-,d,"*"]   &9 &9 \\
R_4&   &   &  &\red5\arrow[-,uu,"*",red,bend right]   &5 &5 \\
\end{tikzcd} &   \begin{tikzcd}[row sep=0.5em,%tiny,
column sep = 0.5em]
       &{C_1}& {C_2}& {C_3}&{C_4}&{C_5}&{C_6}\\
R_1& 1  \arrow[-,r,"1"]&4  \arrow[-,r,"1"]&6  \arrow[-,r,"1"]&7  &10 \arrow[-,d,"*"] \arrow[-,r,"1"]&11 \\
R_2& 2  \arrow[-,r,"1"]&5  &5  \arrow[-,r,"1",red] &8  &\green7 \arrow[-,r,"1"] &12 \\
R_3& 3 &3  & 3  &9\arrow[-,d,"*"]  \arrow[-,ruu,"1",bend left=12 ] &\red8 \arrow[-,uu,"*",red,bend left] &8 \\
R_4&   &   &  &\blue3   &3 &3 \\
\end{tikzcd} \\
\begin{tikzpicture}
 \matrix [matrix of math nodes,left delimiter=(,right delimiter=),row sep=0.2em, column sep=0.1em] (m)
 {
1 & 0 &0   &1 &   &   &   &   &   &   &   &   \\
 0& 1  & 0  &   & 1&   &   &   &   &   &   &   \\
0 &  0 &1   &   &   &   &   &   & 1  &   &   &   \\
 &   &   &1   & 0  &1 &   & \cir{ { }}   &\cir{ { }}    &   &   &   \\
 &   &   & 0  &1   &   &   &  \red\cir{ { $\ast$}}   &\cir{ { $\ast$}}    &   &   & \\
 &   &   &   &   &1   & 1  &\cir{ { }}  &\cir{ { }}     &   &   &   \\
 &   &   &   &   &   &1 & 0  & 0  & \cir{ {$\ast$ }}  &  &1   \\
 &   &   &   &   &   &  0 & 1  &  0 &   \red1 &   &   \\
 &   &   &   &   &   &0 &0 &1   &   & &   \\
 &   &   &   &   &   &   &   &   &1   &  1 &    \\
 &   &   &   &   &   &   &   &   &   &1   & 0  \\
 &   &   &   &   &   &   &   &   &   & 0  &1 \\
};
 %simple rectangle

 \draw (m-3-1.south west) rectangle (m-1-3.north east);
  \draw(m-9-7.south west) rectangle (m-7-9.north east);
 \draw(m-5-4.south west) rectangle (m-4-5.north east);
  \draw (m-12-11.south west) rectangle (m-11-12.north east);
  \draw (m-6-6.south west) rectangle (m-6-6.north east);
    \draw(m-10-10.south west) rectangle (m-10-10.north east);

\end{tikzpicture}& \begin{tikzpicture}
 \matrix [matrix of math nodes,left delimiter=(,right delimiter=),row sep=0.em, column sep=0.1em] (m)
 {
1 & 0 &0   &1 &   &   &   &   &   \cir{ { }} &   &   &   \\
 0& 1  & 0  &   & 1&   &   &   &  \cir{ { }}  &   &   &   \\
0 &  0 &1   &   &   &   &   &   &\cir{ { $\ast$}}  &   &   & \\
 &   &   &1   & 0  &1 &   &   &\cir{ { }}  &   &   &  \\
 &   &   & 0  &1   &   &   &   \red1  &   \cir{ { }}  &   &   & \\
 &   &   &   &   &1   & 1  & &    \cir{ { }}  &   &   &   \\
 &   &   &   &   &   &1 & 0  &   0  &   \cir{ { $\ast$}}  &   &1 \\
 &   &   &   &   &   &  0 & 1  &  0 &   \red\cir{ { $\ast$}}  &   &  \\
 &   &   &   &   &   &0 &0 &1   &  1 & &   \\
 &   &   &   &   &   &   &   &   &1   &  1 &    \\
 &   &   &   &   &   &   &   &   &   &1   & 0  \\
 &   &   &   &   &   &   &   &   &   & 0  &1 \\
};
 %simple rectangle

 \draw (m-3-1.south west) rectangle (m-1-3.north east);
  \draw(m-9-7.south west) rectangle (m-7-9.north east);
 \draw(m-5-4.south west) rectangle (m-4-5.north east);
  \draw (m-12-11.south west) rectangle (m-11-12.north east);
  \draw (m-6-6.south west) rectangle (m-6-6.north east);
    \draw(m-10-10.south west) rectangle (m-10-10.north east);

\end{tikzpicture}
\end{array}$
\caption{ The component tableau $\mathscr T^\mathcal C$ (resp. $\mathscr T^{\mathcal C'}$) on the left (resp. right) is given for the composition $(3,2,1,3,1,2)$, with $5$ (resp. $8$) chosen from $\mathscr B^2_{C_2,C_6}$.  This makes $\ell_{8,10}$ the rightmost line as defined in \ref {7.3}. The Benolo-Sanderson is $I:=I^2_{C_2,C_6}$. The encircled root vector $x_{7,10}$ is not required for $I$ to vanish, as a consequence of Lemma \ref {7.3} verifiable independently in this case. Thus the condition Corollary \ref {7.4} for partial linearity is satisfied. As a consequence the resulting components $\mathscr C$ (resp. $\mathscr C'$) of $\mathscr N$ are distinct.
 } \label{fig15}
\end{center}
 \end{figure}
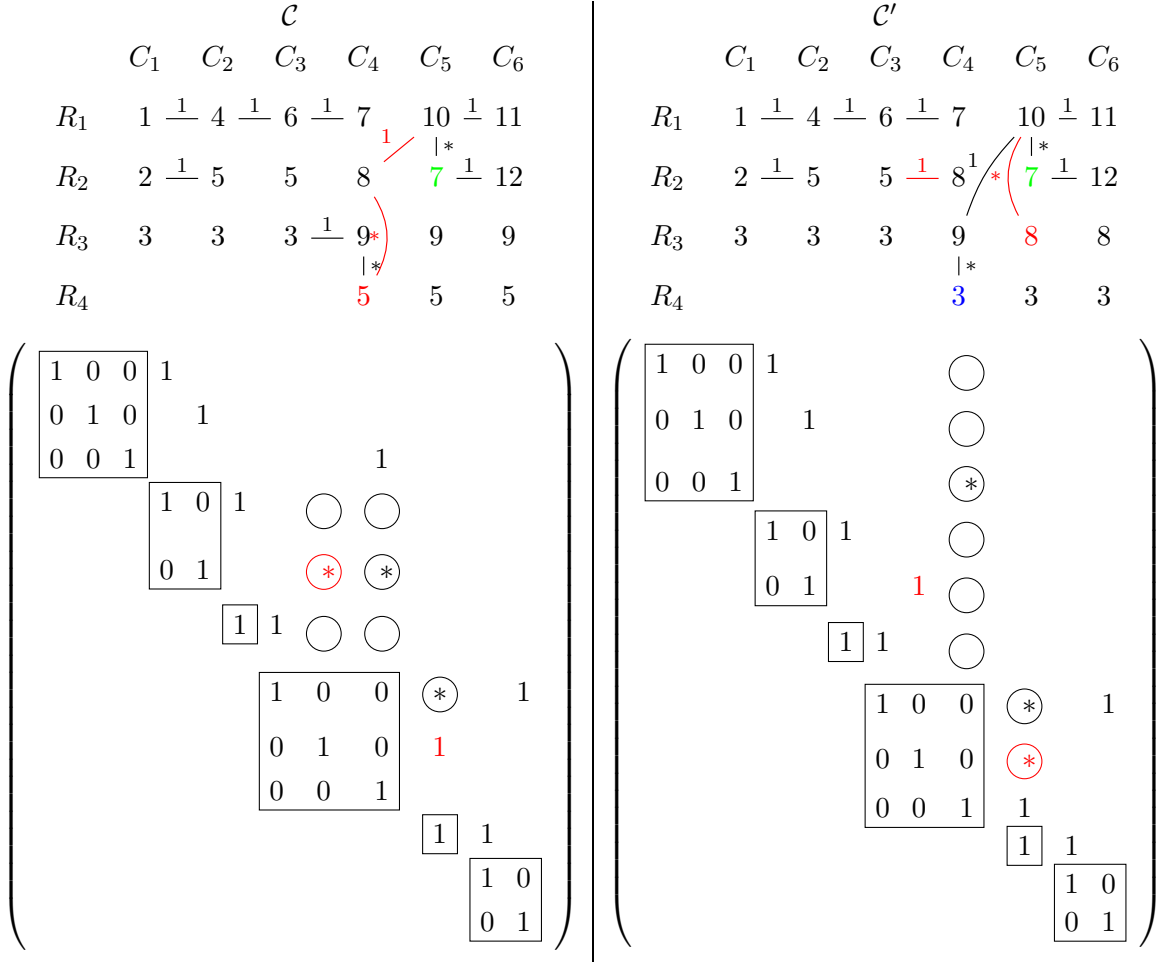

\end{document}